\documentclass[10pt]{article}

\RequirePackage{amsmath}
\RequirePackage{amsfonts}
\RequirePackage{amssymb}
\RequirePackage{enumitem}
\RequirePackage{hyperref}
\RequirePackage{multicol}
\RequirePackage{multirow}
\RequirePackage{algorithm}
\RequirePackage{algpseudocode}
\RequirePackage{graphicx}
\RequirePackage{xcolor}
\RequirePackage{amsthm}
\RequirePackage{tikz}
\RequirePackage{booktabs, multirow}
\RequirePackage{changepage,threeparttable}
\allowdisplaybreaks

\usepackage{dpr_fec,ifthen,dsfont}

\newtheorem{assumption}{Assumption}[section]

\newtheorem{parameterrule}{Parameter Rule}[section]

\newcommand{\edit}[1]{\textcolor{black}{#1}}
\usepackage{isetechreport}
\def\reportyear{23T}
\def\reportno{018}
\def\revisionno{0}
\def\originaldate{\today}
\coraltrue
\cvcrfalse

\begin{document}

\title{A Stochastic-Gradient-based Interior-Point Algorithm for Solving Smooth Bound-Constrained Optimization Problems}

\author{Frank E.~Curtis\thanks{E-mail: frank.e.curtis@lehigh.edu}}
\affil{Department of Industrial and Systems Engineering, Lehigh University}
\author{Vyacheslav Kungurtsev\thanks{E-mail: kunguvya@fel.cvut.cz}}
\affil{Department of Computer Science, Czech Technical University}
\author{Daniel P.~Robinson\thanks{E-mail: daniel.p.robinson@lehigh.edu}}
\affil{Department of Industrial and Systems Engineering, Lehigh University}
\author{Qi Wang\thanks{E-mail: qiw420@lehigh.edu}}
\affil{Department of Industrial and Systems Engineering, Lehigh University}
\titlepage

\maketitle

\begin{abstract}
  A stochastic-gradient-based interior-point algorithm for minimizing a continuously differentiable objective function (that may be nonconvex) subject to bound constraints is presented, analyzed, and demonstrated through experimental results.  The algorithm is unique from other interior-point methods for solving smooth \edit{nonconvex} optimization problems since the search directions are computed using stochastic gradient estimates.  It is also unique in its use of inner neighborhoods of the feasible region---defined by a positive and vanishing neighborhood-parameter sequence---in which the iterates are forced to remain.  It is shown that with a careful balance between the barrier, step-size, and neighborhood sequences, the proposed algorithm satisfies convergence guarantees in both deterministic and stochastic settings.  The results of numerical experiments show that in both settings the algorithm can outperform \edit{projection-based} methods.
\end{abstract}

\newcommand{\barrier}{\tilde\phi}
\newcommand{\buff}{{\rm buff}}
\newcommand{\hood}{\Ncal_{[l,u]}}
\newcommand{\pre}{{\rm pre}}
\renewcommand{\diag}{\Psi}

\section{Introduction}\label{sec.introduction}

The interior-point methodology is one of the most effective approaches for solving continuous constrained optimization problems.  In the context of \edit{deterministic} derivative-based algorithmic strategies, interior-point methods offer convergence guarantees from remote starting points \cite{ByrdHribNoce99,NestNemi94,WaecBieg06}, and in both convex and nonconvex settings such algorithms can offer good worst-case iteration complexity properties \cite{BianChenYe15,NestNemi94}.  Furthermore, many of the most popular software packages for solving large-scale continuous optimization problems are based on interior-point methods \cite{mosek,ByrdHribNoce99,Stur99,TohToddTutu99,VandShan99,WaecBieg06}, and these have been used to great effect for many years.

Despite the extensive literature on theoretical and practical benefits of interior-point methods in the context of \edit{deterministic} derivative-based algorithms for solving \edit{non}convex optimization problems, to the best of our knowledge there has not yet been one that has been shown rigorously to offer convergence guarantees when neither \edit{exact} function nor \edit{exact} derivative evaluations are available, and instead only stochastic gradient estimates are employed.  (An interior-point stochastic-approximation method was proposed and tested in~\cite{NIPS2008_a87ff679}, but as we mention in Remark~\ref{rem.carbonetto} on page~\pageref{rem.carbonetto}, the claimed asymptotic-convergence guarantee in \cite{NIPS2008_a87ff679} overlooks a critical issue related to the step sizes.)  In this paper, we propose, analyze, and present the results of experiments with such an algorithm.  Randomized algorithms for minimizing a linear function over a convex set have been proposed \cite{BadedeKl22,Nara16}, but the setting and the techniques those algorithms employ are distinct from the ones considered in this paper.

For a straightforward presentation of our proposed strategy and its convergence guarantees, we focus on the case of constrained optimization with bound constraints only.  That said, our algorithmic strategies have been designed so that they may be extended for solving problems with continuous \edit{and potentially nonlinear} equality and/or inequality constraints as well.  For example, since interior-point methods handle inequality constraints through the introduction of an additional objective function term that is weighted by a barrier parameter and a continuation strategy that reduces the barrier parameter iteratively, one might consider extending our algorithmic ideas using the recently proposed stochastic algorithms for solving equality-constrained optimization presented in \cite{BeraBollZhou22,BeraCurtOneiRobi21,BeraCurtRobiZhou21,BeraShiYiZhou22,CurtOneiRobi21,CurtRobiZhou21,FangNaMahoKola22,NaAnitKola22,NaMaho22,qiu2023sequential}.  The main challenge to address in such potential extensions is the one that we address in this paper, namely, that derivatives of the barrier function are not Lipschitz continuous in the interior of a set of bound constraints.  We focus primarily on the setting of minimizing an objective that may be nonconvex.  Upon seeing our algorithm, a reader may wonder if a simpler variant has convergence guarantees.  However, we discuss in Section~\ref{sec.convex} why the challenges that we overcome in the general\edit{, potentially nonconvex} setting are not readily avoided with a simpler variant, even in the strongly convex setting.

\subsection{Contributions}\label{sec.contributions}

We propose, analyze, and provide the results of numerical experiments with a \edit{(stochastic-)}gradient-based interior-point method for solving \edit{potentially nonconvex} bound-constrained optimization problems.  \edit{Our contributions are motivated primarily by the stochastic setting, but since our algorithm and analysis are unique for the deterministic setting---and since our experiments show that our algorithm can be competitive with a state-of-the-art projected gradient method in practice---we provide a complete investigation of this setting as well.}

Our algorithm involves multiple \edit{features that are unique from those involved in other deterministic- and stochastic-gradient-based} interior-point methods that have been proposed and analyzed in the literature.  \edit{Our main contributions are as follows.}

\benumerate[leftmargin=20pt]
  \item[(i)] Our algorithm employs a prescribed monotonically nonincreasing and vanishing barrier parameter sequence.  In this manner, the algorithm does not rely on the ability to compute values for derivative-based stationarity tests, as is done in derivative-based interior-point methods for deciding whether to decrease the barrier parameter in a given iteration.  This is significant since stationarity measures cannot be computed accurately in the stochastic setting that we consider.
  \item[(ii)] Our algorithm does not employ a fraction-to-the-boundary rule.  Such a rule is critical for convergence guarantees of other derivative-based interior-point methods for solving nonconvex problems (see, e.g., \cite{ByrdHribNoce99,WaecBieg06}), since it ensures that---with respect to a threshold that depends on the barrier parameter value---the iterates do not get too close to the boundary of the feasible region.  This, in turn, ensures that the iterates remain in a region in which derivatives of the barrier function are Lipschitz continuous.  By contrast, in our proposed algorithm, we employ a unique strategy that involves keeping each iterate within an inner neighborhood of the feasible region.  \edit{Within our algorithm}, these neighborhoods are defined by a monotonically nonincreasing and vanishing sequence.
  \item[(iii)] Our algorithm does not rely on step acceptance criteria (e.g., using line searches, trust regions, regularization, etc.)~that in turn rely on exact objective function evaluations.  This is significant since such evaluations are intractable in various settings of interest; see, e.g., \cite{BottCurtNoce18}.  That said, as is common for other stochastic optimization algorithms, our convergence guarantees rely on knowledge of problem-dependent quantities, including a Lipschitz constant for the gradient of the objective.  (In practice, the problem-dependent quantities that our algorithm requires can be estimated using function and/or derivative values.)
  \item[(iv)] We present general sets of conditions under which the algorithm's barrier, neighborhood, and step size sequences are balanced so as to ensure convergence guarantees in both deterministic and stochastic settings.
  \item[(v)] We show by a representative comparison that, in deterministic and stochastic settings, \edit{our} algorithm can outperform a \edit{projection-based method}.
\eenumerate
One aspect that limits the applicability of our work is that, for the stochastic setting, we assume that the errors of the stochastic gradient estimators are bounded by a known value.  Convergence guarantees have been established under looser assumptions in unconstrained settings \cite{BottCurtNoce18}, but since unique challenges arise in the constrained setting, we consider ours a significant first step in the design and analysis of interior-point algorithms for constrained stochastic optimization.

\subsection{Notation}\label{sec.notation}

We use $\R{}$ to denote the set of real numbers, $\Rext{}$ to denote the set of extended-real numbers (i.e., $\Rext{} := \R{} \cup \{-\infty,\infty\}$), and $\R{}_{\geq a}$ (resp.,~$\R{}_{>a}$, $\R{}_{<a}$, or $\R{}_{\leq a}$) to denote the set of real numbers greater than or equal to (resp.,~greater than, less than, or less than or equal to) $a \in \R{}$.  We append a superscript to such a set to denote the space of vectors or matrices whose elements are restricted to the indicated set; e.g.,~we use $\R{n}$ to denote the set of $n$-dimensional real vectors and~$\R{m \times n}$ to denote the set of $m$-by-$n$-dimensional real matrices.  We use $\Smbb^n \subset \R{n \times n}$ to denote the set of $n$-by-$n$-dimensional real symmetric matrices.  We use $\N{} := \{1,2,\dots\}$ to denote the set of positive integers and, given $n \in \N{}$, we denote~$[n] := \{1,\dots,n\}$.

Given $\Bcal \subseteq \R{n}$, we use $\interior(\Bcal)$ to denote the interior of $\Bcal$.  We use $\ones$ to denote a vector of ones whose length is determined by the context in which it appears.  Given $(A,B) \in \Smbb^n \times \Smbb^n$, we write $A \succeq B$ (resp.,~$A \succ B$) to indicate that $A - B \in \Smbb^n$ is positive semidefinite (resp.,~positive definite).  Given $l \in \Rext{n}$, we use $\diag(l)$ to denote the extended-real-valued diagonal matrix whose $(i,i)$-element is equal to $l_i$ for all $i \in [n]$\edit{, where $l_i$ is the $i$th element of $l$}.  Given a sequence of real-valued vectors $\{\mu_k\}$ and $\Mcal \subseteq \R{n}$\edit{, where $\mu_k \in \R{n}$ for all $k \in \N{}$ with the subscript $k$ indicating the index number in the sequence}, we write $\{\mu_k\} \subset \Mcal$ to indicate that $\mu_k \in \Mcal$ for all $k \in \N{}$.  Moreover, for a real-number sequence, we write $\{\mu_k\} \searrow 0$ to indicate that (a) $\{\mu_k\} \subset \R{}_{>0}$, (b) $\{\mu_k\}$ is monotonically nonincreasing, and (c) the limit of $\{\mu_k\}$ is zero.  Given sequences $\{a_k\} \subset \R{}_{\geq0}$ and $\{b_k\} \subset \R{}_{>0}$, we write $\{a_k\} = \Ocal(b_k)$ (resp., \edit{$\{a_k\} = \Omega(b_k)$}) to indicate that there exists $C \in \R{}_{>0}$ (resp., \edit{$c \in \R{}_{>0}$}) such that $a_k \leq C b_k$ (resp., \edit{$a_k \geq c b_k$}) for all sufficiently large $k \in \N{}$.  Given such sequences, we write $\{a_k\} = o(b_k)$ to indicate that $\{a_k/b_k\} \to 0$.  

The algorithm that we propose is iterative in the sense that any given run of the algorithm produces an iterate sequence (of real-valued vectors) $\{x_k\} \subset \R{n}$.  Like for the iterate sequence, we append a positive integer as a subscript for a quantity to denote its value during an iteration of an algorithm.  \edit{Also, as in the previous paragraph, we sometimes use a subscript to indicate an element index of a vector or matrix.  In all such cases, the meaning of a subscript is clear from the context.}  Multiple subscripts are used in some cases, as needed; e.g., the $i$th element of the $k$th iterate $x_k \in \R{n}$ is denoted as $x_{k,i} \in \R{}$ and the $(i,i)$-element of a matrix $H_k \in \Smbb^n$ is denoted as $H_{k,i,i}$.

At times, we express algebraic operations using quantities with infinite magnitude, namely, $-\infty$ and $\infty$.  In such cases, we adopt natural conventions.  In particular, given $a \in \R{}$, we let $\infty - a = \infty$ and $a - (-\infty) = \infty$, and, given $a \in \R{}_{>0}$, we let $a \cdot \infty = \infty$ and $a / \infty = 0$.  Given a pair of nonnegative extended-real-number-valued vectors $(a,b) \in \Rext{n}_{\geq0} \times \Rext{n}_{\geq0}$, we write $a \perp b$ to indicate that $a_i = 0$ and/or $b_i = 0$ for all $i \in [n]$.

\subsection{Organization}\label{sec.organization}

Our main problem of interest, namely, minimizing a potentially nonconvex continuous function over a set of bound constraints, is stated formally along with a presentation of our main algorithm in Section~\ref{sec.algorithm}.  Our convergence analyses for this algorithm are presented in Section~\ref{sec.analysis}.  In Section~\ref{sec.convex}, we discuss the obstacles of proving a convergence guarantee for a simpler variant of our algorithm, even in the strongly convex setting.  The results of numerical experiments are presented in Section~\ref{sec.numerical} and concluding remarks are provided in Section~\ref{sec.conclusion}.

\section{Algorithm}\label{sec.algorithm}

Our main problem of interest is to minimize an objective function over a feasible region that we denote as $\Bcal := [l,u] \equiv \{x \in \R{n} : l \leq x \leq u\}$, where $(l,u) \in \Rext{n} \times \Rext{n}$ with $l_i < u_i$ for all $i \in [n]$.  We assume that at least one element of $(l,u)$ is finite, so the problem is indeed constrained.  For the sake of generality, we only require that the objective has as its domain an open set $\Bcal^+$ containing $\Bcal$.  Denoting this objective function as $f : \Bcal^+ \to \R{}$, we write our problem of interest as
\bequation\label{prob.opt}
  \min_{x \in \R{n}}\ f(x)\ \ \text{subject to}\ \ x \in \Bcal := [l,u].
\eequation
We make the following assumption pertaining to this \edit{potentially nonconvex} $f$.

\bassumption\label{ass.f}
  The objective function $f : \Bcal^+ \to \R{}$ is continuously differentiable over $\Bcal^+$, bounded below by $f_{\inf} \in \R{}$ over $\Bcal$, and bounded above by $f_{\sup} \in \R{}$ over $\Bcal$.  In addition, its gradient function $\nabla f : \Bcal^+ \to \R{n}$ is Lipschitz continuous with respect to the 2-norm over $\Bcal$ with constant $\ell_{\nabla f,\Bcal} \in \R{}_{>0}$ and is bounded in 2-norm $($resp.,~$\infty$-norm$)$ over $\Bcal$ by $\kappa_{\nabla f,\Bcal, 2} \in \R{}_{>0}$ $($resp.,~$\kappa_{\nabla f, \Bcal, \infty} \in \R{}_{>0}$$)$.
\eassumption

\noindent
Assumption~\ref{ass.f} is mostly standard.  The nonstandard aspect is the assumption that~$f$ is bounded above over $\Bcal$; this is a relatively loose assumption to handle extreme cases in the stochastic setting.  The existence of $\kappa_{\nabla f,\Bcal, \infty}$ follows from that of $\kappa_{\nabla f,\Bcal, 2}$, and vice versa, but we define both for the sake of notational convenience.

Under Assumption~\ref{ass.f}, specifically under the assumption that $f$ is continuously differentiable over an open set containing the feasible region $\Bcal$, it follows that if $x \in \R{n}$ is a minimizer of \eqref{prob.opt}, then there must exist $(y,z) \in \R{n} \times \R{n}$ such that $(x,y,z)$ satisfies the Karush-Kuhn-Tucker (KKT) conditions given by
\bequation\label{eq.KKT}
  \nabla f(x) - y + z = 0,\ \ 0 \leq (x - l) \perp y \geq 0,\ \ \text{and}\ \ 0 \leq (u - x) \perp z \geq 0.
\eequation
Defining the index sets of finite bounds as $\Lcal := \{i \in [n] : l_i > -\infty\}$ and $\Ucal := \{i \in [n] : u_i < \infty\}$, $x \in \R{n}$ implies that $x_i - l_i = \infty > 0$ for all $i \in [n] \setminus \Lcal$ and $u_i - x_i = \infty > 0$ for all $i \in [n] \setminus \Ucal$, meaning that \eqref{eq.KKT} and our definition of the operator $\perp$ in Section~\ref{sec.notation} require that $y_i = 0$ for all $i \in [n] \setminus \Lcal$ and $z_i = 0$ for all $i \in [n] \setminus \Ucal$.

Central ideas of the interior-point methodology are to replace inequality constraints with a parameterized barrier function in the objective, and to solve the original constrained optimization problem through a continuation approach by driving the barrier parameter to zero.  For example, using a so-called log-barrier in the context of \eqref{prob.opt}, this amounts to introducing the barrier parameter $\mu \in \R{}_{>0}$ and using the log-barrier-augmented function $\phi : \interior(\Bcal) \times \R{}_{>0} \to \R{}$ given by
\bequation\label{def.barrier}
  \phi(x,\mu) = f(x) - \mu \sum_{i \in \Lcal} \log(x_i - l_i) - \mu \sum_{i \in \Ucal} \log(u_i - x_i).
\eequation
(We use a log-barrier function throughout the paper, although one might extend our algorithm and analysis for other barrier functions as well.)  Given $\mu \in \R{}_{>0}$ and letting~$\nabla_x \phi$ denote the gradient operator of $\phi$ with respect to its first argument, a minimizer of the barrier-augmented function $\phi(\cdot,\mu)$ over $\interior(\Bcal)$ must satisfy
\bequation\label{eq.barrier_stationary}
  0 = \nabla_x \phi(x,\mu) \equiv \nabla f(x) - \mu \diag(x - l)^{-1} \ones + \mu \diag(u - x)^{-1} \ones.
\eequation

A traditional interior-point method for derivative-based nonconvex optimization would involve fixing the barrier parameter at a value $\mu \in \R{}_{>0}$, employing an unconstrained optimization method to minimize $\phi(\cdot,\mu)$ until an approximate stationarity tolerance is satisfied (with a safeguard such as a fraction-to-the-boundary rule to ensure that the iterates remain within $\interior(\Bcal)$), then reducing the barrier parameter and repeating the procedure in an iterative manner.  However, for our purposes, we avoid the need to check a stationarity tolerance, since this would require an evaluation of a gradient of the objective (see~\eqref{eq.barrier_stationary}), which we presume is intractable to obtain.

Our proposed algorithm, by contrast, employs a prescribed positive barrier parameter sequence $\{\mu_k\} \searrow 0$ and a line-search-free strategy for generating the positive step size sequence $\{\alpha_k\} \subset \R{}_{>0}$.  We state the algorithm in a generic manner, but in our analyses in Section~\ref{sec.analysis} we reveal specific requirements that these sequences must satisfy to yield convergence guarantees in deterministic and stochastic settings.  In addition, rather than rely on a safeguard such as a fraction-to-the-boundary rule---which presents challenges in terms of ensuring convergence guarantees in a stochastic setting since such a rule would enforce an iterate-dependent bound on the steps---our algorithm employs a rule that ensures that, for all $k \in \N{}$, the subsequent iterate remains sufficiently within $\interior(\Bcal)$ by a prescribed margin.  For this, we introduce
\bequation\label{eq.hood}
  \hood(\theta) := \{x \in \R{n} : l + \theta \leq x \leq u - \theta\},
\eequation
and, for all $k \in \N{}$, have the algorithm ensure $x_{k+1} \in \hood(\theta_k)$ for some $\theta_k \in \R{}_{>0}$. The positive sequence $\{\theta_k\} \searrow 0$, base value $\theta_0 \in \R{}_{>\theta_1}$, and initial point $x_1 \in \R{n}$ must be prescribed for each run of the algorithm, and the latter two must satisfy
\bequation\label{eq.theta0}
  x_1 \in \hood(\theta_0)\ \ \text{and}\ \ \theta_0 < \tfrac{\Delta}{2},\ \ \text{where}\ \ \Delta := \min\left\{ \bar\Delta , \min_{i\in[n]} (u_i - l_i) \right\} \in \R{}_{>0}
\eequation
for some $\bar\Delta \in \R{}_{>0}$, where $\bar\Delta$ is introduced merely to ensure that $\Delta$ is finite.

The search direction computation in our proposed algorithm is the main aspect that distinguishes between the deterministic and stochastic settings.  Specifically, letting $g_k$ denote a stochastic gradient estimate computed with respect to $x_k$ (see Section~\ref{sec.stochastic}), we denote the gradient (estimate) for the barrier-augmented function by
\bequation\label{def.q}
  q_k := \left\{ \baligned \nabla f(x_k) - \mu_k \diag(x_k - l)^{-1} \ones + \mu_k \diag(u - x_k)^{-1} \ones && \text{(deterministic)}\phantom{,} \\ g_k - \mu_k \diag(x_k - l)^{-1} \ones + \mu_k \diag(u - x_k)^{-1} \ones && \text{(stochastic)}. \ealigned \right.
\eequation
Then, for $(\lambda_{k,\min},\lambda_{k,\max}) \in \R{}_{>0} \times \R{}_{>0}$ with $\lambda_{k,\min} \leq \lambda_{k,\max}$ and diagonal $H_k \in \Smbb^n$ with $\lambda_{k,\max} I \succeq H_k \succeq \lambda_{k,\min} I$, the search direction $d_k \in \R{n}$ is $d_k = -H_k^{-1} q_k$.

A complete statement of our proposed interior-point method (IPM) for solving problem~\eqref{prob.opt} with prescribed parameter sequences (i.e., barrier-parameter sequence $\{\mu_k\}$, neighborhood-parameter sequence $\{\theta_k\}$, step-size-bound sequences $\{\alpha_{k,\max}\}$ and $\{\gamma_{k,\max}\}$, and eigenvalue-bound sequences $\{\lambda_{k,\min}\}$ and $\{\lambda_{k,\max}\}$) is stated as Algorithm~\ref{alg.ipm}.  We have written Algorithm~\ref{alg.ipm} in a generic manner that demonstrates flexibility in the required parameter sequences.  Our analyses in the next section prescribe additional rules for these sequences that lead to convergence guarantees.

\begin{algorithm}[ht]
  \caption{IPM with Prescribed Parameter Sequences}
  \label{alg.ipm}
  \begin{algorithmic}[1]
    \Require $\{\mu_k\} \searrow 0$; $\{\theta_k\} \searrow 0$; $\{\alpha_{k,\max}\} \subset \Rext{}_{>0}$; $\{\gamma_{k,\max}\} \subset (0,1]$; $\{\lambda_{k,\min}\} \subset \R{}_{>0}$ and $\{\lambda_{k,\max}\} \subset \R{}_{>0}$ such that $\lambda_{k,\min} \leq \lambda_{k,\max}$ for all $k \in \N{}$; and $x_1 \in \hood(\theta_0)$ for some $\theta_0 \in \R{}_{>\theta_1}$ satisfying \eqref{eq.theta0}
    \For{$k = 1,2,\dots$}
      \State choose diagonal $H_k \in \Smbb^n$ such that $\lambda_{k,\max} I \succeq H_k \succeq \lambda_{k,\min} I$\label{line:Hk}
      \State compute $d_k \gets -H_k^{-1} q_k$ \edit{for $q_k$ defined in \eqref{def.q}}
      \State choose $\alpha_k \in (0,\alpha_{k,\max}]$
      \State compute $\gamma_k \gets \max\{\gamma \in (0,\gamma_{k,\max}] : x_k + \gamma \alpha_k d_k \in \hood(\theta_k)\}$\label{line:gamma_def}
      \State set $x_{k+1} \gets x_k + \gamma_k \alpha_k d_k$\label{line:x-update}
    \EndFor
  \end{algorithmic}
\end{algorithm}

\bremark\label{rem.extend}
  One could extend our algorithm and analysis to allow, for all $k \in \N{}$, the employment of non-diagonal $H_k$ and/or the option to set $x_{k+1}$ by searching further along the piecewise linear path defined by the projections of $x_k + \gamma \alpha_k d_k$ onto $\hood(\theta_k)$ over $\gamma \in (0,\gamma_{k,\max}]$.  In such a setting, one needs to ensure that the barrier-augmented function decrease lemmas that appear in our analyses in Section~\ref{sec.deterministic} and \ref{sec.stochastic} $($namely, Lemmas~\ref{lem.decrease} and \ref{lem.decrease_stochastic}, respectively$)$ guarantee decreases of the same order in terms of the algorithmic parameters.  This can be done, for example, by ensuring that the angle between the resulting direction and $-q_k$ is acute and bounded away from 90$^\circ$ by a threshold that is independent of~$k$ and that the norm of the direction and~$-q_k$ are proportional uniformly over all $k \in \N{}$.  However, since such extensions would only obfuscate our analysis $($specifically, Lemma~\ref{lem.gammakbound}$)$ without adding significant value to our conclusions, we consider the simpler procedures in Algorithm~\ref{alg.ipm}, which has the features that would drive convergence in such algorithm variants as well.
\eremark

\bremark
  Algorithm~\ref{alg.ipm} is written as a primal interior-point method in the sense that each search direction is computed from an $n$-by-$n$ ``Newton-like'' system.  One could also consider a primal-dual interior-point method where the sequence of Lagrange multiplier estimates, say $\{(y_k,z_k)\}$ $($see \eqref{eq.KKT}$)$, act as independent components of the iterates.  To ensure convergence guarantees for a deterministic version of such a method, one can employ our strategies for the parameter sequences as long as safeguards are included to ensure that $y_k$ and $z_k$ remain within appropriately defined neighborhoods of $\mu_k \diag(x_k - l)^{-1} \ones$ and $\mu_k \diag(u - x_k)^{-1} \ones$, respectively, for all $k \in \N{}$.  Similar safeguards have been used in the literature; see, e.g., \cite[Section~2.2]{WaecBieg06}.  However, convergence guarantees for such a method in the stochastic setting do not follow readily from our analysis for Algorithm~\ref{alg.ipm}; hence, overall, we do not consider a primal-dual interior-point variant of Algorithm~\ref{alg.ipm} in this paper.
\eremark

\section{Convergence Analyses}\label{sec.analysis}

We analyze the behavior of Algorithm~\ref{alg.ipm} under Assumption~\ref{ass.f} as well as the following assumption.

\bassumption\label{ass.bounded}
  The iterate sequence $\{x_k\}$ of Algorithm~\ref{alg.ipm} is contained in an open set $\Xcal \subseteq \interior(\Bcal)$ over which distances of iterate components to finite bounds are bounded in the sense that, for some $\chi \in \R{}_{>1}$, one has for all $k \in \N{}$ that $x_{k,i} - l_i \leq \chi$ for all $i \in \Lcal$ and $u_i - x_{k,i} \leq \chi$ for all $i \in \Ucal$.
\eassumption

\noindent
The bounds required for Assumption~\ref{ass.bounded} to hold are not restrictive for practical purposes.  Indeed, while Assumption~\ref{ass.bounded} requires that \edit{$\chi \in \R{}_{>1}$} exists, it can be arbitrarily large and \emph{knowledge of it is not required by the algorithm}; see upcoming Lemma~\ref{lem.derivatives}.

\subsection{Preliminary Results}\label{sec.preliminary}

In this subsection, we provide preliminary results that are required for our analyses of Algorithm~\ref{alg.ipm} for the deterministic and stochastic settings, which are considered separately in the subsequent subsections.

Our first lemma essentially shows that derivatives of the barrier-augmented function are unaffected by scaling of the displacements from finite bounds that appear in the barrier function.  For the lemma, we introduce the function $\barrier : \interior(\Bcal) \times \R{}_{>0} \to \R{}$ that one obtains by scaling the barrier terms by $\chi^{-1}$ (see Assumption~\ref{ass.bounded}), namely,
\bequation\label{def.barrier_shifted}
  \barrier(x,\mu) = f(x) - \mu \sum_{i \in \Lcal} \log\((x_i - l_i)/\chi\) - \mu \sum_{i \in \Ucal} \log\((u_i - x_i)/\chi\).
\eequation
Important relationships between $\phi$ and $\barrier$ and the derivatives of these functions with respect to their first argument are the subject of this first lemma.

\blemma\label{lem.derivatives}
  For all $(x,\mu) \in \Xcal \times \R{}_{>0}$, one finds $\barrier(x,\mu) = \phi(x,\mu) + \mu M \geq f_{\inf}$, so $\nabla_x \phi(x,\mu) = \nabla_x \barrier(x,\mu)$, where $M \in \R{}_{>0}$ is independent of $x$ and $\mu$.  Moreover, for any $(\mu,\bar\mu) \in \R{}_{>0} \times \R{}_{>0}$ with $\bar\mu < \mu$, one has that $\barrier(x,\bar\mu) < \barrier(x,\mu)$ for all $x \in \Xcal$.
\elemma
\bproof
  The first desired equation follows from the definitions of $\phi(\cdot,\mu)$ and $\barrier(\cdot,\mu)$ in \eqref{def.barrier} and \eqref{def.barrier_shifted}, respectively, and the fact that, for any $\delta \in \R{}_{>0}$, one finds (since $\chi \in \R{}_{>1}$) that $-\log(\delta/\chi) = -\log(\delta) + \log(\chi)$.  Then, the fact that $\barrier(x,\mu) \geq f_{\inf}$ for all $(x,\mu) \in \Xcal \times \R{}_{>0}$ follows by Assumption~\ref{ass.f} and the fact that $(x_i - l_i)/\chi \in [0,1]$ for all $i \in \Lcal$ and $(u_i - x_i)/\chi \in [0,1]$ for all $i \in \Ucal$.  Next, the desired conclusions pertaining to the derivatives of $\phi(\cdot,\mu)$ and $\barrier(\cdot,\mu)$ follow from the first conclusion.  The final desired conclusion follows from the fact that, for all $x \in \Xcal$, one finds that $(x_i - l_i)/\chi \in [0,1]$ for all $i \in \Lcal$ and $(u_i - x_i)/\chi \in [0,1]$ for all $i \in \Ucal$.
\eproof

A consequence of Lemma~\ref{lem.derivatives} is that, for any $k \in \N{}$ such that the true gradient of the objective $\nabla f(x_k)$ is used in the definition of $q_k$, the search direction computation in Algorithm~\ref{alg.ipm} produces a descent direction for $\barrier(\cdot,\mu_k)$ from $x_k$ (recall that $H_k \succ 0$) \emph{even without explicit knowledge of the bound $\chi$ defined in Assumption~\ref{ass.bounded}}.

We now prove that the scaled barrier-augmented function has a gradient that satisfies a Lipschitz-continuity property over the line segment between any two points in a neighborhood of the type that is defined for the algorithm.  The result also shows that the corresponding Lipschitz constant depends on how close the elements of the points are to their corresponding lower and/or upper bounds, which, as shown in our subsequent analysis, is a fact that can be exploited by our algorithm.

\blemma\label{lem.Lipschitz}
  For any $(\mu,\theta,\bar\theta) \in \R{}_{>0} \times \R{}_{>0} \times \R{}_{>0}$ with $\bar\theta \in (0,\theta]$, $(x,\xbar) \in \hood(\theta) \times \hood(\bar\theta)$, and $\gamma \in [0,1]$, one finds that
  \bequation\label{eq.Lipschitz}
    \|\nabla_x \barrier(x + \gamma(\xbar - x),\mu) - \nabla_x \barrier(x,\mu)\|_2 \leq \gamma \ell_{\nabla f,\Bcal,\mu,x,\xbar} \|\xbar - x\|_2,
  \eequation
  where $\ell_{\nabla f,\Bcal,\mu,x,\xbar} := \ell_{\nabla f,\Bcal} + \mu a(x,\xbar)^{-1} + \mu b(x,\xbar)^{-1}$ with
  \begin{align*}
    a_i(x) &:= x_i - l_i; & a_i^+(x,\xbar) &:= \min\{x_i - l_i, \xbar_i - l_i\}; & a(x,\xbar) &:= \min_{i\in[n]} \{a_i(x) a_i^+(x,\xbar)\} \\
    b_i(x) &:= u_i - x_i; & b_i^+(x,\xbar) &:= \min\{u_i - x_i, u_i - \xbar_i\}; & b(x,\xbar) &:= \min_{i\in[n]} \{b_i(x) b_i^+(x,\xbar)\}.
  \end{align*}
  Moreover, one finds that $\ell_{\nabla f,\Bcal,\mu,x,\xbar} \leq \ell_{\nabla f,\Bcal,\mu,\theta,\bar\theta} := \ell_{\nabla f,\Bcal} + 2\mu\theta^{-1}\bar\theta^{-1} \in \R{}_{>0}$.
\elemma
\bproof
  For arbitrary such $(\mu,\theta,\bar\theta,x,\xbar,\gamma)$, \eqref{eq.barrier_stationary} and Lemma~\ref{lem.derivatives} imply
  \begin{align*}
    &\ \|\nabla_x \barrier(x + \gamma(\xbar - x),\mu) - \nabla_x \barrier(x,\mu)\|_2 \\
    \leq&\ \|\nabla f(x + \gamma(\xbar - x)) - \nabla f(x)\|_2 \\
    &\ + \mu \|((\diag(x - l) + \gamma \diag(\xbar - x))^{-1} - \diag(x - l)^{-1}) \ones\|_2 \\
    &\ + \mu \|((\diag(u - x) - \gamma \diag(\xbar - x))^{-1} - \diag(u - x)^{-1}) \ones\|_2.
  \end{align*}
  Considering the latter two terms, for arbitrary $i \in [n]$, one has
  \bequationNN
    \tfrac{1}{x_i + \gamma(\xbar_i - x_i) - l_i} - \tfrac{1}{x_i - l_i} = \tfrac{\gamma(x_i - \xbar_i)}{(x_i + \gamma(\xbar_i - x_i) - l_i)(x_i - l_i)} \leq \tfrac{\gamma(x_i - \xbar_i)}{a_i(x)a_i^+(x,\xbar)}
  \eequationNN
  and similarly that $(u_i - x_i - \gamma(\xbar_i - x_i))^{-1} - (u_i - x_i)^{-1} \leq \gamma(x_i - \xbar_i)b_i(x)^{-1}b_i^+(x,\xbar)^{-1}$.  By Assumption~\ref{ass.f} and these bounds over all $i \in [n]$, the desired conclusion follows.
\eproof

\bremark
  Lemma~\ref{lem.Lipschitz} and all of our subsequent analysis could be based on the simpler, but more conservative $\ell_{\nabla f,\Bcal,\mu,\theta,\bar\theta}$ rather than $\ell_{\nabla f,\Bcal,\mu,x,\xbar}$ in \eqref{eq.Lipschitz}.  Similarly, the step-size rule that we present in upcoming Parameter Rule~\ref{pm.step_size} could be based on the simpler, but more conservative $2\theta^{-1}\bar\theta^{-1}$ rather than $a(x,\xbar)^{-1} + b(x,\xbar)^{-1}$.  However, these tighter bounds have the effect of allowing larger step sizes, which is beneficial in the numerical experiments presented in Section~\ref{sec.numerical}.  Hence, we make these choices to have consistency between our analysis and numerical experimentation.
\eremark

The following consequence of Lemma~\ref{lem.Lipschitz} is central to our analysis.

\blemma\label{lem.quadratic_upper}
  For any $(\mu,\theta,\bar\theta) \in \R{}_{>0} \times \R{}_{>0} \times \R{}_{>0}$ with $\bar\theta \in (0,\theta]$ and $(x,\xbar) \in \hood(\theta) \times \hood(\bar\theta)$, one finds with $\ell_{\nabla f,\Bcal,\mu,x,\xbar} \in \R{}_{>0}$ defined in Lemma~\ref{lem.Lipschitz} that
  \bequationNN
    \barrier(\xbar,\mu) \leq \barrier(x,\mu) + \nabla_x \barrier(x,\mu)^T(\xbar - x) + \thalf \ell_{\nabla f,\Bcal,\mu,x,\xbar} \|\xbar - x\|_2^2.
  \eequationNN
\elemma
\bproof
  For arbitrary $(\mu,\theta,\bar\theta,x,\xbar)$ satisfying the conditions of the lemma, it follows with the Fundamental Theorem of Calculus and the Cauchy-Schwarz inequality that
  \begin{align*}
    &\ \barrier(\xbar, \mu) - \barrier(x, \mu) \\
    =&\ \int_0^1 \tfrac{\partial \barrier(x + \gamma (\xbar - x), \mu)}{\partial \gamma} \text{d}\gamma = \int_0^1 \nabla_x \barrier(x + \gamma(\xbar - x), \mu)^T(\xbar - x) \text{d}\gamma \\
    =&\ \nabla_x \barrier(x,\mu)^T(\xbar - x) + \int_0^1 (\nabla_x \barrier(x + \gamma(\xbar - x), \mu) - \nabla_x \barrier (x, \mu))^T(\xbar - x) \text{d}\gamma \\
    \leq&\ \nabla_x \barrier(x, \mu)^T(\xbar - x) + \|\xbar - x\|_2 \int_0^1 \|\nabla_x \barrier(x + \gamma(\xbar - x), \mu) - \nabla_x \barrier (x, \mu)\|_2 \text{d}\gamma.
  \end{align*}
  Hence, the desired conclusion follows along with Lemma~\ref{lem.Lipschitz} and since $\int_0^1 \gamma \text{d} \gamma = \thalf$.
\eproof

The prior lemma motivates the following parameter rule that we make going forward.  We remark at this stage that, for our analysis of the deterministic setting in the next subsection, one can consider $\alpha_{k,\max} \gets \infty$ for all $k \in \N{}$ so that the step size is always set to be the first term in the minimum in \eqref{eq.parameter_example}.  However, for the stochastic setting, our analysis requires a more conservative choice for $\{\alpha_{k,\max}\}$; see Section~\ref{sec.stochastic}.  Hence, we introduce $\{\alpha_{k,\max}\}$ at this stage, and carry it through our analysis, to maintain consistency between the deterministic and stochastic settings.  \edit{We also remark that while we aim to prove each of the results in our analysis in manners that keep them relatively generic, we ultimately show specific convergence guarantees for our algorithm when the barrier- and neighborhood-parameter sequences are defined respectively by $\mu_k = \mu_1 k^{t_\mu}$ and $\theta_{k-1} = \theta_0 k^{t_\theta}$ for all $k \in \N{}$ for some $(t_\mu,t_\theta) \in (-\infty,0) \times (-\infty,0)$.  Moreover, as seen in Parameter Rule~\ref{pm.step_size} below, we define the step-size sequence in terms of these sequences and a user-defined parameter $t_\alpha \in (-\infty,0]$.  Our analysis reveals critical distinctions between the deterministic and stochastic settings in terms of restrictions on the choice of $(t_\mu,t_\theta,t_\alpha)$; e.g., in the deterministic setting, our analysis ultimately requires $t_\mu = t_\theta$ and $t_\mu + t_\alpha \in [-1,0)$ (e.g., these are satisfied by $t_\mu = t_\theta = -1$ and $t_\alpha = 0$), and for the stochastic setting the choice needs to be restricted further.  We provide additional commentary on these distinctions in our analysis of the stochastic setting in Section~\ref{sec.stochastic}.}

\begin{parameterrule}\label{pm.step_size}
  \edit{Given $t_\alpha \in (-\infty,0]$, for} all $k \in \N{}$ the algorithm has
  \bequationNN
    \alpha_{k,\max} \geq \alpha_{k,\min},\ \ \text{where}\ \ \alpha_{k,\min} :=  \tfrac{\lambda_{k,\min} \edit{k^{t_\alpha}}}{\ell_{\nabla f,\Bcal} + 2\mu_k\theta_k^{-2}},
  \eequationNN
  and the algorithm sets
  \begin{align}
    \alpha_k &\gets \min\left\{\tfrac{\lambda_{k,\min} \edit{k^{t_\alpha}}}{\ell_{\nabla f,\Bcal,k}}, \alpha_{k,\max}\right\}, \label{eq.parameter_example} \\
    \text{where}\ \ \alpha_{k,\pre} &\gets \tfrac{\lambda_{k,\min} \edit{k^{t_\alpha}}}{\ell_{\nabla f,\Bcal} + \mu_k a(x_k,x_k)^{-1} + \mu_k b(x_k,x_k)^{-1}}, \nonumber \\
    \bar\gamma_k &\gets \max\{\gamma \in (0,\gamma_{k,\max}] : x_k + \gamma \alpha_{k,\pre} d_k \in \hood(\theta_k)\}, \nonumber \\
    \text{and}\ \ \ell_{\nabla f,\Bcal,k} &\gets \ell_{\nabla f,\Bcal} + \mu_k a(x_k,x_k + \bar\gamma_k \alpha_{k,\pre} d_k)^{-1} + \mu_k b(x_k,x_k + \bar\gamma_k \alpha_{k,\pre} d_k)^{-1}. \nonumber
  \end{align}
\end{parameterrule}

The following lemma shows that the step-size rule in Parameter Rule~\ref{pm.step_size} employs a denominator, namely, $\ell_{\nabla f,\Bcal,k}$, that serves as an upper bound for the Lipschitz constant seen in Lemmas~\ref{lem.Lipschitz} and \ref{lem.quadratic_upper}.  An implication of this fact is that the inequalities in these lemmas hold with that constant replaced by $\ell_{\nabla f,\Bcal,k}$, and another implication, stated in the lemma, is that the step size is contained in a prescribed interval.  \edit{The proof reveals that the important property of $\alpha_{k,\pre}$ is that it can be computed prior to~$\alpha_k$, yet is ensured to be an upper bound for the value of $\alpha_k$ that will be computed.}

\blemma\label{lem.Lipschitz_bound}
  For all $k \in \N{}$, with $\ell_{\nabla f,\Bcal,\mu_k,x_k,x_{k+1}}$ defined as in Lemma~\ref{lem.Lipschitz} and $\ell_{\nabla f,\Bcal,k}$ defined as in Parameter Rule~\ref{pm.step_size}, one finds that
  \bequation\label{eq.Lipschtiz_inequalities}
    \ell_{\nabla f, \Bcal, \mu_k, x_k, x_{k+1}} \leq \ell_{\nabla f,\Bcal,k} \leq \ell_{\nabla f,\Bcal,\mu_k,\theta_{k-1},\theta_k} \leq \ell_{\nabla f,\Bcal} + 2\mu_k\theta_k^{-2},
  \eequation
  from which it follows that the step size in Parameter Rule~\ref{pm.step_size} has $\alpha_k \in [\alpha_{k,\min},\alpha_{k,\max}]$.
\elemma
\bproof 
  Consider arbitrary $k \in \N{}$.  To prove the first inequality in \eqref{eq.Lipschtiz_inequalities}, it follows from the definitions in Lemma~\ref{lem.Lipschitz} that it is sufficient to prove that, for all $i \in [n]$,
  \bsubequations
    \begin{align}
      a_i^+(x_k, x_{k+1}) &\geq a_i^+(x_k, x_k + \bar\gamma_k \alpha_{k,\pre} d_k) \label{eq.compare_a_i_plus} \\
      \text{and}\ \ b_i^+(x_k, x_{k+1}) &\geq b_i^+(x_k, x_k + \bar\gamma_k \alpha_{k,\pre} d_k) \label{eq.compare_b_i_plus}.
    \end{align}
  \esubequations
  Toward this end, let us first show that $\gamma_k \alpha_k \leq \bar{\gamma}_k \alpha_{k,\pre}$.  Denoting (with $\min \emptyset = \infty$)
  \begin{align*}
    \delta_k^l(\alpha_k) &:= \min\left\{\tfrac{l_i + \theta_k - x_{k,i}}{\alpha_k d_{k,i}} : d_{k,i} < 0, i \in [n] \right\} \\ \text{and}\ \ 
    \delta_k^u(\alpha_k) &:= \min\left\{\tfrac{u_i - \theta_k - x_{k,i}}{\alpha_k d_{k,i}} : d_{k,i} > 0, i \in [n] \right\},
  \end{align*}
  the definition of $\gamma_k$ in line \ref{line:gamma_def} of~Algorithm \ref{alg.ipm} yields $\gamma_k = \min\{\gamma_{k,\max},\delta_k^l(\alpha_k),\delta_k^u(\alpha_k)\}$.  Thus, $\gamma_k \alpha_k = \min\{\gamma_{k,\max} \alpha_k,\delta_k^l(1),\delta_k^u(1)\}$.  Similarly, by the definition of $\bar{\gamma}_k$ in Parameter Rule~\ref{pm.step_size}, one finds that $\bar{\gamma}_k \alpha_{k,\pre} = \min\{\gamma_{k,\max} \alpha_{k,\pre}, \delta_k^l(1),\delta_k^u(1)\}$.  Hence, since $a(x_k, x_k + \bar\gamma_k \alpha_{k,\pre} d_k) \le a(x_k, x_k)$ and $b(x_k, x_k + \bar\gamma_k \alpha_{k,\pre}d_k) \le b(x_k, x_k)$ imply $\alpha_k \leq \alpha_{k,\pre}$, one finds $\gamma_k \alpha_k \leq \bar{\gamma}_k \alpha_{k,\pre}$, as desired.  Now, for $i \in [n]$ with $d_{k,i} < 0$,
  \begin{align*}
    a_i^+(x_k,x_{k+1})
      &=    x_{k+1,i} - l_i = x_{k,i} + \gamma_k \alpha_k d_{k,i} - l_i \\
      &\geq x_{k,i} + \bar{\gamma}_k \alpha_{k,\pre} d_{k,i} - l_i = a_i^+(x_k, x_k + \bar\gamma_k \alpha_{k,\pre} d_k),
  \end{align*}
  while for $i \in [n]$ with $d_{k,i} \geq 0$, $a_i^+(x_k, x_{k+1}) = x_{k, i} - l_i = a_i^+(x_k, x_k + \bar\gamma_k \alpha_{k,\pre} d_k)$.  Therefore, \eqref{eq.compare_a_i_plus} holds.  One finds that \eqref{eq.compare_b_i_plus} holds with a similar derivation.  Consequently, the first desired inequality in \eqref{eq.Lipschtiz_inequalities} holds.  The remaining desired inequalities in \eqref{eq.Lipschtiz_inequalities} follow by the definitions in Lemma~\ref{lem.Lipschitz}, Parameter Rule~\ref{pm.step_size}, and $\{\theta_k\} \searrow 0$.
\eproof

The step-size choice in Parameter Rule~\ref{pm.step_size} depends on the Lipschitz constant $\ell_{\nabla f,\Bcal}$ (amongst other quantities prescribed and/or computed in Algorithm~\ref{alg.ipm}), which can be any Lipschitz constant for~$\nabla f$ over~$\Bcal$ (i.e., it does not need to be the minimal such Lipschitz constant).  This choice is reasonable for practical purposes since such a value can be computed or estimated in practice.  Overall, this choice of step size in Parameter Rule~\ref{pm.step_size} can be viewed as a generalization of the \edit{$\ell_{\nabla f,\Bcal}^{-1}$}-type step-size rules common in unconstrained \edit{deterministic and stochastic} optimization.

\bremark\label{rem.careful}
  Observe that Lemma~\ref{lem.Lipschitz} reveals that convergence guarantees in the context of an interior-point method do not follow readily from the standard arguments for a gradient-based method for unconstrained optimization.  In particular, even though the lemma shows that the gradient of the $($scaled$)$ barrier-augmented function is Lipschitz continuous over $\hood(\theta)$ for any $\theta \in \R{}_{>0}$, the Lipschitz constant $($for a given $\mu \in \R{}_{>0}$$)$ can diverge as $\theta \searrow 0$.  Hence, ensuring convergence requires a careful balance between the barrier-parameter sequence, step-size sequence, and neighborhood-parameter sequence, as revealed in our subsequent analyses.
\eremark

\subsection{Deterministic Setting}\label{sec.deterministic}

We now focus on convergence guarantees that can be shown for Algorithm~\ref{alg.ipm} in the deterministic setting, i.e., when $q_k$ is computed using $\nabla f(x_k)$ for all $k \in \N{}$.  We begin by proving a set of generic results, then conclude with observations about specific choices for the parameter sequences that yield convergence guarantees with respect to stationarity measures.

We first provide a decrease lemma for the shifted barrier-augmented function, which is reminiscent of a standard decrease lemma for a gradient-based method in the context of unconstrained optimization.  The particular form of this result is a consequence of the choice of the step size stated in Parameter Rule~\ref{pm.step_size}.

\blemma\label{lem.decrease}
  For all $k \in \N{}$, one finds that
  \bequationNN
    \barrier(x_{k+1},\mu_{k+1}) - \barrier(x_k,\mu_k) \leq -\thalf \gamma_k \alpha_k \|\nabla_x \barrier(x_k,\mu_k)\|_{H_k^{-1}}^2.
  \eequationNN
\elemma
\bproof
  For all $k \in \N{}$, \edit{Lemmas~\ref{lem.derivatives}, \ref{lem.quadratic_upper}, and \ref{lem.Lipschitz_bound},} the value of $x_{k+1}$ in line~\ref{line:x-update} of Algorithm~\ref{alg.ipm}, and the conditions on $H_k$ in line~\ref{line:Hk} of Algorithm~\ref{alg.ipm} imply
  \begin{align*}
    &\ \barrier(x_{k+1},\mu_k) - \barrier(x_k,\mu_k) \\
      \leq&\ \nabla_x \barrier(x_k,\mu_k)^T(x_{k+1} - x_k) + \thalf \ell_{\nabla f,\Bcal,\mu_k,x_k,x_{k+1}} \|x_{k+1} - x_k\|_2^2 \\
      \leq&\ -\nabla_x \barrier(x_k,\mu_k)^T(\gamma_k \alpha_k H_k^{-1}\nabla_x \barrier(x_k,\mu_k)) + \thalf \ell_{\nabla f,\Bcal,k} \|\gamma_k \alpha_k H_k^{-1} \nabla_x \barrier(x_k,\mu_k)\|_2^2 \\
      \leq&\ - \gamma_k \alpha_k \|\nabla_x \barrier(x_k,\mu_k)\|_{H_k^{-1}}^2 + \thalf \gamma_k^2 \alpha_k^2 \lambda_{k,\min}^{-1} \ell_{\nabla f,\Bcal,k} \|\nabla_x \barrier(x_k,\mu_k)\|_{H_{k}^{-1}}^2 \\
      =&\ - \gamma_k \alpha_k (1- \thalf \gamma_k \alpha_k \lambda_{k,\min}^{-1} \ell_{\nabla f,\Bcal,k}) \|\nabla_x \barrier(x_k,\mu_k)\|_{H_k^{-1}}^2.
  \end{align*}
  Now, one finds under Parameter Rule~\ref{pm.step_size} that the parameter sequences yield
  \begin{align*}
    \alpha_k \ell_{\nabla f,\Bcal,k} \leq \lambda_{k,\min} \edit{k^{t_\alpha}} &\implies \edit{\gamma_k \alpha_k \ell_{\nabla f,\Bcal,k} \leq \lambda_{k,\min} k^{t_\alpha}} \\
    &\implies \gamma_k \alpha_k \ell_{\nabla f,\Bcal,k} \leq \lambda_{k,\min} \\
    &\iff \thalf \leq 1 - \thalf \gamma_k \alpha_k \lambda_{k,\min}^{-1} \ell_{\nabla f,\Bcal,k}.
  \end{align*}
  Thus, one finds from above that $\barrier(x_{k+1},\mu_k) - \barrier(x_k,\mu_k) \leq -\thalf \gamma_k \alpha_k \|\nabla_x \barrier(x_k,\mu_k)\|_{H_k^{-1}}^2$.
  Combining this inequality with Lemma~\ref{lem.derivatives}, which since $\mu_{k+1} < \mu_k$ shows that $\barrier(x_{k+1},\mu_{k+1}) < \barrier(x_{k+1},\mu_k)$, one reaches the desired conclusion.
\eproof

We now prove a critical lower bound on each element of the sequence $\{\gamma_k\}$.

\blemma\label{lem.gammakbound}
  For all $k \in \N{}$, define
  \bequation\label{eq.gamma_min}
    \gamma_{k,\min} := \min\left\{ 1, \tfrac{\lambda_{k,\min} \(\tfrac{\thalf \mu_k \Delta}{\mu_k + \thalf \kappa_{\nabla f,\Bcal,\infty} \Delta} - \theta_k\)}{\alpha_k (\kappa_{\nabla f,\Bcal,\infty} + \mu_k \theta_{k-1}^{-1})} \right\}
  \eequation
  and suppose $\gamma_{k,\max} \in [\gamma_{k,\min},1]$.  Then, for all $k \in \N{}$, $\gamma_k \geq \gamma_{k,\min}$.
\elemma
\bproof
  Recall that the algorithm ensures, for all $k \in \N{}$, that
  \bequation\label{eq:theta-feasible}
    x_k \in \hood(\theta_{k-1}) \iff l + \theta_{k-1} \leq x_k \leq u - \theta_{k-1}.
  \eequation
  For all $k \in \N{}$ and $i \in [n]$, let $\gamma_{k,i} := \max\{\gamma \in (0,\gamma_{k,\max}] : x_{k,i} + \gamma \alpha_k d_{k,i} \in [l_i+\theta_k,u_i - \theta_k]\}$ so that $\gamma_k \gets \min_{i\in[n]} \gamma_{k,i}$.  Considering arbitrary $k \in \N{}$ and $i \in [n]$, let us suppose that $d_{k,i} < 0$ and prove a lower bound on $\gamma_{k,i}$ that is independent of the index $i \in [n]$. One would find---although we omit details for brevity---that the same lower bound for $\gamma_{k,i}$ can be proved when $d_{k,i} > 0$ in a similar manner.  All of this, along with the fact that $\gamma_{k,i} = \gamma_{k,\max}$ when $d_{k,i} = 0$, leads to the desired conclusion.

  Consider arbitrary $k \in \N{}$ and $i \in [n]$ and, as previously stated, suppose that $d_{k,i} < 0$.  If~$i \notin \Lcal$, then it follows that $\gamma_{k,i} = \gamma_{k,\max}$.  Hence, we  proceed under the assumption that $i \in \Lcal$.  If $x_{k,i} + \gamma_{k,\max} \alpha_k d_{k,i} \geq l_i + \theta_k$, then $\gamma_{k,i} = \gamma_{k,\max}$.  Otherwise, the algorithm ensures $x_{k,i} + \gamma_{k,i} \alpha_k d_{k,i} = l_i + \theta_k$, and \eqref{eq:theta-feasible} and Assumption~\ref{ass.f} give
  \begin{align}
    \gamma_{k,i} &= \tfrac{x_{k,i} - l_i - \theta_k}{\alpha_k [H_k]_{i,i}^{-1} \(\nabla_i f(x_k) - \mu_k(x_{k,i} - l_i)^{-1} + \mu_k(u_i - x_{k,i})^{-1}\)} \nonumber \\
    &\geq \tfrac{x_{k,i} - l_i - \theta_k}{\alpha_k [H_k]_{i,i}^{-1} \(\nabla_i f(x_k) + \mu_k(u_i - x_{k,i})^{-1}\)} \geq \tfrac{\lambda_{k,\min} (x_{k,i} - l_i - \theta_k)}{\alpha_k (|\nabla_i f(x_k)| + \mu_k \theta_{k-1}^{-1})}. \label{eq.gamma_below}
  \end{align}
  The remainder of our analysis in this case hinges on providing a positive lower bound for $x_{k,i} - l_i$.  First, if $i \in \Lcal$ and $i \notin \Ucal$, then one finds that $d_{k,i} < 0$ means
  \begin{align*}
    [H_k]_{i,i}^{-1} \(-\nabla_i f(x_k) + \tfrac{\mu_k}{x_{k,i} - l_i}\) &< 0 \\
    \iff \(-\nabla_i f(x_k) + \tfrac{\mu_k}{x_{k,i} - l_i}\) &< 0 \iff \nabla_i f(x_k) > 0\ \ \text{and}\ \ x_{k,i} - l_i > \tfrac{\mu_k}{|\nabla_i f(x_k)|}.
  \end{align*}
  Second, if $i \in \Lcal \cup \Ucal$, then one finds with $\Delta_i := u_i - l_i$ that
  \begin{align*}
    && d_{k,i} = [H_k]_{i,i}^{-1} \(-\nabla_i f(x_k) + \tfrac{\mu_k}{x_{k,i} - l_i} - \tfrac{\mu_k}{u_i - x_{k,i}}\) &< 0 \\
    \iff && -\nabla_i f(x_k) + \tfrac{\mu_k}{x_{k,i} - l_i} - \tfrac{\mu_k}{\Delta_i - (x_{k,i} - l_i)} &< 0 \\
    \iff && \mu_k \(\tfrac{\Delta_i - 2(x_{k,i} - l_i)}{(x_{k,i} - l_i)(\Delta_i - (x_{k,i} - l_i))} \) &< \nabla_i f(x_k) \\
    \iff && \nabla_i f(x_k) (x_{k,i} - l_i)^2 - (2\mu_k + \nabla_i f(x_k) \Delta_i)(x_{k,i} - l_i) + \mu_k \Delta_i &< 0.
  \end{align*}
  Given this inequality, there are three subcases to consider depending on $\nabla_i f(x_k)$.
  \benumerate
    \item[(i)] Suppose $\nabla_i f(x_k) = 0$.  Then, $x_{k,i} - l_i > \thalf \Delta_i$.
    \item[(ii)] Suppose $\nabla_i f(x_k) < 0$.  Then, by the quadratic formula, it follows that
    \bsubequations
      \begin{align}
      x_{k,i} - l_i &< \tfrac{\mu_k}{\nabla_i f(x_k)} + \tfrac{\Delta_i}{2} - \sqrt{\tfrac{\mu_k^2}{(\nabla_i f(x_k))^2} + \tfrac{\Delta_i^2}{4}} \label{eq.gradpos_upper} \\ \text{or}\ \ 
      x_{k,i} - l_i &> \tfrac{\mu_k}{\nabla_i f(x_k)} + \tfrac{\Delta_i}{2} + \sqrt{\tfrac{\mu_k^2}{(\nabla_i f(x_k))^2} + \tfrac{\Delta_i^2}{4}}. \label{eq.gradpos_lower}
      \end{align}
    \esubequations
    In fact, the upper bound on $x_{k,i} - l_i$ stated in \eqref{eq.gradpos_upper} is not possible since
    \begin{align*}
      &\ \tfrac{\mu_k}{\nabla_i f(x_k)} + \tfrac{\Delta_i}{2} - \sqrt{\tfrac{\mu_k^2}{(\nabla_i f(x_k))^2} + \tfrac{\Delta_i^2}{4}} \\
      \leq&\ \tfrac{\mu_k}{\nabla_i f(x_k)} + \tfrac{\Delta_i}{2} - \sqrt{\tfrac{\mu_k^2}{(\nabla_i f(x_k))^2} + 2 \tfrac{\mu_k}{\nabla_i f(x_k)}\tfrac{\Delta_i}{2} + \tfrac{\Delta_i^2}{4}} \\
      \leq&\ \tfrac{\mu_k}{\nabla_i f(x_k)} + \tfrac{\Delta_i}{2} - \left| \tfrac{\mu_k}{\nabla_i f(x_k)} + \tfrac{\Delta_i}{2} \right| \leq 0
    \end{align*}
    while the algorithm ensures $x_{k,i} - l_i > 0$.  Hence, \eqref{eq.gradpos_lower} must hold in this case, from which it follows (by dropping the $\tfrac14 \Delta_i^2$ term) that $x_{k,i} - l_i > \thalf \Delta_i$.
    \item[(iii)] Suppose $\nabla_i f(x_k) > 0$.  Then, by the quadratic formula, it follows that
    \bequation\label{eq.gradneg_lower}
      x_{k,i} - l_i > \tfrac{\mu_k}{\nabla_i f(x_k)} + \tfrac{\Delta_i}{2} - \sqrt{\tfrac{\mu_k^2}{(\nabla_i f(x_k))^2} + \tfrac{\Delta_i^2}{4}}.
    \eequation
    Define $a_{k,i} := \tfrac{\mu_k^2}{(\nabla_i f(x_k))^2} + 2 \tfrac{\mu_k}{\nabla_i f(x_k)} \tfrac{\Delta_i}{2} + \tfrac{\Delta_i^2}{4}$ and $b_{k,i} := \tfrac{\mu_k^2}{(\nabla_i f(x_k))^2} + \tfrac{\Delta_i^2}{4}$, and observe that $a_{k,i} > b_{k,i} > 0$ while the right-hand side of \eqref{eq.gradneg_lower} is equal to $s(a_{k,i}) - s(b_{k,i})$, where $s(\cdot) := \sqrt{\cdot}$ is the square root function.  By the mean value theorem, there exists a real number $c \in [b_{k,i},a_{k,i}]$ such that one finds
    \bequationNN
      s(a_{k,i}) - s(b_{k,i}) = s'(c)(a_{k,i} - b_{k,i}),\ \ \text{where}\ \ s'(c) = \tfrac{1}{2\sqrt{c}} \geq \tfrac{1}{2\sqrt{a_{k,i}}}.
    \eequationNN
    Hence, one has from \eqref{eq.gradneg_lower} that
  \eenumerate
  \begin{align*}
    x_{k,i} - l_i > s(a_{k,i}) - s(b_{k,i}) \geq \tfrac{a_{k,i} - b_{k,i}}{2 \sqrt{a_{k,i}}} &= \tfrac{\tfrac{\mu_k}{\nabla_i f(x_k)} \tfrac{\Delta_i}{2}}{\tfrac{\mu_k}{\nabla_i f(x_k)} + \tfrac{\Delta_i}{2}} = \tfrac{\thalf \mu_k \Delta_i}{\mu_k + \thalf \nabla_i f(x_k) \Delta_i}.
  \end{align*}
  Combining the results above when $i \in \Lcal$, one finds that $d_{k,i} < 0$ implies
  \bequationNN
    x_{k,i} - l_i \geq \min \left\{ \tfrac{\mu_k}{|\nabla_i f(x_k)|}, \tfrac{\Delta_i}{2}, \tfrac{\thalf \mu_k \Delta_i}{\mu_k + \thalf |\nabla_i f(x_k)| \Delta_i} \right\} = \tfrac{\thalf \mu_k \Delta_i}{\mu_k + \thalf |\nabla_i f(x_k)| \Delta_i}.
  \eequationNN
  Combining this inequality, \eqref{eq.gamma_below}, the facts that $\max_{i\in[n]} |\nabla_i f(x_k)| \leq \kappa_{\nabla f,\Bcal,\infty}$ and $\Delta \leq \min_{i\in[n]} \Delta_i$, and the monotonicity of $\frac{\rho z}{\tau + \omega z}$ with respect to $z$ when $\rho$, $\tau$, and $\omega$ are positive, one reaches the desired conclusion.
\eproof

The prior lemma motivates the following rule that we make going forward.  Similarly as for the choice of the step size (recall Parameter Rule~\ref{pm.step_size}), the remainder of our analysis for the deterministic setting can use $\gamma_{k,\max} \gets 1$ for all $k \in \N{}$, but for the stochastic setting our analysis requires a more conservative choice for $\{\gamma_{k,\max}\}$.

\begin{parameterrule}\label{pm.gamma}
  For all $k \in \N{}$, with $\gamma_{k,\min}$ from \eqref{eq.gamma_min}, $\gamma_{k,\max} \in [\gamma_{k,\min},1]$.
\end{parameterrule}

We now prove a generic convergence theorem for Algorithm~\ref{alg.ipm} in a deterministic setting.  We follow this theorem with a corollary that provides specific choices of the parameter sequences that ensure that the conditions of the theorem hold.

\btheorem\label{th.deterministic}
  Suppose that Assumptions~\ref{ass.f} and \ref{ass.bounded} and Parameter Rules \ref{pm.step_size} and \ref{pm.gamma} hold.  If, further, the parameter sequences of Algorithm~\ref{alg.ipm} yield
  \bequation\label{eq.infinite_series}
    \sum_{k=1}^\infty \gamma_k \alpha_k = \infty,
  \eequation
  then
  \bequation\label{eq.liminf_Hinv}
    \liminf_{k\to\infty} \|\nabla_x \barrier(x_k,\mu_k)\|_{H_k^{-1}}^2 = 0,
  \eequation
  meaning that if there exists $\overline{r} \in \R{}_{>0}$ such that $\lambda_{k,\max} \leq \overline{r}$ for all $k \in \N{}$, then
  \bequation\label{eq.liminf}
    \liminf_{k\to\infty} \|\nabla_x \barrier(x_k,\mu_k)\|_2^2 = 0.
  \eequation
  \edit{Now suppose in addition that such $\overline{r}$ exists, the sequence $\{\mu_k \theta_{k-1}^{-1}\}$ is bounded, and there exists a set $\Kcal \subseteq \N{}$ of infinite cardinality such that $\{\nabla \barrier(x_k,\mu_k)\}_{k\in\Kcal} \to 0$ and $\{x_k\}_{k\in\Kcal} \to \xbar$ for some $\xbar \in \Bcal$.  Then,} the limit point $\xbar$ is a KKT point for \eqref{prob.opt} in the sense that there exists $(\ybar,\zbar) \in \R{n} \times \R{n}$ such that $(\xbar,\ybar,\zbar)$ satisfies \eqref{eq.KKT}.
\etheorem
\bproof
  It follows from Lemma~\ref{lem.derivatives} that $\barrier$ is bounded below by $f_{\inf}$ over $\Xcal \times \R{}_{>0}$. Then, one finds by summing the expression in Lemma~\ref{lem.decrease} over $k \in \N{}$ that
  \begin{align*}
    \infty > \barrier(x_1,\mu_1) - f_{\inf} &\geq \sum_{k=1}^\infty (\barrier(x_k,\mu_k) - \barrier(x_{k+1},\mu_{k+1})) \geq \sum_{k=1}^\infty \tfrac{\gamma_k \alpha_k}{2} \|\nabla_x \barrier(x_k,\mu_k)\|_{H_k^{-1}}^2.
  \end{align*}
  If there exists $\epsilon \in \R{}_{>0}$ and $k_\epsilon \in \N{}$ such that $\|\nabla_x \barrier(x_k,\mu_k)\|_{H_k^{-1}}^2 \geq \epsilon$ for all $k \in \N{}$ with $k \geq k_\epsilon$, then the conclusion above contradicts \eqref{eq.infinite_series}.  Hence, it follows that such~$\epsilon$ and~$k_\epsilon$ do not exist, meaning that \eqref{eq.liminf_Hinv} holds, as desired.  Now, if there exists $\overline{r} \in \R{}_{>0}$ such that $\lambda_{k,\max} \leq \overline{r}$ for all $k \in \N{}$, then $0 = \liminf_{k\to\infty} \|\nabla_x \barrier(x_k,\mu_k)\|_{H_k^{-1}}^2 \geq \liminf_{k\to\infty} \overline{r}^{-1}\|\nabla_x \barrier(x_k,\mu_k)\|_{2}^2$, from which \eqref{eq.liminf} holds, as desired.  Now suppose that such $\overline{r}$ exists, $\{\mu_k \theta_{k-1}^{-1}\}$ is bounded, and there exists an infinite-cardinality set $\Kcal \subseteq \N{}$ as described in the theorem.  By Lemma~\ref{lem.derivatives}, it follows that $\{\nabla \phi(x_k,\mu_k)\}_{k\in\Kcal} \to 0$ as well.   Using this limit and, for all $k \in \Kcal$, defining the auxiliary sequences 
  \bequation\label{def:yk}
    y_k := \mu_k \diag(x_k - l)^{-1} \ones\ \ \text{and}\ \ z_k := \mu_k \diag(u - x_k)^{-1} \ones,
  \eequation
  it follows that $\{(x_k,y_k,z_k)\}_{k\in\Kcal} \subset \R{n} \times \R{n} \times \R{n}$ satisfies
  \bequation\label{eq:xk-to-xbar}
      \{x_k\}_{k\in\Kcal} \to \xbar \ \ \text{and} \ \
      \{\|\nabla f(x_k) - y_k + z_k\|_2\}_{k\in\Kcal} \to 0.
  \eequation
  Next, for all $k\in\N{}$, it follows from $x_k\in \hood(\theta_{k-1})$ (see Algorithm~\ref{alg.ipm}) that one has $0 \leq y_{k,i} = \tfrac{\mu_k}{x_{k,i} - l_i} \leq \tfrac{\mu_k}{\theta_{k-1}}$ and $0 \leq z_{k,i} = \tfrac{\mu_k}{u_i-x_{k,i}} \leq \tfrac{\mu_k}{\theta_{k-1}}$. Since $\{\mu_k \theta_{k-1}^{-1}\}$ is bounded by assumption, it follows that $\{y_k\}_{k\in\Kcal}$ and $\{z_k\}_{k\in\Kcal}$ are bounded.  Then, the Bolzano-Weierstrass Theorem gives the existence of an infinite subsequence of indices $K_{y,z} \subseteq K$ and vectors $\ybar \in \R{n}$ and $\zbar \in \R{n}$ such that 
  \bequation\label{eq:y-conv}
    \{y_k\}_{k \in K_{y,z}} \to \ybar\ \ \text{and} \ \ \{z_k\}_{k \in K_{y,z}} \to \zbar.
  \eequation
  Using these limits, \eqref{def:yk}, and $\{\mu_k\}\searrow 0$, it follows that 
  \bequation\label{eq:y-is-zero}
    \baligned
      \ybar_i &= 0\ \ \text{for all $i \in [n]$ with $\xbar_i \neq l_i$} \\
      \text{and}\ \ \zbar_i &= 0\ \ \text{for all $i \in [n]$ with $\xbar_i \neq u_i$}.
    \ealigned
  \eequation  
  Combining $\xbar\in\Bcal$, $K_{y,z} \subseteq K$, and \eqref{def:yk}--\eqref{eq:y-is-zero}, it follows that $\xbar$ is a KKT point for~\eqref{prob.opt} since the tuple $(\xbar,\ybar,\zbar)$ satisfies~\eqref{eq.KKT}, thus completing the proof.
\eproof

The following corollary shows that there exist choices of the parameter sequences such that the conditions of Theorem~\ref{th.deterministic} hold.  \edit{For the sake of generality, we introduce the distinct constants $t_\mu$ and $t_\theta$ for the barrier and neighborhood parameter sequences, respectively.  However, the proof shows that we require $t_\mu = t_\theta$, as stated in the corollary.  The proof could be written more concisely using this equation, but in order to show why this equation is needed, we use extra steps in the proof to show why in certain places $t_\mu \geq t_\theta$ is required, whereas $t_\mu \leq t_\theta$ is required in others.}

\bcorollary\label{co.sequence_choice}
  Suppose that Assumptions~\ref{ass.f} and \ref{ass.bounded} and Parameter Rules~\ref{pm.step_size} and \ref{pm.gamma} hold.  Then, there exist parameter choices for Algorithm~\ref{alg.ipm} such that the infinite series in \eqref{eq.infinite_series} is unbounded and $\{\mu_k \theta_{k-1}^{-1}\}$ is bounded; e.g., these consequences follow if for some $\underline{r} \in \R{}_{>0}$, \edit{$(t_\mu,t_\theta,t_\alpha) \in (-\infty,0) \times (-\infty,0) \times (-\infty,0]$ with $t_\mu = t_\theta$ and $t_\mu + t_\alpha \in [-1,0)$}, and $\mu_1 \in \R{}_{>0}$ with $\mu_1 > \tfrac{\thalf \theta_0 \kappa_{\nabla f,\Bcal,\infty} \Delta}{\thalf \Delta - \theta_0}$ the algorithm has $\mu_k = \mu_1 k^{\edit{t_\mu}}$, $\theta_{k-1} = \theta_0 k^{\edit{t_\theta}}$, $\alpha_{k,\max} \gets \infty$, $\gamma_{k,\max} \gets 1$, and $\underline{r} \leq \lambda_{k,\min} \leq \lambda_{k,\max}$ for all $k \in \N{}$.  Thus, with these choices, the lower limit in \eqref{eq.liminf_Hinv} holds, and if there exists $\overline{r} \in \R{}_{\geq\underline{r}}$ such that $\lambda_{k,\max} \leq \overline{r}$ for all $k \in \N{}$, then the lower limit in \eqref{eq.liminf} holds.  Finally, if all of the aforementioned choices of the parameter sequences are made and there exists an infinite-cardinality set $\Kcal \subseteq \N{}$ such that $\{\nabla \barrier(x_k,\mu_k)\}_{k\in\Kcal} \to 0$ and $\{x_k\}_{k\in\Kcal} \to \xbar$ for some $\xbar \in \Bcal$, then the limit point $\xbar$ is a KKT point for \eqref{prob.opt}.
\ecorollary
\bproof
  Under Parameter Rules~\ref{pm.step_size} and \ref{pm.gamma}, Lemmas~\ref{lem.Lipschitz_bound} and \ref{lem.gammakbound} imply that with the parameter choices given in the corollary, one finds that
  \begin{align}
    \gamma_k \alpha_k
    &\geq \underline{r} \min\left\{ \tfrac{\edit{k^{t_\alpha}}}{\ell_{\nabla f,\Bcal} + 2 \mu_1 \theta_0^{-2} k^{\edit{t_\mu}} (k+1)^{-2\edit{t_\theta}}} , \tfrac{ \tfrac{\thalf \mu_1 \Delta k^{\edit{t_\mu}}}{\mu_1 k^{\edit{t_\mu}} + \thalf \kappa_{\nabla f, \Bcal,\infty} \Delta} - \theta_0 (k+1)^{\edit{t_\theta}} }{ \kappa_{\nabla f,\Bcal,\infty} + \mu_1 \theta_0^{-1} \edit{k^{t_\mu} k^{-t_\theta}} } \right\} \nonumber \\
    &=: \underline{r}\min\{ \beta_k, \eta_k \}. \label{def:betak-alphak}
  \end{align}
  With respect to the sequence $\{\beta_k\}$, one finds for all $k \in \N{}$ that
  \bequation\label{betak-bound}
    k^{\edit{t_\mu}} (k+1)^{-2\edit{t_\theta}} \leq k^{\edit{t_\mu}} (2k)^{-2\edit{t_\theta}} = 2^{-2\edit{t_\theta}} k^{\edit{t_\mu - 2t_\theta}}.
  \eequation
  \edit{Since $t_\theta \leq t_\mu < 0$ implies $-t_\mu + 2t_\theta \leq 0$, it follows with \eqref{betak-bound} that there exists $c \in \R{}_{>0}$ such that $\beta_k \geq c k^{-t_\mu + 2t_\theta + t_\alpha}$ for all $k \in \N{}$, which in turn implies with $t_\theta \geq t_\mu$ that $\beta_k \geq c k^{t_\mu + t_\alpha}$ for all $k \in \N{}$.}  For the sequence $\{\eta_k\}$, one finds with \eqref{eq.theta0} and since \edit{$t_\mu \leq 0$, $t_\theta \leq 0$, and $t_\mu \geq t_\theta$} that, for all $k \in \N{}$,
  \begin{align*}
    \tfrac{\thalf \mu_1 \Delta k^{\edit{t_\mu}}}{\mu_1 k^{\edit{t_\mu}} + \thalf \kappa_{\nabla f, \Bcal,\infty} \Delta} - \theta_0 (k+1)^{\edit{t_\theta}}
    &\geq \tfrac{\thalf \mu_1 \Delta k^{\edit{t_\mu}}}{\mu_1 k^{\edit{t_\mu}} + \thalf \kappa_{\nabla f, \Bcal,\infty} \Delta} - \theta_0 k^{\edit{t_\theta}} \nonumber \\
    &= \( \tfrac{\thalf \mu_1 \Delta}{\mu_1 k^{\edit{t_\mu}} + \thalf \kappa_{\nabla f,\Bcal,\infty} \Delta} - \theta_0 k^{\edit{t_\theta - t_\mu}} \) k^{\edit{t_\mu}} \nonumber \\
    &\geq  \( \tfrac{\thalf \mu_1 \Delta}{\mu_1 + \thalf \kappa_{\nabla f,\Bcal,\infty} \Delta} - \theta_0 \) k^{\edit{t_\mu}}, 
  \end{align*}
  where one finds that \eqref{eq.theta0} (i.e., $\theta_0 < \tfrac\Delta2$) and
  \bequationNN
    \mu_1 > \tfrac{\thalf \theta_0 \kappa_{\nabla f,\Bcal,\infty} \Delta}{\thalf \Delta - \theta_0}\ \ \text{imply}\ \ \tfrac{\thalf \mu_1 \Delta}{\mu_1 + \thalf \kappa_{\nabla f,\Bcal,\infty} \Delta} > \theta_0\edit{,}
  \eequationNN
  \edit{which since $t_\theta \leq t_\mu$ means there exists $\cbar \in \R{}_{>0}$ such that $\eta_k \geq \cbar k^{t_\mu - (t_\mu - t_\theta)} = \cbar k^{t_\theta}$ for all $k \in \N{}$.  Combining \eqref{def:betak-alphak}, these lower bounds for $\{\beta_k\}$ and $\{\eta_k\}$, and the facts that $t_\theta \geq t_\mu$ and $t_\alpha \leq 0$, it follows that there exists $\overline{\overline{c}} \in \R{}_{>0}$ such that}
  \bequation\label{eq.alpha_rate}
    \gamma_k \alpha_k \geq \underline{r}\min\{\beta_k,\eta_k\} \geq \underline{r} \edit{\overline{\overline{c}} k^{t_\mu+t_\alpha}},
  \eequation
  and since \edit{$t_\mu + t_\alpha \in [-1,0)$}, one concludes that the infinite series in \eqref{eq.infinite_series} is unbounded.
  
  Finally, for all $k \in \N{}$, it follows from the given parameter choices that \edit{$\mu_k \theta_{k-1}^{-1} = \mu_1 \theta_0^{-1} k^{t_\mu - t_\theta}$ for all $k \in \N{}$, so $t_\theta \geq t_\mu$ ensures $\{\mu_k \theta_{k-1}^{-1}\}$ is bounded, as claimed.}
\eproof

\bremark\label{rem.carbonetto}
  Related to Remark~\ref{rem.careful}, we observe that the claimed asymptotic convergence guarantee in \cite{NIPS2008_a87ff679} overlooks a critical issue in terms of the step sizes.  The method in \cite{NIPS2008_a87ff679} employs step sizes that, amongst other considerations, ensure that each iterate remains feasible.  The claimed convergence guarantee then assumes that the step sizes are unsummable.  However, when the step sizes need to be reduced to maintain feasibility, one cannot presume that the resulting step sizes are unsummable.  This issue is overlooked in \cite{NIPS2008_a87ff679}, but---through a careful balance between the barrier-parameter, step-size, and neighborhood parameter sequences in our algorithm---we find in Corollary~\ref{co.sequence_choice} that there exist parameter choices such that \eqref{eq.infinite_series} holds.  We also carry this property forward for the stochastic setting, as seen in Section~\ref{sec.stochastic}.
\eremark

We conclude this subsection by observing that if $\Bcal$ is bounded, then $\{x_k\}$ has a convergent subsequence and, under the stated conditions in Corollary~\ref{co.sequence_choice}, an infinite-cardinality set~$\Kcal$ of the type described in the corollary is guaranteed to exist.

\subsection{Stochastic Setting}\label{sec.stochastic}

We now provide a convergence guarantee for Algorithm~\ref{alg.ipm} in a stochastic setting when, in every run for all $k \in \N{}$, $q_k$ is computed using an unbiased stochastic gradient estimate $g_k$ with bounded error; see upcoming Assumption~\ref{ass.stochastic}.  Formally, we consider the stochastic process defined by the algorithm, namely, $\{(X_k, G_k, Q_k, H_k, D_k, \overline\Gamma_k, A_k, \Gamma_k)\}$, where, for all $k \in \N{}$, the random variables correspond to the iterate~$X_k$, stochastic gradient estimator~$G_k$, stochastic barrier-augmented function gradient~$Q_k$, scaling matrix~$H_k$, direction~$D_k$, neighborhood enforcement parameter~$\overline\Gamma_k$, step size~$A_k$, and neighborhood enforcement parameter $\Gamma_k$.  A realization of this process is $\{(x_k, g_k, q_k, H_k, d_k, \bar\gamma_k, \alpha_k, \gamma_k)\}$, as in Algorithm~\ref{alg.ipm}.  (Here, we have introduced a slight abuse of notation in terms of $H_k$, which acts as both a random variable and its realization.  We prefer this slightly abused notation rather than introduce additional notation; it should not lead to confusion since, for our analysis in this subsection---which considers $H_k$ as a random variable for all $k \in \N{}$---ultimately relies on the fact that the eigenvalues of the elements of $\{H_k\}$ can be bounded by the prescribed bound sequences $\{\lambda_{k,\min}\}$ and $\{\lambda_{k,\max}\}$.)  The behavior of any run of the algorithm is determined by the initial conditions (including that $X_1 = x_1$) and the sequence of stochastic gradient estimators $\{G_k\}$.  Let~$\Fcal_1$ denote the $\sigma$-algebra defined by the initial conditions and, for all $k \in \N{}$ with $k \geq 2$, let $\Fcal_k$ denote the $\sigma$-algebra defined by the initial conditions and the random variables $\{G_1,\dots,G_{k-1}\}$, a realization of which determines the realizations of $\{X_j\}_{j=1}^k$ and $\{(G_j,Q_j,D_j,\overline\Gamma_j,A_j,\Gamma_j)\}_{j=1}^{k-1}$.  In this manner, the sequence $\{\Fcal_k\}$ is a filtration.

For our analysis in this subsection, we continue to make Assumptions~\ref{ass.f} and \ref{ass.bounded}, where for Assumption~\ref{ass.bounded} we assume that the set $\Xcal$ and real number~$\chi$ are uniform over all possible realizations of the stochastic process.  In terms of the stochastic gradient estimators and scaling matrices, we make the following assumption.

\begin{assumption}\label{ass.stochastic}
  For all $k \in \N{}$, one has $\E[G_k | \Fcal_k] = \nabla f(X_k)$.  In addition, there exists $(\sigma_2,\sigma_\infty) \in \R{}_{\geq0} \times \R{}_{\geq0}$ such that, for all $k \in \N{}$, one has
  \bequationNN
    \|G_k - \nabla f(X_k)\|_2 \leq \sigma_2\ \ \text{and}\ \ \|G_k - \nabla f(X_k)\|_\infty \leq \sigma_\infty.
  \eequationNN
  Finally, for all $k \in \N{}$, the matrix $H_k \in \mathbb{S}^n$ is $\Fcal_k$-\edit{measurable}.
\end{assumption}

\noindent
In Assumption~\ref{ass.stochastic}, the existence of $\sigma_\infty$ follows from that of $\sigma_2$, and vice versa, but we introduce both of these values for the sake of notational convenience.  It follows under Assumption~\ref{ass.stochastic} that, for all $k \in \N{}$, one has $\|G_k\|_\infty \leq \kappa_{\nabla f,\Bcal,\infty} + \sigma_\infty$.

In terms of the algorithmic choices that we consider in our present analysis, we state the following rule that can be seen as a slightly modified combination of Parameter Rules~\ref{pm.step_size} and \ref{pm.gamma}.  Overall, unlike in the deterministic setting where one can prove convergence guarantees (see Theorem~\ref{th.deterministic} and Corollary~\ref{co.sequence_choice}) using $\alpha_{k,\max} \gets \infty$ and $\gamma_{k,\max} \gets 1$ for all $k \in \N{}$ and without knowledge of $\kappa_{\nabla f,\Bcal,\infty}$ appearing in the definition of $\gamma_{k,\min}$ in Lemma~\ref{lem.gammakbound}, for the stochastic setting our analysis requires more conservative choices and information.  \edit{The rule employs what we refer to as \emph{buffer} sequences, namely, $\{\alpha_{k,\buff}\} \subset \R{}_{\geq0}$ and $\{\gamma_{k,\buff}\} \subset \R{}_{\geq0}$, to allow the algorithm to choose $\alpha_k > \alpha_{k,\min}$ and $\gamma_k > \gamma_{k,\min}$ for all $k \in \N{}$, which may lead to better practical performance.  A user may choose to set $\alpha_{k,\buff} = \gamma_{k,\buff} = 0$ for all $k \in \N{}$ for simplicity, if desired.  For consistency, like in the deterministic setting, we introduce the tuple of parameters $(t_\mu, t_\theta, t_\alpha)$, but since it is ultimately required for the stochastic setting like it is in the deterministic setting, we impose upfront that $t_\mu = t_\theta$.}

\begin{parameterrule}\label{pm.step_size_gamma}
  \edit{Given prescribed $(t_\mu,t_\theta,t_\alpha) \in (-\infty,-\thalf) \times (-\infty,-\thalf) \times (-\infty,0)$ such that $t_\mu = t_\theta$, $t_\mu + t_\alpha \in [-1,0)$, and $t_\mu + 2t_\alpha \in (-\infty,-1)$ along with prescribed $\alpha_{\buff} \in \R{}_{\geq0}$, $\{\alpha_{k,\buff}\} \subset \R{}_{\geq0}$, $\gamma_{\buff} \in \R{}_{\geq0}$, and $\{\gamma_{k,\buff}\} \subset \R{}_{\geq0}$ such that $\alpha_{k,\buff} \leq \alpha_{\buff} k^{2t_\mu}$ and $\gamma_{k,\buff} \leq \gamma_{\buff} k^{t_\mu}$ for all $k \in \N{}$, the algorithm employs}
  \begin{align*}
    \alpha_{k,\min} &:= \tfrac{\lambda_{k,\min}\edit{k^{t_\alpha}}}{\ell_{\nabla f,\Bcal} + 2\mu_k\theta_k^{-2}}, && \gamma_{k,\min} := \min\left\{ 1, \tfrac{\lambda_{k,\min} \(\tfrac{\thalf \mu_k \Delta}{\mu_k + \thalf (\kappa_{\nabla f,\Bcal,\infty} + \sigma_\infty) \Delta} - \theta_k\)}{\alpha_{k,\max} (\kappa_{\nabla f,\Bcal,\infty} + \sigma_\infty + \mu_k \theta_{k-1}^{-1})} \right\}, \\
    \alpha_{k,\max} &:= \alpha_{k,\min} + \alpha_{k,\buff}, && \text{and}\ \ \gamma_{k,\max} := \min\{1, \gamma_{k,\min} + \gamma_{k,\buff}\}
  \end{align*}
  and makes the $($run-and-iterate-dependent$)$ choice $\alpha_k \gets \min\left\{ \tfrac{\lambda_{k,\min}\edit{k^{t_\alpha}}}{\ell_{\nabla f,\Bcal,k}}, \alpha_{k,\max} \right\}$ \edit{for all $k \in \N{}$, where $\ell_{\nabla f,\Bcal,k} \in \R{}_{>0}$ is defined as in Parameter Rule~\ref{pm.step_size}.}
\end{parameterrule}

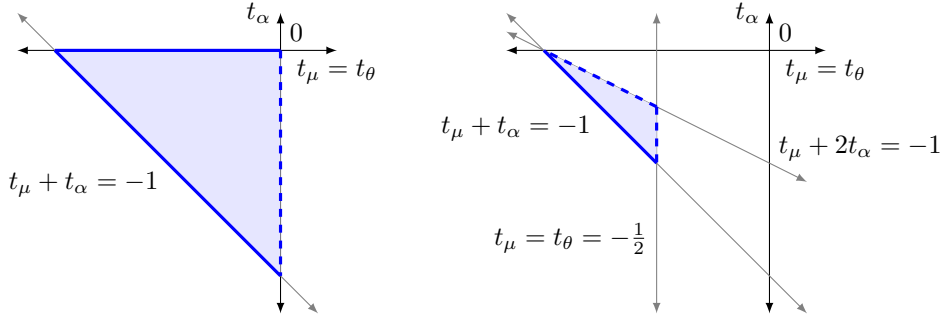
\begin{figure}[ht]
  \centering
  \begin{tikzpicture}

  \coordinate [label=above right:{$0$}] (o) at (0.0,0.0);
  \coordinate [label=below:{$t_\mu = t_\theta$}] (t1) at ( 0.75, 0.0);
  \coordinate [label=left:$t_\alpha$] (t2) at ( 0.0, 0.5);
  \coordinate [label=below left:{$t_\mu + t_\alpha = -1$}] (c1) at (-1.5,-1.5);
  \coordinate (e1) at (-3.5, 0.0);
  \coordinate (e2) at ( 0.0,-3.5);
  \draw[black,latex-latex] (t1) -- (e1);
  \draw[black,latex-latex] (t2) -- (e2);
  \draw[gray,latex-latex] (-3.5,0.5) -- (0.5,-3.5);
  \fill[blue!10] ( 0.0,0.0) -- (-3.0, 0.0) -- ( 0.0,-3.0);
  \draw[blue,very thick] ( 0.0,0.0) -- (-3.0, 0.0);
  \draw[blue,very thick,dashed] ( 0.0,0.0) -- ( 0.0,-3.0);
  \draw[blue,very thick] (-3.0,0.0) -- ( 0.0,-3.0);

  \begin{scope}[shift={(6.5,0.0)}]

  \coordinate [label=above right:{$0$}] (o) at (0.0,0.0);
  \coordinate [label=below:{$t_\mu = t_\theta$}] (t1) at ( 0.75, 0.0);
  \coordinate [label=left:$t_\alpha$] (t2) at ( 0.0, 0.5);
  \coordinate [label=below left:{$t_\mu + t_\alpha = -1$}] (c1) at (-2.25,-0.75);
  \coordinate [label=left:{$t_\mu = t_\theta = -\thalf$}] (c2) at (-1.5,-2.5);
  \coordinate [label=right:{$t_\mu + 2t_\alpha = -1$}] (c3) at (0.0,-1.3);
  \coordinate (e1) at (-3.5, 0.0);
  \coordinate (e2) at ( 0.0,-3.5);
  \draw[black,latex-latex] (t1) -- (e1);
  \draw[black,latex-latex] (t2) -- (e2);
  \draw[gray,latex-latex] (-3.5,0.5) -- (0.5,-3.5);
  \draw[gray,latex-latex] (-1.5,0.5) -- (-1.5,-3.5);
  \draw[gray,latex-latex] (-3.5,0.25) -- (0.5,-1.75);
  \fill[blue!10] (-1.50,-0.75) -- (-3.00, 0.00) -- (-1.50,-1.50);
  \draw[blue,very thick,dashed] (-1.50,-0.75) -- (-3.00, 0.00);
  \draw[blue,very thick,dashed] (-1.50,-0.75) -- (-1.50,-1.50);
  \draw[blue,very thick] (-1.50,-1.50) -- (-3.00,-0.00);

  \end{scope}

\end{tikzpicture}
  \caption{\edit{Allowable values for $(t_\mu,t_\theta,t_\alpha)$ for deterministic (left) and stochastic (right) settings.}}
  \label{fig.t_values}
\end{figure}

\edit{We illustrate in Figure~\ref{fig.t_values} the allowable values of $(t_\mu,t_\theta,t_\alpha)$ for the deterministic and stochastic settings.  Observe in particular that in each setting the set of allowable values is nonempty; e.g., in the stochastic setting, one can choose $t_\mu = t_\theta = -\tfrac34$ and $t_\alpha = -\tfrac14$.  Note, however, that the choice $t_\mu = t_\theta = -1$ is not allowed in the stochastic setting, even though it is allowed in the deterministic setting.}

\edit{
\bremark
  As in the deterministic setting, our main result for the stochastic setting (Theorem~\ref{th.main} on page~\pageref{th.main}) employs the choices $\mu_k = \mu_1 k^{t_\mu}$ and $\theta_{k-1} = \theta_0 k^{t_\theta}$ for all $k \in \N{}$, where $t_\mu = t_\theta$.  Our subsequent analysis shows that Parameter Rule~\ref{pm.step_size_gamma} ensures $\alpha_{k,\min} = \Omega(k^{t_\mu + t_\alpha})$, which along with $t_\mu + t_\alpha \in [-1,0)$ ensures that $\{\alpha_{k,\min}\}$ is unsummable.  Moreover, our subsequent analysis also shows that $\{\gamma_{k,\min} \alpha_{k,\min}\}$ is unsummable, which is critical since we show that the expected decrease in $\barrier$ in each iteration $k \in \N{}$ is bounded above by $-\tfrac14 \gamma_{k,\min} \alpha_{k,\min} \|\nabla_x \barrier(X_k,\mu_k)\|_{H_k^{-1}}^2$ plus a term due to the noise in the stochastic gradient estimates; see upcoming Lemma~\ref{lem.decrease_stochastic}.  At the same time, our analysis shows that, under Parameter Rule~\ref{pm.step_size_gamma}, this noise can be bounded as $\Ocal(k^{2t_\mu}) + \Ocal(k^{t_\mu + 2t_\alpha})$, which along with $t_\mu \in (-1,-\thalf)$ and $t_\mu + 2t_\alpha \in (-\infty,-1)$ ensures that this noise term is summable.  (It is to ensure this fact about the noise that the choice $t_\mu = t_\theta = -1$ and $t_\alpha = 0$---which would give $t_\mu + 2t_\alpha = -1$---is not allowed, even though it is in the deterministic setting.)  It follows overall that the expected decrease sufficiently dominates the noise so that our analysis can conclude that a subsequence of $\{\|\nabla_x \barrier(X_k,\mu_k)\|_{H_k^{-1}}\}$ is driven to zero almost surely.
\eremark
}

Since the barrier parameter sequence $\{\mu_k\}$ and neighborhood parameter sequence $\{\theta_k\}$ are prescribed and Lemmas~\ref{lem.derivatives}, \ref{lem.Lipschitz}, and \ref{lem.quadratic_upper} hold independently of any algorithm, it follows that the result of Lemma~\ref{lem.Lipschitz_bound} holds in the present setting. We formalize this fact in the following lemma, the proof of which is omitted since it would follow the same line of argument as the proof of Lemma~\ref{lem.Lipschitz_bound} stated previously.

\blemma\label{lem.Lipschitz_bound_stochastic}
  For all $k \in \N{}$, with $\ell_{\nabla f,\Bcal,\mu_k,X_k,X_{k+1}}$ defined as in Lemma~\ref{lem.Lipschitz} and $\ell_{\nabla f,\Bcal,k}$ defined as in Parameter Rule~\ref{pm.step_size}, one finds that
  \bequationNN
    \ell_{\nabla f, \Bcal, \mu_k, X_k, X_{k+1}} \leq \ell_{\nabla f,\Bcal,k} \leq \ell_{\nabla f,\Bcal,\mu_k,\theta_{k-1},\theta_k} \leq \ell_{\nabla f,\Bcal} + 2\mu_k\theta_k^{-2},
  \eequationNN
  from which it follows that Parameter Rule~\ref{pm.step_size_gamma} guarantees $A_k \in [\alpha_{k,\min},\alpha_{k,\max}]$.
\elemma

We also have that the following result, similar to Lemma~\ref{lem.gammakbound}, holds in the present setting.  The proof is omitted since it would follow the same line of argument as the proof of Lemma~\ref{lem.gammakbound}, the primary differences being that, for all $k \in \N{}$, one has $A_k \leq \alpha_{k,\max}$ and in place of $|\nabla_i f(x_k)| \leq \kappa_{\nabla f, \Bcal, \infty}$ one can employ $|G_{k,i}| \leq \kappa_{\nabla f,\Bcal,\infty} + \sigma_\infty$.

\blemma\label{lem.gammakbound_stochastic}
  For all $k \in \N{}$, Parameter Rule~\ref{pm.step_size_gamma} guarantees $\Gamma_k \in [\gamma_{k,\min},\gamma_{k,\max}]$.
\elemma

Our next lemma provides a preliminary upper bound on the expected per-iteration change in the shifted barrier-augmented function.

\blemma\label{lem.decrease_stochastic_prelim}
  For all $k \in \N{}$, one finds that
  \begin{align*}
    &\ \barrier(X_{k+1},\mu_{k+1}) - \barrier(X_k,\mu_k) \\
    \leq&\ -\Gamma_k A_k \|\nabla_x \barrier(X_k,\mu_k)\|_{H_k^{-1}}^2 + \Gamma_k A_k \nabla_x \barrier(X_k,\mu_k)^T H_k^{-1} (\nabla_x \barrier(X_k,\mu_k) - Q_k) \\
    &\ + \thalf \Gamma_k^2 A_k^2 \lambda_{k,\min}^{-1} \ell_{\nabla f,\Bcal,k} \|Q_k\|_{H_k^{-1}}^2.
  \end{align*}
\elemma
\bproof
  Similarly as in the proof of Lemma~\ref{lem.decrease}, for all $k \in \N{}$, one finds from Lemmas~\ref{lem.quadratic_upper} and \ref{lem.Lipschitz_bound_stochastic}, line~\ref{line:x-update} of Algorithm~\ref{alg.ipm}, and line~\ref{line:Hk} of Algorithm~\ref{alg.ipm} that
  \begin{align*}
    &\ \barrier(X_{k+1},\mu_k) - \barrier(X_k,\mu_k) \\
    \leq&\ \nabla_x \barrier(X_k,\mu_k)^T(X_{k+1} - X_k) + \thalf \ell_{\nabla f,\Bcal,\mu_k,X_k,X_{k+1}} \|X_{k+1} - X_k\|_2^2 \\
    \leq&\ -\nabla_x \barrier(X_k,\mu_k)^T(\Gamma_k A_k H_k^{-1} Q_k) + \thalf \ell_{\nabla f,\Bcal,k} \|\Gamma_k A_k H_k^{-1} Q_k \|_2^2 \\
    \leq&\ -\Gamma_k A_k \nabla_x \barrier(X_k,\mu_k)^T H_k^{-1} Q_k + \thalf \Gamma_k^2 A_k^2 \lambda_{k,\min}^{-1} \ell_{\nabla f,\Bcal,k} \|Q_k\|_{H_k^{-1}}^2.
  \end{align*}
  Adding and subtracting $-\Gamma_k A_k \|\nabla_x \barrier(X_k,\mu_k)\|_{H_k^{-1}}^2$ on the right-hand side and using the fact that Lemma~\ref{lem.derivatives} and $\mu_{k+1} < \mu_k$ imply that $\barrier(X_{k+1},\mu_{k+1}) < \barrier(X_{k+1},\mu_k)$, one reaches the desired conclusion.
\eproof

Our next lemma provides an upper bound on the conditional expectation of the middle term on the right-hand side of the inequality in Lemma~\ref{lem.decrease_stochastic_prelim}.

\blemma\label{eq.inner_product_bound}
  For all $k \in \N{}$, one finds that
  \begin{align*}
    &\ \E[\Gamma_k A_k \nabla_x \barrier(X_k,\mu_k)^T H_k^{-1} (\nabla_x \barrier(X_k,\mu_k) - Q_k) | \Fcal_k] \\
    \leq&\ (\gamma_{k,\min} \alpha_{k,\buff} + \gamma_{k,\buff} \alpha_{k,\min} + \gamma_{k,\buff} \alpha_{k,\buff}) \lambda_{k,\min}^{-1} (\kappa_{\nabla f,\Bcal,2} + 2 \sqrt{n} \mu_k \theta_{k-1}^{-1}) \sigma_2.
  \end{align*}
\elemma
\bproof
  Let $\Ical_k$ be the event that $P_k := \nabla_x \barrier(X_k,\mu_k)^T H_k^{-1} (\nabla_x \barrier(X_k,\mu_k) - Q_k) \geq 0$ and let~$\Ical_k^c$ be the complementary event that $P_k < 0$.  By Assumption~\ref{ass.stochastic}, Parameter Rule~\ref{pm.step_size_gamma}, the Law of Total Expectation, $\E[P_k | \Fcal_k] = 0$, and Lemmas~\ref{lem.Lipschitz_bound_stochastic} and \ref{lem.gammakbound_stochastic},
  \begin{align*}
    &\ \E[\Gamma_k A_k P_k | \Fcal_k] \\
    =&\ \E[\Gamma_k A_k P_k | \Fcal_k \wedge \Ical_k] \P[\Ical_k | \Fcal_k] + \E[\Gamma_k A_k P_k | \Fcal_k \wedge \Ical_k^c] \P[\Ical_k^c | \Fcal_k] \\
    \leq&\ \gamma_{k,\max} \alpha_{k,\max} \E[P_k | \Fcal_k \wedge \Ical_k] \P[\Ical_k | \Fcal_k] + \gamma_{k,\min} \alpha_{k,\min} \E[P_k | \Fcal_k \wedge \Ical_k^c] \P[\Ical_k^c | \Fcal_k] \\
    \leq&\ \gamma_{k,\min} \alpha_{k,\min} (\E[P_k | \Fcal_k \wedge \Ical_k] \P[\Ical_k | \Fcal_k] + \E[P_k | \Fcal_k \wedge \Ical_k^c] \P[\Ical_k^c | \Fcal_k]) \\
    &\ + (\gamma_{k,\min} \alpha_{k,\buff} + \gamma_{k,\buff} \alpha_{k,\min} + \gamma_{k,\buff} \alpha_{k,\buff}) \E[P_k | \Fcal_k \wedge \Ical_k] \P[\Ical_k | \Fcal_k] \\
    =&\ (\gamma_{k,\min} \alpha_{k,\buff} + \gamma_{k,\buff} \alpha_{k,\min} + \gamma_{k,\buff} \alpha_{k,\buff}) \E[P_k | \Fcal_k \wedge \Ical_k] \P[\Ical_k | \Fcal_k].
  \end{align*}
  By the Cauchy-Schwarz inequality and Assumptions~\ref{ass.f} and~\ref{ass.stochastic}, one has
  \begin{align*}
    &\ \E[P_k | \Fcal_k \wedge \Ical_k] \P[\Ical_k | \Fcal_k] \\
    \leq&\ \E[ \|H_k^{-1} \nabla_x \barrier(X_k,\mu_k)\|_2 \|\nabla_x \barrier(X_k,\mu_k) - Q_k\|_2 | \Fcal_k \wedge \Ical_k] \P[\Ical_k | \Fcal_k] \\
    =&\ \E[ \|H_k^{-1} \nabla_x \barrier(X_k,\mu_k)\|_2 \|G_k - \nabla f(X_k)\|_2 | \Fcal_k \wedge \Ical_k] \P[\Ical_k | \Fcal_k] \\
    \leq&\ \lambda_{k,\min}^{-1} (\kappa_{\nabla f,\Bcal,2} + 2 \sqrt{n} \mu_k \theta_{k-1}^{-1}) \sigma_2,
  \end{align*}
  which combined with the result above yields the desired conclusion.
\eproof

Our next lemma provides an upper bound on the last term on the right-hand side of the inequality in Lemma~\ref{lem.decrease_stochastic_prelim}.

\blemma\label{lem.noise_bound}
  For all $k \in \N{}$, one finds that
  \begin{align*}
    &\ \thalf \Gamma_k^2 A_k^2 \lambda_{k,\min}^{-1} \ell_{\nabla f,\Bcal,k} \|Q_k\|_{H_k^{-1}}^2 \\
    \leq&\ \tfrac{3}{4} \Gamma_k^2 A_k^2 \lambda_{k,\min}^{-1} \ell_{\nabla f,\Bcal,k} \|\nabla_x \barrier(X_k,\mu_k)\|_{H_k^{-1}}^2 + \edit{\tfrac{3}{2} ( \tfrac{k^{2t_\alpha}}{\ell_{\nabla f, \Bcal} + 2\mu_k \theta_k^{-2}} + \alpha_{\buff} \lambda_{k,\min}^{-1} k^{2t_\mu+t_\alpha} ) \sigma_2^2}.
  \end{align*}
\elemma
\bproof
  Consider arbitrary $k \in \N{}$.  Since for any $(a,b) \in \R{n} \times \R{n}$, the fact that $\|\thalf a - b\|_{H_k^{-1}}^2 \geq 0$ implies that $\|a + b\|_{H_k^{-1}}^2 \leq \tfrac{3}{2}\|a\|_{H_k^{-1}}^2 + 3\|b\|_{H_k^{-1}}^2$, it follows that
  \begin{align*}
    \thalf \|Q_k\|_{H_k^{-1}}^2
    =&\ \thalf \|\nabla_x \barrier(X_k,\mu_k) + Q_k - \nabla_x \barrier(X_k,\mu_k)\|_{H_k^{-1}}^2 \\
    \leq&\ \tfrac{3}{4}\|\nabla_x \barrier(X_k,\mu_k)\|_{H_k^{-1}}^2 + \tfrac{3}{2}\|Q_k - \nabla_x \barrier(X_k,\mu_k)\|_{H_k^{-1}}^2 \\
    \leq&\ \tfrac{3}{4} \|\nabla_x \barrier(X_k,\mu_k)\|_{H_k^{-1}}^2 + \tfrac{3}{2} \lambda_{k,\min}^{-1} \sigma_2^2.
  \end{align*}
  \edit{Thus, under Parameter Rule~\ref{pm.step_size_gamma}, specifically the facts that $\Gamma_k \leq 1$, $A_k \leq \tfrac{\lambda_{k,\min}k^{t_\alpha}}{\ell_{\nabla f, \Bcal, k}}$, and $A_k \leq \alpha_{k,\max} = \alpha_{k,\min} + \alpha_{k,\buff} \leq \tfrac{\lambda_{k,\min}k^{t_\alpha}}{\ell_{\nabla f, \Bcal} + 2 \mu_k \theta_k^{-2}} + \alpha_{\buff} k^{2t_\mu}$, one finds that
  \begin{align*}
    &\ \thalf \Gamma_k^2 A_k^2 \lambda_{k,\min}^{-1} \ell_{\nabla f,\Bcal,k} \|Q_k\|_{H_k^{-1}}^2 \\
    \leq&\ \tfrac{3}{4}\Gamma_k^2 A_k^2 \lambda_{k,\min}^{-1} \ell_{\nabla f,\Bcal,k} \|\nabla_x \barrier(X_k,\mu_k)\|_{H_k^{-1}}^2 + \tfrac{3}{2}\Gamma_k^2 A_k^2 \lambda_{k,\min}^{-2}\ell_{\nabla f, \Bcal, k} \sigma_2^2 \\
    \leq&\ \tfrac{3}{4}\Gamma_k^2 A_k^2 \lambda_{k,\min}^{-1} \ell_{\nabla f,\Bcal,k} \|\nabla_x \barrier(X_k,\mu_k)\|_{H_k^{-1}}^2 \\
    &\qquad + \tfrac{3}{2} \tfrac{\lambda_{k,\min}k^{t_\alpha}}{\ell_{\nabla f, \Bcal, k}} ( \tfrac{\lambda_{k,\min}k^{t_\alpha}}{\ell_{\nabla f, \Bcal} + 2 \mu_k \theta_k^{-2}} + \alpha_{\buff} k^{2t_\mu} ) \lambda_{k,\min}^{-2} \ell_{\nabla f, \Bcal, k} \sigma_2^2 \\
    =&\ \tfrac{3}{4} \Gamma_k^2 A_k^2 \lambda_{k,\min}^{-1} \ell_{\nabla f,\Bcal,k} \|\nabla_x \barrier(X_k,\mu_k)\|_{H_k^{-1}}^2 + \tfrac{3}{2} ( \tfrac{k^{2t_\alpha}}{\ell_{\nabla f, \Bcal} + 2\mu_k \theta_k^{-2}} + \alpha_{\buff} \lambda_{k,\min}^{-1} k^{2t_\mu+t_\alpha} ) \sigma_2^2,
  \end{align*}
  as desired.}
\eproof

Combining the prior three lemmas, we obtain the following result.  This result is reminiscent of Lemma~\ref{lem.decrease}, but accounts for the stochastic gradient errors.

\blemma\label{lem.decrease_stochastic}
  For all $k \in \N{}$, one finds that
  \begin{align*}
    &\ \E[\barrier(X_{k+1},\mu_{k+1}) | \Fcal_k] - \barrier(X_k,\mu_k) \\
    \leq&\ - \tfrac{1}{4} \gamma_{k,\min} \alpha_{k,\min} \|\nabla_x \barrier(X_k,\mu_k)\|_{H_k^{-1}}^2 \\
    &\ + (\gamma_{k,\min} \alpha_{k,\buff} + \gamma_{k,\buff} \alpha_{k,\min} + \gamma_{k,\buff} \alpha_{k,\buff}) \lambda_{k,\min}^{-1} (\kappa_{\nabla f,\Bcal,2} + 2 \sqrt{n} \mu_k \theta_{k-1}^{-1}) \sigma_2 \\
    &\ + \edit{\tfrac{3}{2} ( \tfrac{k^{2t_\alpha}}{\ell_{\nabla f, \Bcal} + 2\mu_k \theta_k^{-2}} + \alpha_{\buff} \lambda_{k,\min}^{-1} k^{2t_\mu+t_\alpha} ) \sigma_2^2}.
  \end{align*}
\elemma
\bproof
  Consider arbirary $k \in \N{}$.  Combining Lemmas~\ref{lem.decrease_stochastic_prelim}, \ref{eq.inner_product_bound}, and \ref{lem.noise_bound},
  \begin{align*}
    &\ \E[\barrier(X_{k+1},\mu_{k+1}) | \Fcal_k] - \barrier(X_k,\mu_k) \\
    \leq&\ - \E[\Gamma_k A_k (1 - \tfrac{3}{4} \Gamma_k A_k \lambda_{k,\min}^{-1} \ell_{\nabla f,\Bcal,k}) \|\nabla_x \barrier(X_k,\mu_k)\|_{H_k^{-1}}^2 | \Fcal_k] \\
    &\ + (\gamma_{k,\min} \alpha_{k,\buff} + \gamma_{k,\buff} \alpha_{k,\min} + \gamma_{k,\buff} \alpha_{k,\buff}) \lambda_{k,\min}^{-1} (\kappa_{\nabla f,\Bcal,2} + 2 \sqrt{n} \mu_k \theta_{k-1}^{-1}) \sigma_2 \\
    &\ +  \edit{\tfrac{3}{2} ( \tfrac{k^{2t_\alpha}}{\ell_{\nabla f, \Bcal} + 2\mu_k \theta_k^{-2}} + \alpha_{\buff} \lambda_{k,\min}^{-1} k^{2t_\mu+t_\alpha} ) \sigma_2^2}.
  \end{align*}
  Now, one finds under Parameter Rule~\ref{pm.step_size_gamma} \edit{(so $\Gamma_k \leq 1$ and $k^{t_\alpha} \leq 1$)} that
  \begin{align*}
    A_k \leq \tfrac{\lambda_{k,\min}\edit{k^{t_\alpha}}}{\ell_{\nabla f, \Bcal, k}} \implies 1 - \tfrac{3}{4} \Gamma_kA_k\lambda_{k,\min}^{-1} \ell_{\nabla f,\Bcal,k} \geq 1 - \tfrac{3}{4}\Gamma_k\edit{k^{t_\alpha}} \ge \tfrac{1}{4}.
  \end{align*}  
  Thus, from above, Parameter Rule~\ref{pm.step_size_gamma}, and Lemma~\ref{lem.Lipschitz_bound_stochastic}, the conclusion follows.
\eproof

We now show that if the parameter sequences are chosen similarly as in Corollary~\ref{co.sequence_choice}, then the coefficients in the upper bound proved in Lemma~\ref{lem.decrease_stochastic} satisfy critical properties for proving our ultimate convergence guarantee.  \edit{We recall that in our proof of Corollary~\ref{co.sequence_choice}, we conducted extra steps in order to show why it is needed to have $t_\mu = t_\theta$.  This restriction is again required in the proof of the following lemma, and for conciseness in certain situations we immediately employ this equation in the proof.}

\blemma\label{lem.sequences_stochastic}
  If \edit{for $(t_\mu,t_\theta,t_\alpha) \in (-1,-\thalf) \times (-1,-\thalf) \times (-\infty,0)$ in Parameter Rule~\ref{pm.step_size_gamma} such that $t_\mu = t_\theta$, $t_\mu + t_\alpha \in [-1,0)$, and $t_\mu + 2t_\alpha \in (-\infty,-1)$} and for some $\underline{r} \in \R{}_{>0}$ and $\mu_1 \in \R{}_{>0}$ with $\mu_1 > \tfrac{\thalf \theta_0 (\kappa_{\nabla f,\Bcal,\infty} + \sigma_\infty) \Delta}{\thalf \Delta - \theta_0}$ the algorithm has $\mu_k = \mu_1 k^{\edit{t_\mu}}$, $\theta_{k-1} = \theta_0 k^{\edit{t_\theta}}$, and $\underline{r} \leq \lambda_{k,\min} \leq \lambda_{k,\max}$ for all $k \in \N{}$, then there \edit{exists} $(\ktilde,c,C) \in \N{} \times \R{}_{>0} \times \R{}_{>0}$ with
  \bequationNN
    \E[\barrier(X_{k+1},\mu_{k+1}) | \Fcal_k] - \barrier(X_k,\mu_k) \leq - \edit{c k^{t_\mu+t_\alpha}} \|\nabla_x \barrier(X_k,\mu_k)\|_{H_k^{-1}}^2 + C k^{\edit{\max\{2t_\mu,t_\mu+2t_\alpha\}}}
  \eequationNN
  for all $k \in \N{}$ with $k \geq \ktilde$.
\elemma
\bproof
  \edit{We prove the result using Lemma~\ref{lem.decrease_stochastic}; we show that there exists $c \in \R{}_{>0}$ such that $\tfrac{1}{4} \gamma_{k,\min} \alpha_{k,\min} \geq c k^{t_\mu+t_\alpha}$ for all sufficiently large $k \in \N{}$, then show that there exists $C \in \R{}_{>0}$ such that the latter two terms on the right-hand side of the upper bound in Lemma~\ref{lem.decrease_stochastic} are bounded by $C k^{\max\{2t_\mu, t_\mu + 2t_\alpha\}}$ for all sufficiently large $k \in \N{}$.} Under the conditions and Parameter Rule~\ref{pm.step_size_gamma}, one has for all $k \in \N{}$ that
  \begin{align*}
    \gamma_{k,\min} \alpha_{k,\min}
    &\edit{=} \min\left\{ \alpha_{k,\min} , \tfrac{\alpha_{k,\min} \lambda_{k,\min} \(\tfrac{\thalf \mu_k \Delta}{\mu_k + \thalf (\kappa_{\nabla f,\Bcal,\infty} + \sigma_\infty) \Delta} - \theta_k\)}{(\alpha_{k,\min} + \alpha_{k,\buff}) (\kappa_{\nabla f,\Bcal,\infty} + \sigma_\infty + \mu_k \theta_{k-1}^{-1})} \right\} \\
    &=: \min\{\alpha_{k,\min}, \hat\eta_k\}.
  \end{align*}
  The proof can now proceed similarly as in the proof of Corollary~\ref{co.sequence_choice}.  In particular, with the given parameter choices and by \eqref{betak-bound}, one finds for all $k \in \N{}$ that
  \bequation\label{eq.stoch_1}
    \baligned
      \alpha_{k,\min} &\geq \tfrac{\underline{r} \edit{k^{t_\alpha}}}{\ell_{\nabla f,\Bcal} + 2 \mu_1 \theta_0^{-2} k^{\edit{t_\mu}} (k+1)^{-2\edit{t_\theta}}} \\
      &\geq \tfrac{\underline{r}\edit{k^{t_\alpha}}}{\ell_{\nabla f,\Bcal} + 2 \mu_1 \theta_0^{-2} 2^{-2\edit{t_\mu}} k^{\edit{t_\mu - 2t_\theta}}} \edit{\geq \chat k^{-t_\mu + 2t_\theta + t_\alpha} = \chat k^{t_\mu + t_\alpha}}
    \ealigned
  \eequation
  \edit{for some $\chat \in \R{}_{>0}$. Now, for the sequence $\{\hat{\eta}_k\}$, first observe that under the conditions of the lemma one has $2t_\mu < t_\mu + t_\alpha < 0$ since $t_\mu < -\thalf$ and $-1 \leq t_\mu + t_\alpha < 0$.  Thus, by $\{\alpha_{k,\buff}\} = \Ocal(k^{2\edit{t_\mu}})$ and \eqref{eq.stoch_1},} it follows that there exists $\khat \in \N{}$ such that $\alpha_{k,\buff} \leq \alpha_{k,\min}$ for all $k \in \N{}$ with $k \geq \khat$, which in turn means that
  \bequation\label{eq.stoch_2}
    \tfrac{\alpha_{k,\min}}{\alpha_{k,\min} + \alpha_{k,\buff}} \geq \thalf\ \ \text{for all}\ \ k \in \N{}\ \ \text{with}\ \ k \geq \khat.
  \eequation
  \edit{One finds with \eqref{eq.theta0}, $t_\mu \leq 0$, $t_\theta \leq 0$, and $t_\mu \geq t_\theta$ that a similar derivation as in the proof of Corollary~\ref{co.sequence_choice} shows that there exists $\ctilde \in \R{}_{>0}$ such that, for all $k \in \N{}$,}
  \bequation\label{eq.stoch_3}
    \tfrac{\thalf \mu_k \Delta}{\mu_k + \thalf (\kappa_{\nabla f,\Bcal,\infty} + \sigma_\infty) \Delta} - \theta_k \geq \edit{\ctilde k^{t_\theta} = \ctilde k^{t_\mu}}.
  \eequation
  Combining \eqref{eq.stoch_1}--\eqref{eq.stoch_3} \edit{and $t_\alpha \leq 0$}, there exists $c \in \R{}_{>0}$ such that $\tfrac{1}{4}\gamma_{k,\min} \alpha_{k,\min} \geq \edit{ck^{t_\mu+t_\alpha}} $ for all $k \in \N{}$ with $k \geq \khat$.  On the other hand, the conditions of the lemma and Parameter Rule~\ref{pm.step_size_gamma} imply \edit{that there exists $\cbar \in \R{}_{>0}$ such that, for all $k \in \N{}$,}
  \bequation\label{eq.alpha_rate_2}
    \edit{\alpha_{k,\min} \lambda_{k,\min}^{-1} = \tfrac{k^{t_\alpha}}{\ell_{\nabla f,\Bcal} + 2 \mu_1 \theta_0^{-2} k^{t_\mu} (k+1)^{-2t_\theta}} \leq \tfrac{1}{\ell_{\nabla f,\Bcal} + 2 \mu_1 \theta_0^{-2} k^{t_\mu} k^{-2t_\theta}} \leq \cbar k^{t_\mu}.}
  \eequation
  \edit{Combining this fact, Parameter Rule~\ref{pm.step_size_gamma}, and $\lambda_{k,\min}^{-1} \leq \underline{r}^{-1}$ due to the conditions of the lemma, one finds that $\{\gamma_{k,\min}\alpha_{k,\buff}\lambda_{k,\min}^{-1}\} = \Ocal(k^{2t_\mu})$, $\{\gamma_{k,\buff}\alpha_{k,\min}\lambda_{k,\min}^{-1}\} = \Ocal(k^{2t_\mu})$, and $\{\gamma_{k,\buff}\alpha_{k,\buff}\lambda_{k,\min}^{-1}\} = o(k^{2t_\mu})$.  Combining these facts and $\mu_k \theta_{k-1}^{-1} = \mu_1 \theta_0^{-1} k^{t_\mu - t_\theta} = \mu_1 \theta_0^{-1}$ for all $k \in \N{}$, it follows that the second term on the right-hand side of the upper bound in Lemma~\ref{lem.decrease_stochastic} is $\Ocal(k^{2t_\mu})$.  Hence, all that remains is to consider the final term in that upper bound.  Toward this end, observe that by Parameter Rule~\ref{pm.step_size_gamma} there exists $\overline{\overline{c}} \in \R{}_{>0}$ such that, for all $k \in \N{}$, one finds that
  \bequationNN
    \tfrac{k^{2t_\alpha}}{\ell_{\nabla f,\Bcal} + 2 \mu_k \theta_k^{-2}} = \tfrac{k^{2t_\alpha}}{\ell_{\nabla f,\Bcal} + 2 \mu_1 \theta_0^{-2} k^{t_\mu} (k+1)^{-2t_\theta}} \le \tfrac{k^{2t_\alpha}}{\ell_{\nabla f,\Bcal} + 2 \mu_1 \theta_0^{-2} k^{t_\mu} k^{-2t_\theta}} \le \overline{\overline{c}} k^{t_\mu + 2t_\alpha},
  \eequationNN
  and hence, with $\lambda_{k,\min}^{-1} \leq \underline{r}^{-1}$ due to the conditions of the lemma, one has
  \bequationNN
    \tfrac{3}{2} ( \tfrac{k^{2t_\alpha}}{\ell_{\nabla f, \Bcal} + 2\mu_k \theta_k^{-2}} + \alpha_{\buff} \lambda_{k,\min}^{-1} k^{2t_\mu+t_\alpha} ) \sigma_2^2 = \Ocal(k^{t_\mu + 2t_\alpha}) + \Ocal(k^{2t_\mu+t_\alpha}).
  \eequationNN
  Combined with the bound proved previously, it follows that the sum of the last two terms on the right-hand side of the upper bound in Lemma~\ref{lem.decrease_stochastic} is $\Ocal(k^{2t_\mu}) + \Ocal(k^{t_\mu + 2t_\alpha})$ (where $k^{2t_\mu+t_\alpha} = \Ocal(k^{2t_\mu})$ since $2t_\mu + t_\alpha < 2t_\mu$), so the desired conclusion follows.}
\eproof

We now prove our main convergence theorem for the stochastic setting.

\btheorem\label{th.main}
  Suppose that Assumptions~\ref{ass.f} and \ref{ass.bounded} and Parameter Rule~\ref{pm.step_size_gamma} hold, and that the parameter sequences are chosen as in Lemma~\ref{lem.sequences_stochastic}.  Then,
  \bequationNN
    \liminf_{k \to \infty} \|\nabla_x \barrier(X_k,\mu_k)\|_{H_k^{-1}}^2 = 0\ \ \text{almost surely},
  \eequationNN
  meaning that if there exists $\overline{r} \in \R{}_{>0}$ such that $\lambda_{k,\max} \leq \overline{r}$ for all $k \in \N{}$, then
  \bequationNN
    \liminf_{k\to\infty} \|\nabla_x \barrier(X_k,\mu_k)\|_2^2 = 0\ \ \text{almost surely}.
  \eequationNN
  Consequently, if all of the aforementioned choices of the parameter sequences are made and, in a given run of the algorithm generating a realization of the iterate sequence $\{x_k\}$ there exists an infinite-cardinality set $\Kcal \subseteq \N{}$ such that $\{\nabla \barrier(x_k,\mu_k)\}_{k\in\Kcal} \to 0$ and $\{x_k\}_{k\in\Kcal} \to \xbar$ for some $\xbar \in \Bcal$, then the limit point $\xbar$ is a KKT point for \eqref{prob.opt}.
\etheorem
\bproof
  By the Law of Total Expectation, it follows from Lemma~\ref{lem.sequences_stochastic} that there exists $(\ktilde,c,C) \in \R{}_{>0} \times \R{}_{>0}$ such that, for all $k \in \N{}$ with $k \geq \ktilde$, one has
  \begin{multline*}
    \E[\barrier(X_{k+1},\mu_{k+1})] - \E[\barrier(X_k,\mu_k)] \\ \leq - \edit{c k^{t_\mu + t_\alpha}} \E[\|\nabla_x \barrier(X_k,\mu_k)\|_{H_k^{-1}}^2] + C k^{\edit{\max\{2t_\mu,t_\mu+2t_\alpha\}}}.
  \end{multline*}
  \edit{Now consider arbitrary $K \in \N{}$. Summing the expression above} for $k \in \{\ktilde,\dots,\ktilde + K\}$, it follows along with Lemma~\ref{lem.derivatives} that
  \begin{align*}
    f_{\inf} - \E[\barrier(x_{\ktilde},\mu_{\ktilde})]
      &\leq \E[\barrier(X_{k+1},\mu_{k+1})] - \E[\barrier(x_{\ktilde},\mu_{\ktilde})] \\
      &\leq - c \sum_{k=\ktilde}^{\ktilde + K} k^{\edit{t_\mu + t_\alpha}} \E[\|\nabla_x \barrier(X_k,\mu_k)\|_{H_k^{-1}}^2] + C \sum_{k=\ktilde}^{\ktilde+K} k^{\edit{\max\{2t_\mu,t_\mu+2t_\alpha\}}},
  \end{align*}
  which after rearrangement yields
  \bequationNN
    \sum_{k=\ktilde}^{\ktilde+K} k^{\edit{t_\mu+t_\alpha}} \E[\|\nabla_x \barrier(X_k,\mu_k)\|_{H_k^{-1}}^2] \leq \edit{\tfrac{1}{c}} (\E[\barrier(x_{\ktilde},\mu_{\ktilde})] - f_{\inf}) + \edit{\tfrac{C}{c}} \sum_{k=\ktilde}^{\ktilde+K} k^{\edit{\max\{2t_\mu,t_\mu+2t_\alpha\}}}.
  \eequationNN
  Under Assumptions~\ref{ass.f} and \ref{ass.bounded} and since \edit{$2t_\mu \in (-\infty,-1)$ and $t_\mu + 2t_\alpha \in (-\infty,-1)$}, the right-hand side of this inequality converges to a finite limit as $K \to \infty$.  Since $\sum_{k=1}^\infty k^{\edit{t_\mu + t_\alpha}} = \infty$, one finds along with the nonnegativity of $\|\nabla_x \barrier (X_k, \mu_k)\|_{H_{k}^{-1}}^2$ and Fatou's lemma that
  \bequationNN
    0 = \liminf_{k \to \infty} \E[\|\nabla_x \barrier (X_k, \mu_k)\|_{H_{k}^{-1}}^2] \geq \E[\liminf_{k \to \infty} \|\nabla_x \barrier (X_k, \mu_k)\|_{H_{k}^{-1}}^2] = 0.
  \eequationNN
  Consider the random variable $L := \liminf_{k \to \infty} \|\nabla_x \barrier (X_k, \mu_k)\|_{H_{k}^{-1}}^2$.  By nonnegativity of $\|\nabla_x \barrier (X_k,\mu_k)\|_{H_{k}^{-1}}^2$ and the Law of Total Expectation, it follows from above that $0 = \E[L] \geq \P[L > 0] \E[L | L > 0]$, so $0 = \P[L > 0] = \P[\liminf_{k \to \infty} \|\nabla_x \barrier (X_k, \mu_k)\|_{H_{k}^{-1}}^2 > 0]$, which is the first desired conclusion.  The second desired conclusion follow from the fact that if $\overline{r}$ exists as stated, then $\|\cdot\|_{H_k^{-1}}^2 \geq \overline{r}^{-1} \|\cdot\|_2^2$ for all $k \in \N{}$. The  last desired conclusion follows using the same argument as in the proof of Theorem~\ref{th.deterministic}.
\eproof

Results similar to Theorem~\ref{th.main} have been proved for the stochastic gradient method in the unconstrained setting \cite{BottCurtNoce18}.  In fact, for the stochastic gradient method in an unconstrained setting, one can prove under conditions that are similar to ours that the gradient of the objective function \edit{converges} to zero almost surely; this is stronger than a $\liminf$ result of the type in Theorem~\ref{th.main}.  However, such results rely on the gradient of the objective function being Lipschitz continuous.  In our setting, we have the Lipschitz-continuity-type property in Lemma~\ref{lem.Lipschitz}, but the gradient of the (shifted) barrier function is not globally Lipschitz over $\Bcal$, meaning that such a result in the unconstrained setting does not carry over to our present setting.

\section{Obstacles for a Simplified Algorithm}\label{sec.convex}

A reader may wonder if convergence guarantees can be proved for a simpler variant of our algorithm, namely, one that simply employs projections onto the inner neighborhoods of the feasible region.  After all, in the setting of convex optimization, convergence guarantees exist for \edit{projected-gradient-based} methods; see, e.g., \cite{NemiJudiLanShap09}.  In this section, we argue that the situation is not straightforward, and even though one may extend our algorithm and analysis to consider a less conservative strategy (recall Remark~\ref{rem.extend}), one runs into obstacles when trying to provide a convergence guarantee for a variant of our algorithm that simply uses a projection of $x_k + \alpha_k d_k$ for all $k \in \N{}$.

Let us consider the setting when $f$ is $\psi$-strongly convex for some $\psi \in \R{}_{>0}$, and let us consider a simplified variant of our algorithm, where, for all $k \in \N{}$,
\bequation\label{eq.convex_update}
  H_k = I,\ \alpha_k = \tfrac{1}{\ell_{\nabla f,\Bcal} + 2\mu_k\theta_k^{-2}},\ \text{and}\ x_{k+1} \gets \proj_{\hood(\theta_k)} (x_k - \alpha_k \nabla_x \phi(x_k,\mu_k)),
\eequation
where $\proj_{\Scal}(\cdot)$ denotes the orthogonal projection operator on convex $\Scal \subseteq \R{n}$.  For reasons seen in our subsequent analysis, suppose for all $k \in \N{}$ that
\bequation\label{eq.convex_step_size}
  \theta_k \leq c\mu_k,\ \ \text{where}\ \ c := \min_{i \in \Lcal \cup \Ucal} \left\{ \(\kappa_{\nabla f,\Bcal,\infty} + \tfrac{2 \mu_1}{u_i - l_i} \)^{-1} \right\}.
\eequation

Following a standard approach in the convex optimization literature, let us attempt to prove convergence of the algorithm defined by \eqref{eq.convex_update}--\eqref{eq.convex_step_size} by showing that the distance to the unique solution of \eqref{prob.opt}, call it $x_* \in \R{n}$, vanishes.  Since each iteration of the algorithm makes a step toward minimizing $\barrier(\cdot,\mu_k)$, it is natural to approach the analysis by considering the distance between the unique minimizer of this function to $x_*$.  A typical result in the literature is that, for sufficiently small barrier parameter values, this distance is proportional to $\mu_k$.  For concreteness, we state the following result; for further discussion and a proof, see \cite{wright2002properties}.

\bproposition\label{prop.distance_to_solution}
  Suppose that Assumptions~\ref{ass.f} and \ref{ass.bounded} hold and the objective~$f$ is $\psi$-strongly convex.  Let $x_* \in \R{n}$ be the unique point such that \eqref{eq.KKT} holds for some $(y_*,z_*) \in \R{n} \times \R{n}$ and let $\xbar_k \in \R{n}$ denote the unique minimizer of $\barrier(\cdot,\mu_k) : \R{n} \to \R{}$ for all $k \in \N{}$.  Then, for all sufficiently large $k \in \N{}$, it holds that $\|\xbar_k - (x_* + \mu_k \zeta)\|_2 = \Ocal(\mu_k^2)$, where the vector $\zeta \in \R{n}$ is defined independently from $\{\mu_k\}$.
\eproposition

Let us now show the expected result that with the step-size choice in \eqref{eq.convex_step_size}, the iterate update in \eqref{eq.convex_update} corresponds to a step toward the unique minimizer of $\barrier(\cdot,\mu_k)$.

\bproposition\label{prop.contraction_to_xbar}
  Suppose that Assumptions~\ref{ass.f} and \ref{ass.bounded} hold, the objective~$f$ is $\psi$-strongly convex, and the algorithm employs the update in \eqref{eq.convex_update}--\eqref{eq.convex_step_size}.  Let $\xbar_k \in \R{n}$ denote the unique minimizer of $\barrier(\cdot,\mu_k) : \R{n} \to \R{}$ for all $k \in \N{}$.  Then, for all $k \in \N{}$, one finds that $x_k \in \hood(\theta_{k-1}) \subset \hood(\theta_k)$, $x_{k+1} \in \hood(\theta_k)$, and $\xbar_k \in \hood\(c \mu_k\)$, where $c \in \R{}_{>0}$ is defined as in \eqref{eq.convex_step_size}.  Consequently, it follows for all $k \in \N{}$ that $(x_k,x_{k+1},\xbar_k) \in \hood(\theta_k) \times \hood(\theta_k) \times \hood(\theta_k)$ and
  \bequationNN
    \|x_{k+1} - \xbar_k\|_2^2 \leq (1 - \alpha_k \psi) \|x_k - \xbar_k\|_2^2.
  \eequationNN
\eproposition
\bproof
  That $x_k \in \hood(\theta_{k-1}) \subset \hood(\theta_k)$ and $x_{k+1} \in \hood(\theta_k)$ hold for all $k \in \N{}$ follows from \eqref{eq.convex_update} and the fact that $\{\theta_k\} \searrow 0$.  Now consider arbitrary $k \in \N{}$ and observe from the definition of $\xbar_k$ and Lemma~\ref{lem.derivatives} that it is the unique vector such that \eqref{eq.barrier_stationary} holds with $x = \xbar_k$ and $\mu = \mu_k$.  Let us prove the desired conclusion that $\xbar_k \in \hood(c \mu_k)$ by proving that $c\mu_k$ is a lower bound for $\sbar_{k,i} := \min\{\xbar_{k,i} - l_i, u_i - \xbar_{k,i}\}$ for all $i \in [n]$.  Consider arbitrary $i \in [n]$.  If $i \notin \Lcal \cup \Ucal$, then $\sbar_{k,i} = \infty \geq c\mu_k$, as desired.  If $i \in \Lcal \setminus \Ucal$, then \eqref{eq.barrier_stationary} and Assumption~\ref{ass.f} imply that $\xbar_{k,i} - l_i = \mu_k/\nabla_i f(\xbar_k) \geq \mu_k / \kappa_{\nabla f,\Bcal,\infty}$, which again yields that $\sbar_{k,i} \geq c\mu_k$.  Using a similar argument, the bound also holds for $i \in \Ucal \setminus \Lcal$.  Finally, if $i \in \Lcal \cap \Ucal$, then assuming without loss of generality that $\xbar_{k,i} \leq (u_i + l_i)/2$, it follows that $1/(u_i - \xbar_{k,i}) \leq 2/(u_i - l_i)$, so \eqref{eq.barrier_stationary} yields
  \bequationNN
    \tfrac{1}{\xbar_{k,i} - l_i} = \tfrac{\nabla_i f(\xbar_k)}{\mu_k} + \tfrac{1}{u_i - \xbar_{k,i}} \leq \tfrac{\nabla_i f(\xbar_k)}{\mu_k} + \tfrac{2}{u_i - l_i} = \tfrac{\nabla_i f(\xbar_k) (u_i - l_i) + 2\mu_k}{\mu_k(u_i - l_i)},
  \eequationNN
  which along with Assumption~\ref{ass.f} and $\sbar_{k,i} = \xbar_{k,i} - l_i$ yields the desired conclusion.  (If $\xbar_{k,i} \geq (u_i + l_i)/2$, the conclusion follows using a similar argument with $\sbar_{k,i} = u_i - \xbar_{k,i}$.)

  It has been shown that $(x_k,x_{k+1},\xbar_k) \in \hood(\theta_k) \times \hood(\theta_k) \times \hood(\theta_k)$, as desired.  Consequently, using an argument similar to the proof of Lemma~\ref{lem.Lipschitz}, $\nabla \barrier (\cdot,\mu_k)$ is Lipschitz over $\hood(\theta_k)$ with constant $\ell_{\nabla f,\Bcal} + 2\mu_k \theta_k^{-2}$, which with the choice of $\alpha_k$, the fact that $\barrier(\cdot,\mu_k)$ is $\psi$-strongly convex for all $k \in \N{}$, and \cite[Eq.~(3.14)]{bubeck2015convex} yields
  \begin{align*}
    0
    &\leq \barrier(x_{k+1},\mu_k) - \barrier(\xbar_k,\mu_k) \\
    &\leq \tfrac{1}{\alpha_k} (x_k - x_{k+1})^T(x_k - \xbar_k) - \tfrac{1}{2 \alpha_k} \|x_k - x_{k+1}\|_2^2 - \tfrac{\psi}{2} \|x_k - \xbar_k\|_2^2.
  \end{align*}
  Therefore, it follows that
  \begin{align*}
    \|x_{k+1} - \xbar_k\|_2^2
      &= \|x_{k+1} - x_k + x_k - \xbar_k\|_2^2 \\
      &= \|x_{k+1} - x_k\|_2^2 + 2(x_{k+1} - x_k)^T(x_k - \xbar_k) + \|x_k - \xbar_k\|_2^2 \\
      &\leq \|x_{k+1} - x_k\|_2^2 - \|x_k - x_{k+1}\|_2^2 - \alpha_k \psi \|x_k - \xbar_k\|_2^2 + \|x_k - \xbar_k\|_2^2 \\
      &= (1 - \alpha_k \psi) \|x_k - \xbar_k\|_2^2,
  \end{align*}
  which is the final desired conclusion.
\eproof

Let us now use the prior two results to show a relationship between consecutive distances from an iterate to the solution of \eqref{prob.opt} that holds for all $k \in \N{}$.

\bproposition\label{prop.final}
  Suppose that Assumptions~\ref{ass.f} and \ref{ass.bounded} hold, the objective~$f$ is $\psi$-strongly convex, and the algorithm employs the updates in \eqref{eq.convex_update}--\eqref{eq.convex_step_size}.  Then, for all sufficiently large $k \in \N{}$ and $\zeta \in \R{n}$ defined as in Proposition~\ref{prop.distance_to_solution}, one has
  \bequation\label{eq.final}
    \|x_{k+1} - x_*\|_2 \leq \sqrt{1 - \alpha_k \psi} \|x_k - x_*\|_2 + 2 \mu_k \|\zeta\|_2 + \Ocal(\mu_k^2).
  \eequation
\eproposition
\bproof
  Combining Propositions~\ref{prop.distance_to_solution} and \ref{prop.contraction_to_xbar} and the triangle inequality, one has
  \begin{align*}
    \|x_{k+1} - x_*\|_2
      &\leq \|x_{k+1} - \xbar_k\|_2 + \|\xbar_k - x_*\|_2 \\
      &\leq \sqrt{1 - \alpha_k \psi} \|x_k - \xbar_k\|_2 + \|\xbar_k - x_*\|_2 \\
      &\leq \sqrt{1 - \alpha_k \psi} (\|x_k - x_*\|_2 + \|\xbar_k - x_*\|_2) + \|\xbar_k - x_*\|_2 \\
      &\leq \sqrt{1 - \alpha_k \psi} \|x_k - x_*\|_2 + 2 \mu_k \|\zeta\|_2 + \Ocal(\mu_k^2),
  \end{align*}
  as desired.
\eproof

At first glance, the result of Proposition~\ref{prop.final} might appear to be useful since $\{\mu_k\} \searrow 0$ also implies that $\{\alpha_k\} \searrow 0$.  Unfortunately, however, the last two terms on the right-hand side of \eqref{eq.final} obstruct an ability to prove that $\{x_k\} \to x_*$, even in this strongly convex and deterministic setting.  To see this, note that the recurrence defined by \eqref{eq.final} has the form of sequences $\{u_k\}$, $\{v_k\}$, and $\{e_k\}$ such that
\bequationNN
  u_{k+1} \leq v_k u_k + e_k,\ \ u_k \in \R{}_{\geq0},\ \ v_k \in [0,1),\ \ \text{and}\ \ e_k \in \R{}_{>0}\ \ \text{for all}\ \ k \in \N{}.
\eequationNN
Such a recurrence yields $\{u_k\} \to 0$ if $\sum_{k=1}^\infty (1 - v_k) = \infty$ and $\lim_{k\to\infty} e_k/(1 - v_k) = 0$; see, e.g., \cite[pg.~45]{Poly87}.  However, with $v_k := \sqrt{1 - \alpha_k \psi}$ and $e_k := C \mu_k$ for some $C \in \R{}_{>0}$ (for simplicity), and for example supposing that $\theta_k = c\mu_k$ for all $k \in \N{}$ (see \eqref{eq.convex_step_size}), one indeed finds $\sum_{k=1}^\infty (1 - v_k) = \infty$, but
\bequationNN
  \tfrac{C\mu_k}{1 - \sqrt{1 - \alpha_k \psi}} = \tfrac{C\mu_k}{1 - \sqrt{1 - \tfrac{\psi}{\ell_{\nabla f,\Bcal} + 2 c^{-2} \mu_k^{-1}}}} \xrightarrow{k\to\infty} \tfrac{C4c^{-2}}{\psi} > 0.
\eequationNN
Consequently, Proposition~\ref{prop.final} does not readily lead to a convergence guarantee for the simplified algorithm stated in \eqref{eq.convex_update}--\eqref{eq.convex_step_size}.  One might modify the algorithm and/or analysis in this section to reach such a guarantee in the deterministic setting upon which one might build a convergence theory for a stochastic algorithm, but we contend that the ultimate conclusions would only be comparable to those for the algorithm analyzed in Section~\ref{sec.analysis}, perhaps with the extensions mentioned in Remark~\ref{rem.extend}.

\section{Numerical Results}\label{sec.numerical}

Our numerical experiments serve two main purposes: (1) we demonstrate that our interior-point method, which we refer to as SIPM, is reliable over a well known set of test problems, and (2) we compare the performance of SIPM with a \edit{projected-stochastic-gradient} method, which we refer to as PSGM.  \edit{(We also ran experiments with the algorithm discussed in Section~\ref{sec.convex}, namely, the one specified in \eqref{eq.convex_update}.  The relative performance of SIPM in comparison to that algorithm was qualitatively identical to the relative performance of SIPM versus PSGM.  Therefore, for the sake of brevity, in this section we only present results for SIPM and PSGM.)}  We implemented a set of test problems and the algorithms in Matlab.  The experiments were conducted on the High Performance Computing cluster at Lehigh University with Matlab R2021b using the Deep Learning Toolbox.

\subsection{Test problems}

We tested the algorithms by training prediction models for binary classification using data from LIBSVM \cite{CC01a}.  From LIBSVM, we selected the 43 binary classification datasets with training data file size at most 8 GB; for these, the numbers of features are in the range $[2, 47263]$ and the numbers of data points are in the range $[44, 5000000]$. We provide results pertaining to training data, and for those datasets with corresponding testing data, we provide results pertaining to that data as well.  Each dataset from LIBSVM consists of $A \in \R{m \times n_f}$ and $b \in \{-1,1\}^m$, where $m$ is the number of data points and $n_f$ is the number of features.

To cover both a convex and a nonconvex objective, we consider two models: (1) logistic regression and (2) a neural network with one hidden layer and a cross-entropy loss function.  For training a \edit{convex} logistic regression model, the number of optimization variables is the number of features plus one for the bias term, i.e., $n = n_f + 1$.  The \edit{nonconvex} neural network model consists of a fully connected hidden layer with $h$ neurons and tanh activation and a fully connected output layer with sigmoid activation.  The number of optimization variables is the number of weights plus bias terms at each node in the hidden and output layers, so $n = (n_f + 2) h + 1$, where $h := \max\{2, \min\{ \lceil \tfrac{n_f}{2} \rceil, 100\}\}$. For both models, we set $l = -1 \times \ones$ and $u = 1 \times \ones$, which causes many bounds to be active at a solution.  Table~\ref{tab: problem parameters} (pg.~\pageref{tab: problem parameters}) shows the number of variables for each problem, i.e., objective and dataset pair.

\subsection{Implementation details}

We generated $x_1$ for each problem with elements drawn from a uniform distribution over $[-0.01, 0.01]$.  This point was fixed for all runs.

SIPM requires the problem-dependent parameters $\kappa_{\nabla f, \Bcal, \infty}$, $\ell_{\nabla f, \Bcal}$, and $\sigma_\infty$.  For consistency across our experiments in both the deterministic and stochastic settings, we employed estimates $\overline{\kappa_{\nabla f, \Bcal, \infty}}$ and $\overline{\ell_{\nabla f, \Bcal}}$, which were set by: (1) temporarily setting these values to~1; (2) running 500 iterations of SIPM using true gradients, these temporary values, and the remaining parameters set as in the next paragraph; and (3) setting, at termination, the values as $\overline{\kappa_{\nabla f, \Bcal, \infty}} \gets \max_{k \in [500]} \left\{\|\nabla f(x_k)\|_\infty\right\}$ and $\overline{\ell_{\nabla f, \Bcal}} \gets \max_{k \in [500] \setminus \{1\}} \left\{\|\nabla f(x_{k-1}) - \nabla f(x_k)\|_2/\|x_{k-1} - x_k\|_2\right\}$.  For the deterministic setting, we set $\overline{\sigma_\infty} \gets 0$, whereas for the stochastic setting we employed an estimate $\overline{\sigma_\infty}$ for each dataset, which was set by: (1) generating 100 stochastic gradients at~$x_1$ with a mini-batch size of $\lceil 0.01 m \rceil$ (which was also the mini-batch size used in all of our experiments for the stochastic setting) and (2) setting $\overline{\sigma_\infty}$ as the maximum $\infty$-norm difference between each of these stochastic gradients and $\nabla f(x_1)$.  Once all of these values were computed---see Table~\ref{tab: problem parameters}---they were fixed for all of our experiments.

Since the performance of an interior-point method is affected by the initial and final values of the barrier parameter, we recommend choosing $\{\mu_k\}$ and $\{\theta_k\}$ based on the computational budget.  Hence, let us define \texttt{maxiter} as an iteration limit for the deterministic setting and $\texttt{maxiter} = (\text{number of epochs})/0.01$ for the stochastic setting, where for the former our experiments consider $\texttt{maxiter} \in \{100, 1000\}$ and for the latter our experiments consider the number of epochs in $\{1, 1000\}$.  (Note that setting \texttt{maxiter} as above for the stochastic setting is consistent with the mini-batch size of $\lceil 0.01 m \rceil$.)  Using this value, and letting $g(x_1)$ denote the gradient (estimate) at the initial point in a run, the parameters for SIPM were set as
\begin{align*}
  \bar{\Delta} &\gets 100,\ \ \Delta \gets \min\{\bar{\Delta}, \min_{i \in [n]}\{u_i - l_i\}\}, \\
  \mu_1 &\gets \max\left\{10^{-5}, \min\left\{\tfrac{10^{-3}\|g(x_1)\|_2 }{\|\diag(u-x_1)^{-1} - \diag(x_1-l)^{-1} \|_2}, 1\right\}\right\},\ \ \text{and} \\
  \theta_0 &\gets \min\left\{ \min_{i \in [n]}\{x_i - l_i\},  \min_{i \in [n]}\{u_i - x_i\}, \bar{\theta}_0\right\},\ \ \text{where}\ \ \bar{\theta}_0\gets \frac{1}{\frac{2}{\Delta} + \frac{\overline{\kappa_{\nabla f, \Bcal,\infty}} + \overline{\sigma_\infty}}{\mu_1}},
\end{align*}
as well as $\mu_k \gets \mu_1 s_k$ and $\theta_k \gets \theta_0 s_k$ for all $k \in [\texttt{maxiter}]$, where $\{s_k\}$ is composed of equal-length repetitions of the elements in \edit{$\{1, 10^{-1}, \dots, 10^{-\nu}, 10^{-8}/\mu_1\}$, where $\nu$ is the largest integer such that $10^{-\nu} > 10^{-8}/\mu_1$.  That is, $$\{s_k\} = \{1,\dots,1,10^{-1},\dots,10^{-1},\dots,10^{-\nu},\dots,10^{-\nu},10^{-8}/\mu_1,\dots,10^{-8}/\mu_1\}.$$
}In this manner, $\mu_1$ ensures that the initial search direction is not dominated by the log-barrier term, whereas $\mu_{\texttt{maxiter}} = 10^{-8}$ ensures that SIPM terminates with a prescribed small barrier parameter.  For the remaining parameters, the implementation used for all $k \in \N{}$: $H_k \gets \overline{\ell_{\nabla f, \Bcal}} I + \mu_k (\diag(x_k - l))^{-2} + \mu_k (\diag(u - x_k))^{-2}$, $\lambda_{k,\min} \geq \overline{\ell_{\nabla f, \Bcal}}$ as the smallest eigenvalue of $H_k$, $\alpha_{k,\buff} \gets (\texttt{maxiter}/k)^{1.1} \geq 1$, $\gamma_{k, \buff} \gets (\texttt{maxiter}/k)^{0.55} \geq 1$, and $(\alpha_k,\gamma_{k,\max})$ as in Assumption~\ref{ass.stochastic}.  By choosing $\alpha_{k,\buff}$ and $\gamma_{k,\buff}$ in this manner, SIPM employs $\alpha_k = \lambda_{k,\min}/\ell_{\nabla f, \Bcal, k}$ and $\gamma_{k, \max} = 1$ for all $k \in [\texttt{maxiter}]$.  Other formulas for $\alpha_{k,\buff}$ and $\gamma_{k,\buff}$ were tested; the values above worked best for our experiments.

For PSGM, in order to have a direct comparison with SIPM, the step sizes were also set using the sequence $\{s_k\}$ defined in the previous paragraph in such a manner that the initial and final step sizes for PSGM were the same as those used by SIPM.  We remark in passing that PSGM has convergence\edit{-to-neighborhood} guarantees when a fixed step size is used, but we did not experiment with such a choice since our aim is to compare with SIPM, which is only defined for diminishing step sizes.

\begin{table}[!htp]\centering
\caption{Problem sizes and algorithmic parameters.}\label{tab: problem parameters}
\scriptsize
\texttt{
\begin{tabular}{lrrrrrrrrrr}\toprule
\multirow{2}{*}{dataset} &\multicolumn{4}{c}{logistic regression} & &\multicolumn{4}{c}{neural network + cross-entropy loss} \\\cmidrule{2-5}\cmidrule{7-10}
&$n$ &$\overline{\ell_{\nabla f,\Bcal}}$ &$\overline{\kappa_{\nabla f, \Bcal, \infty}}$ &$\overline{\sigma_\infty}$ & &$n$ &$\overline{\ell_{\nabla f,\Bcal}}$ &$\overline{\kappa_{\nabla f, \Bcal, \infty}}$ &$\overline{\sigma_\infty}$ \\\midrule
a1a &124 &1.36 &0.14 &0.35 & &7751 &3.12 &0.25 &0.34 \\
a2a &124 &1.38 &0.14 &0.22 & &7751 &3.05 &0.24 &0.21 \\
a3a &124 &1.36 &0.13 &0.33 & &7751 &3.01 &0.20 &0.29 \\
a4a &124 &1.37 &0.14 &0.16 & &7751 &2.95 &0.22 &0.15 \\
a5a &124 &1.37 &0.14 &0.16 & &7751 &2.99 &0.19 &0.15 \\
a6a &124 &1.36 &0.13 &0.11 & &7751 &3.00 &0.21 &0.11 \\
a7a &124 &1.36 &0.14 &0.11 & &7751 &2.97 &0.20 &0.09 \\
a8a &124 &1.36 &0.13 &0.09 & &7751 &2.98 &0.22 &0.08 \\
a9a &124 &1.36 &0.13 &0.07 & &7751 &3.00 &0.23 &0.07 \\
australian &15 &0.22 &0.14 &0.44 & &113 &0.25 &0.07 &0.44 \\
breast-cancer &11 &0.22 &0.15 &0.51 & &61 &0.28 &0.11 &0.51 \\
cod-rna &9 &0.88 &0.02 &0.05 & &41 &0.25 &0.12 &0.05 \\
colon-cancer &2001 &1.66 &0.10 &0.68 & &200201 &2.09 &0.11 &0.65 \\
covtype &55 &0.34 &0.03 &0.02 & &1513 &4.00 &0.45 &0.02 \\
diabetes &9 &0.55 &0.08 &0.40 & &41 &0.25 &0.12 &0.40 \\
duke &7130 &4.48 &0.14 &0.65 & &713101 &21.69 &0.75 &0.53 \\
fourclass &3 &0.37 &0.09 &0.42 & &9 &0.25 &0.11 &0.42 \\
german &25 &1.63 &0.14 &0.50 & &313 &3.90 &0.25 &0.50 \\
gisette\_scale &5001 &6.27 &0.50 &0.26 & &500201 &25.25 &0.50 &0.18 \\
heart &14 &0.15 &0.08 &0.56 & &106 &0.25 &0.06 &0.56 \\
ijcnn1 &23 &0.32 &0.28 &0.03 & &265 &0.25 &0.30 &0.03 \\
ionosphere\_scale &35 &1.22 &0.11 &0.64 & &613 &0.66 &0.10 &0.64 \\
leu &7130 &1.08 &0.06 &0.76 & &713101 &6.69 &0.89 &0.72 \\
liver-disorders &6 &0.58 &0.06 &0.62 & &22 &0.25 &0.09 &0.62 \\
madelon &501 &74.76 &0.29 &0.25 & &50201 &9.08 &0.50 &0.25 \\
mushrooms &113 &1.30 &0.13 &0.17 & &6385 &0.60 &0.07 &0.16 \\
phishing &69 &3.34 &0.48 &0.14 & &2381 &4.52 &0.54 &0.13 \\
rcv1\_train &47237 &0.05 &0.02 &0.12 & &4.72e+6 &0.23 &0.02 &0.12 \\
real-sim &20959 &0.25 &0.14 &0.04 & &2.10e+6 &0.25 &0.14 &0.04 \\
skin\_nonskin &4 &0.44 &0.18 &0.02 & &11 &0.25 &0.22 &0.02 \\
sonar\_scale &61 &2.80 &0.41 &0.54 & &1861 &4.17 &0.41 &0.53 \\
splice &61 &5.57 &0.52 &0.52 & &1861 &6.46 &0.52 &0.52 \\
SUSY &19 &0.34 &0.03 &0.01 & &181 &0.25 &0.03 &0.01 \\
svmguide1 &5 &0.35 &0.11 &0.26 & &13 &0.25 &0.11 &0.26 \\
svmguide3 &23 &0.61 &0.11 &0.38 & &265 &0.25 &0.20 &0.38 \\
w1a &301 &0.59 &0.24 &0.14 & &30201 &0.34 &0.35 &0.13 \\
w2a &301 &0.60 &0.23 &0.12 & &30201 &0.34 &0.35 &0.11 \\
w3a &301 &0.60 &0.24 &0.10 & &30201 &0.34 &0.35 &0.09 \\
w4a &301 &0.60 &0.24 &0.10 & &30201 &0.34 &0.35 &0.05 \\
w5a &301 &0.60 &0.24 &0.08 & &30201 &0.34 &0.35 &0.08 \\
w6a &301 &0.61 &0.23 &0.05 & &30201 &0.34 &0.35 &0.03 \\
w7a &301 &0.61 &0.23 &0.04 & &30201 &0.34 &0.35 &0.03 \\
w8a &301 &0.61 &0.23 &0.03 & &30201 &0.34 &0.35 &0.02 \\
\bottomrule
\end{tabular}
}
\end{table}

\subsection{Comparison of SIPM and PSGM}\label{sec.compare sipm with psg} 

All runs of both algorithms terminated when the iteration limit was reached.  To compare performance, we considered two measures at the final iterate: the objective value $f(x_{\texttt{maxiter}})$ (computed over the training set and testing set, when available) and the norm of a projected gradient $\|\proj_{\hood(0)}(x_{\texttt{maxiter}} - \nabla f(x_{\texttt{maxiter}})) - x_{\texttt{maxiter}}\|_{\infty}$; see, e.g., \cite{birgin2000nonmonotone}.  (Even for the stochastic setting, the projected gradient at the final iterate was computed using the true gradient for the purpose of our comparison.)  In particular, for all runs and each measure, we computed a relative performance measure; e.g., in terms of $f$, we use
\bequationNN \label{eq.relative_performance_measure}
  r_{p}:=\tfrac{f(x_{\texttt{maxiter}}^\text{SIPM}) - f(x_{\texttt{maxiter}}^\text{PSGM})}{\max\left\{f(x_{\texttt{maxiter}}^\text{SIPM}), f(x_{\texttt{maxiter}}^\text{PSGM}), 1\right\}} \in [-1,1] \quad \text{ for $p \in$ set of problems},
\eequationNN
and likewise for the norm of the projected gradient.  The values are within $[-1,1]$ since the final objective values and projected-gradient norms are nonnegative.

\begin{figure}
     \centering
     \includegraphics[width=3.7cm]{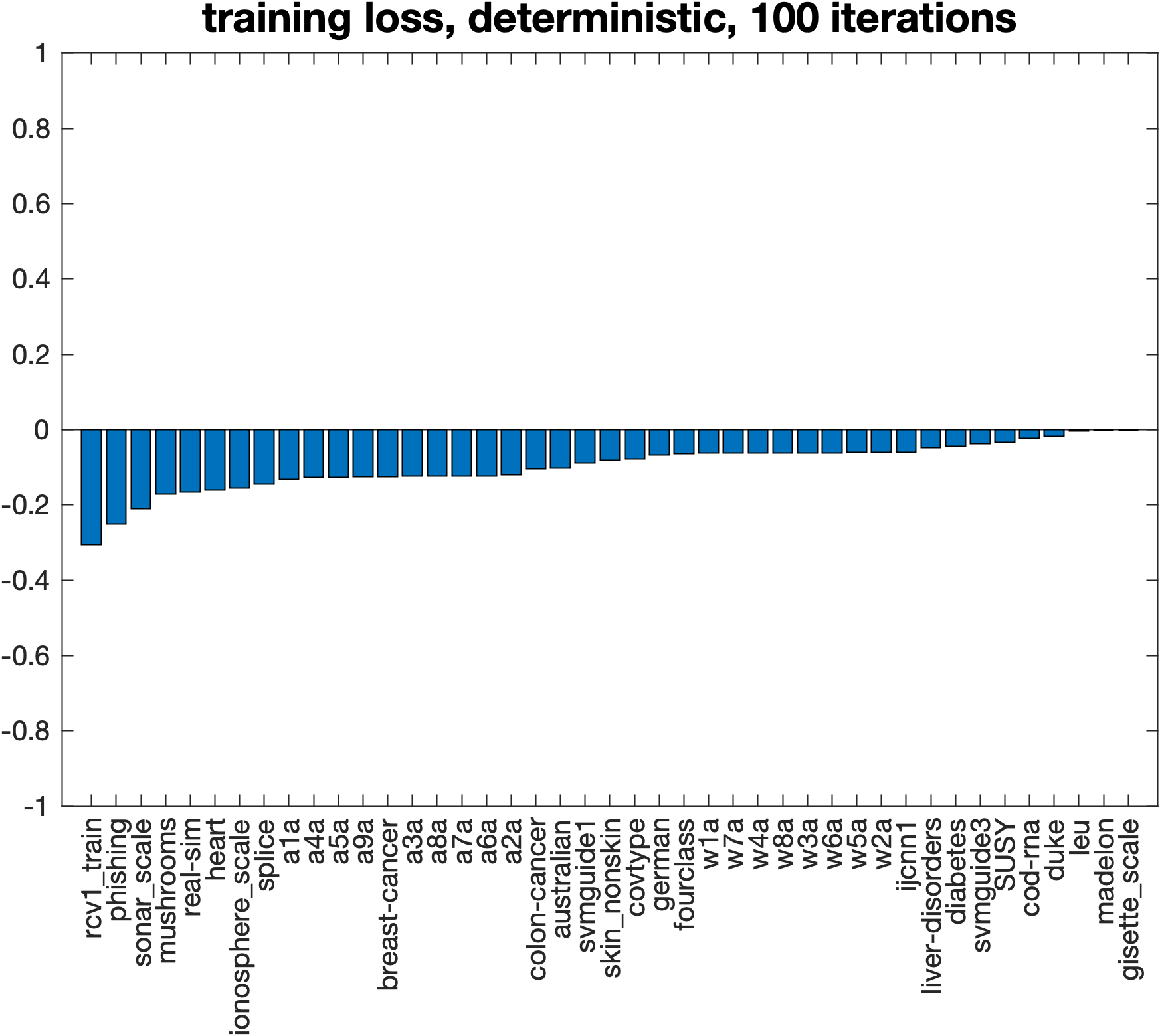}
     \includegraphics[width=3.7cm]{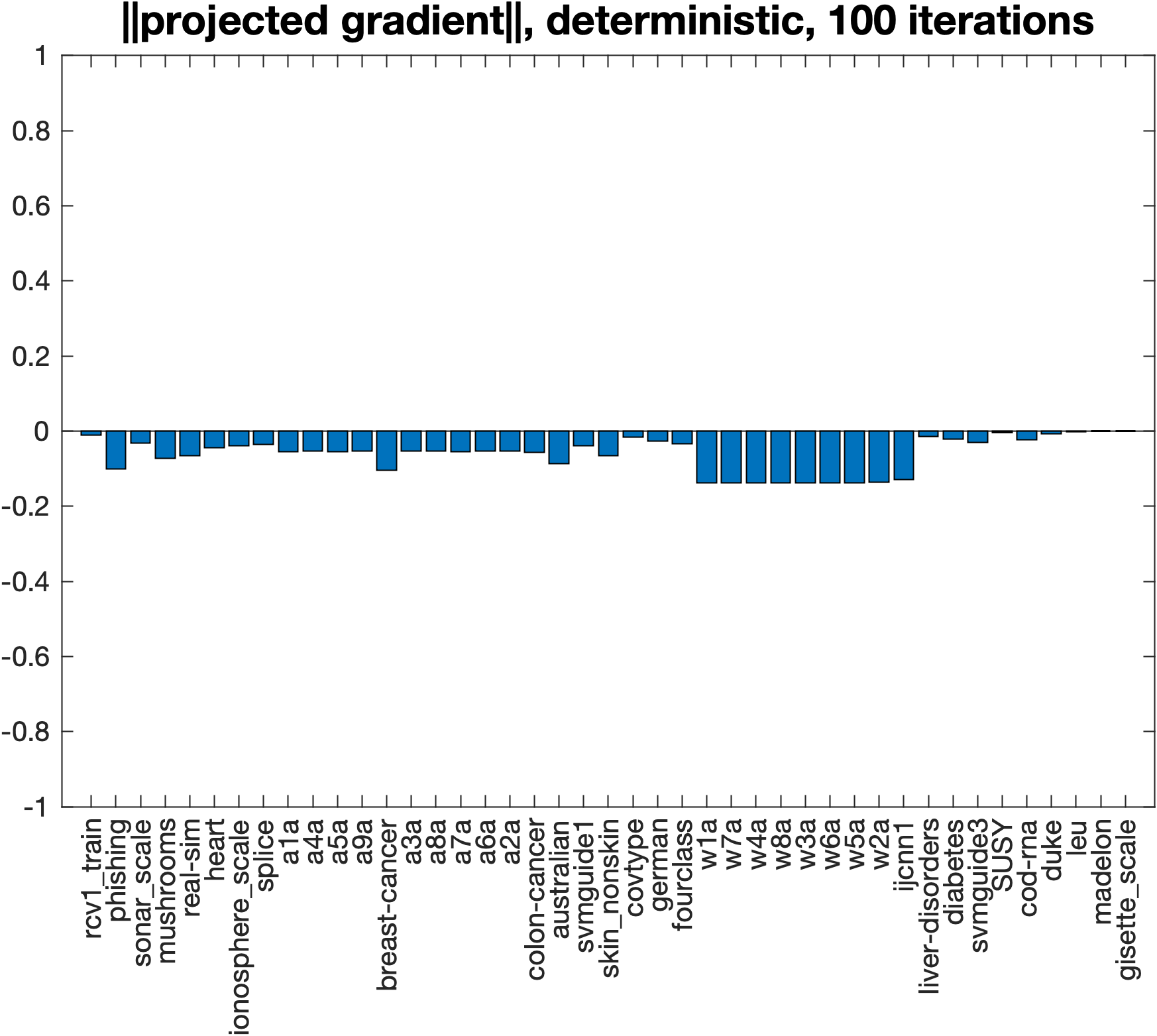}
     \includegraphics[width=3.7cm]{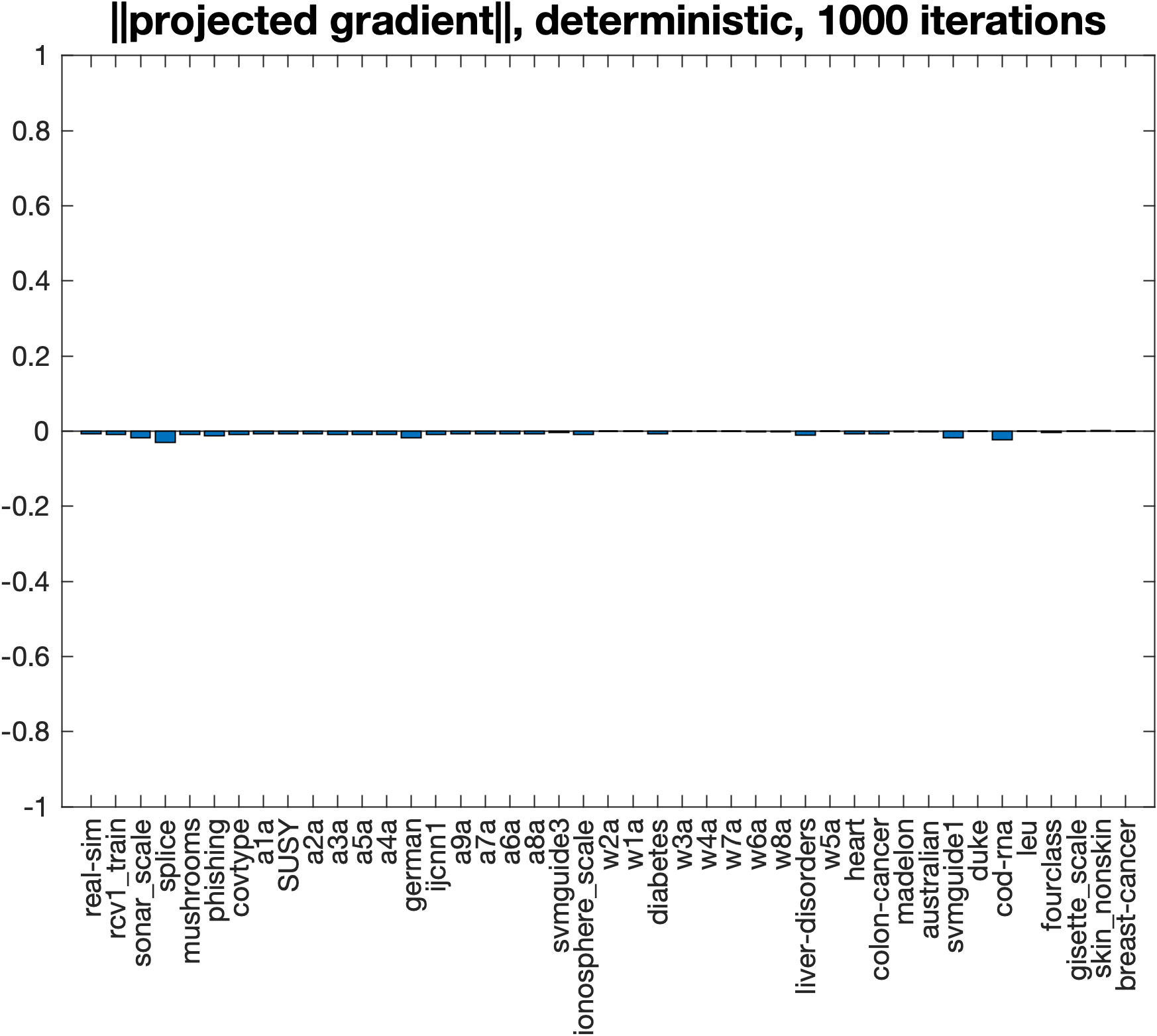}
     \caption{Relative performance of SIPM and PSGM in the deterministic setting when solving logistic regression problems.}
     \label{fig:deterministic relative performance logistic}
 \end{figure}

\begin{figure}
     \centering
     \includegraphics[width=3.7cm]{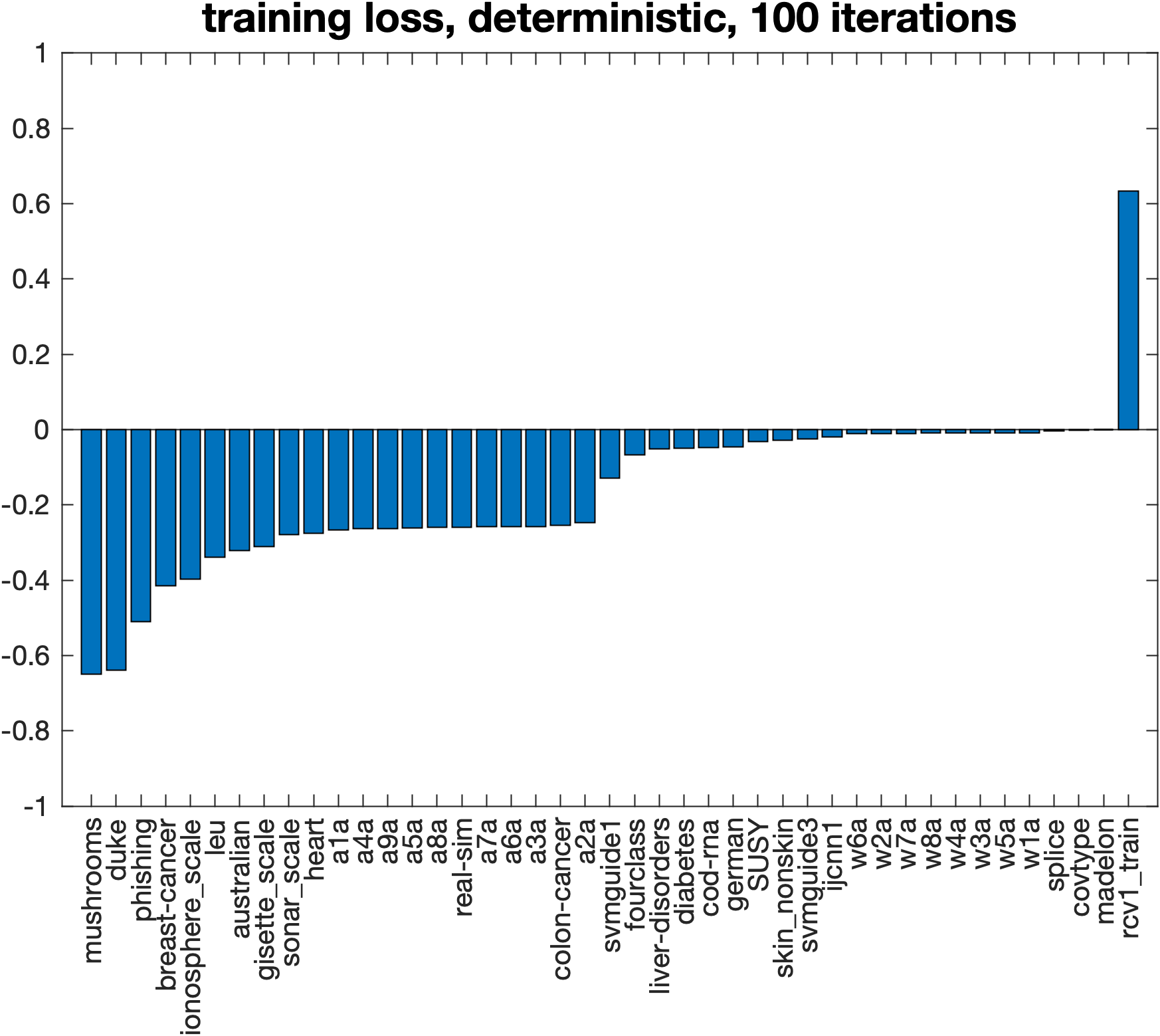}
     \includegraphics[width=3.7cm]{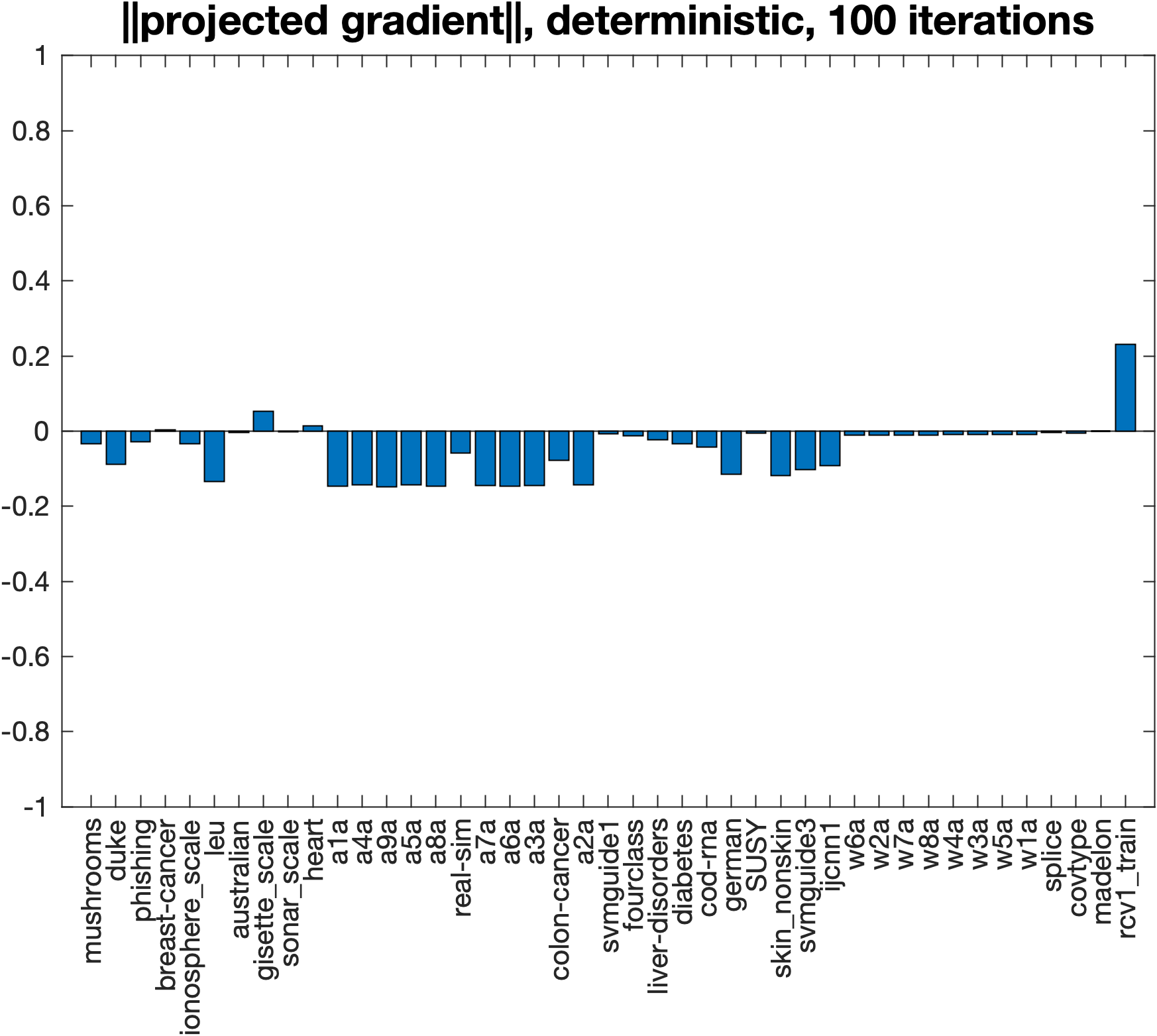}
     \includegraphics[width=3.7cm]{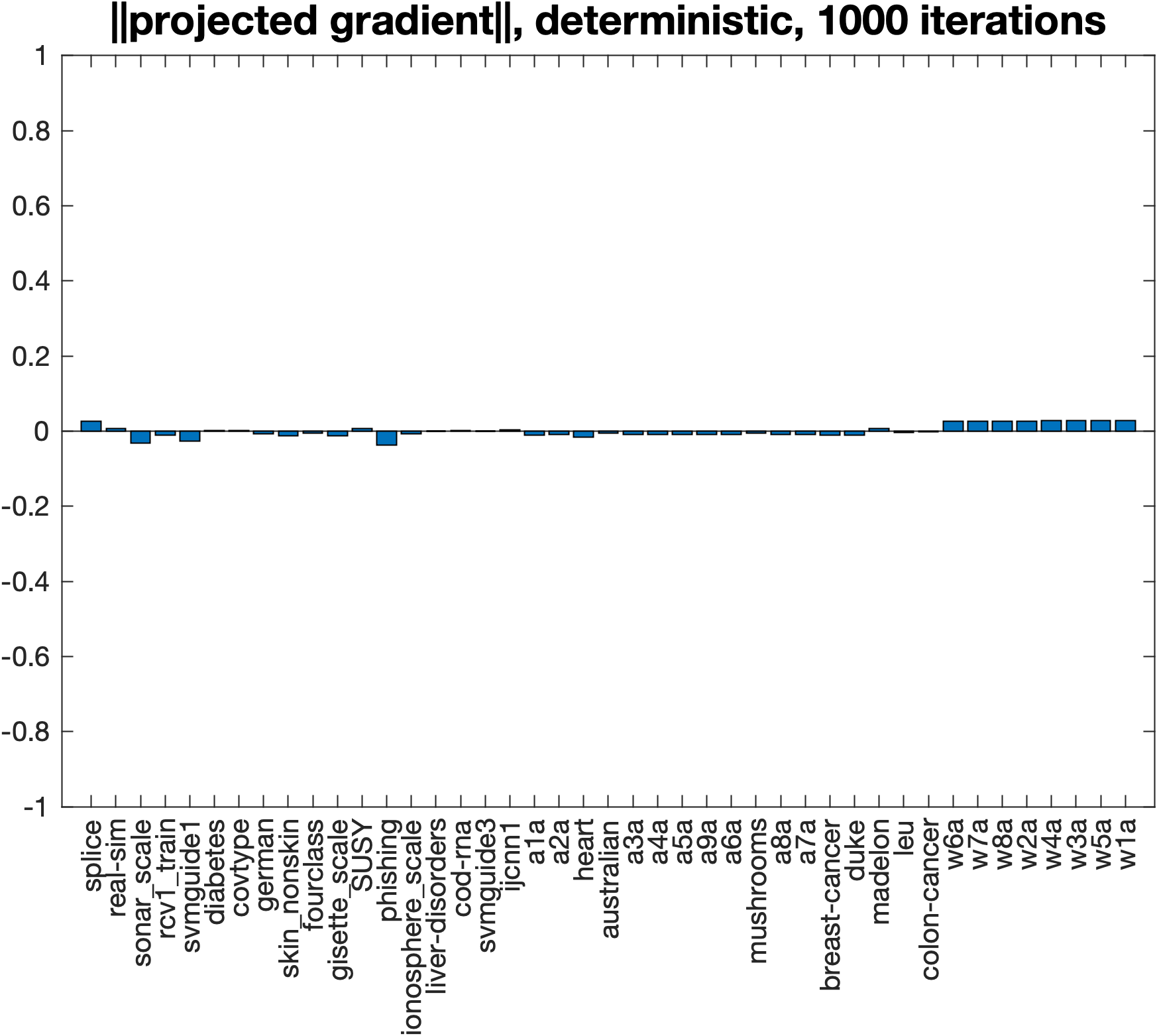}
     \caption{Relative performance of SIPM and PSGM in the deterministic setting when training neural network models $($with one hidden layer$)$ with cross-entropy loss.}
     \label{fig:deterministic relative performance twolayerdnn}
 \end{figure}

Figure~\ref{fig:deterministic relative performance logistic} provides relative performance measures for runs for solving \edit{convex} logistic regression problems.  The bar plot in the first column is for final objective values with respect to the training data when $\texttt{maxiter}=100$ and the plot in the middle column is for projected-gradient norms with respect to the training data when $\texttt{maxiter}=100$.  These show that, within a relatively small iteration limit, SIPM can outperform PSGM.  The bar plot in the third column is for projected-gradient norms with respect to the training data when $\texttt{maxiter}=1000$.  This plot shows that, with a more substantial budget, both algorithms reach points that are nearly stationary, which shows in the deterministic setting that SIPM is as reliable as PSGM. Figure~\ref{fig:deterministic relative performance twolayerdnn} provides similar results for the \edit{nonconvex} neural-network-training problems.  \edit{In these experiments, there was one problem on which PSGM significantly outperformed SIPM when the iteration budget was $\texttt{maxiter}=100$.  That said, for the larger iteration budget of $\texttt{maxiter}=1000$, the performance of both algorithms was again comparable.}

\begin{figure}
     \centering
     \includegraphics[width=3.7cm]{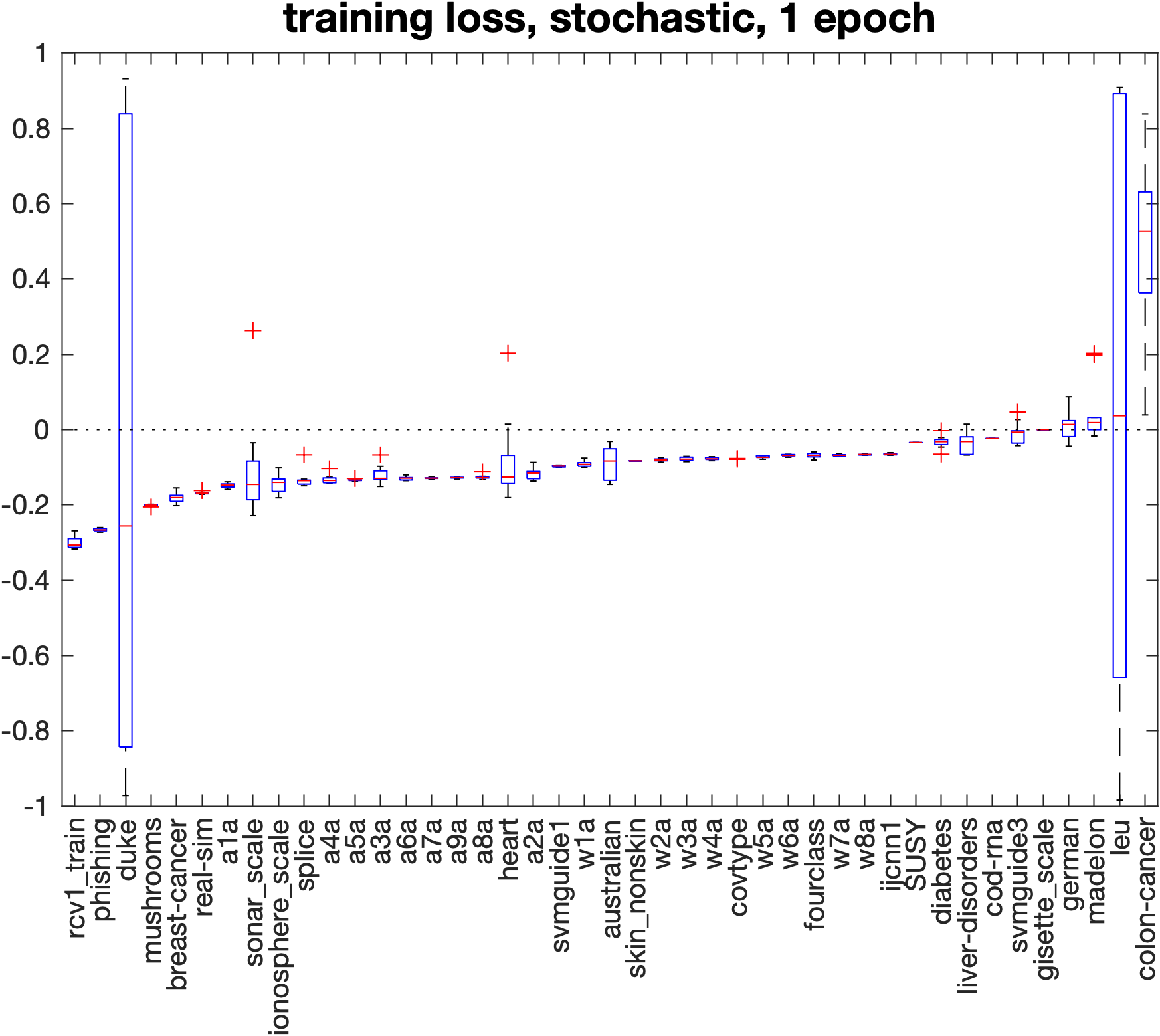}
     \includegraphics[width=3.7cm]{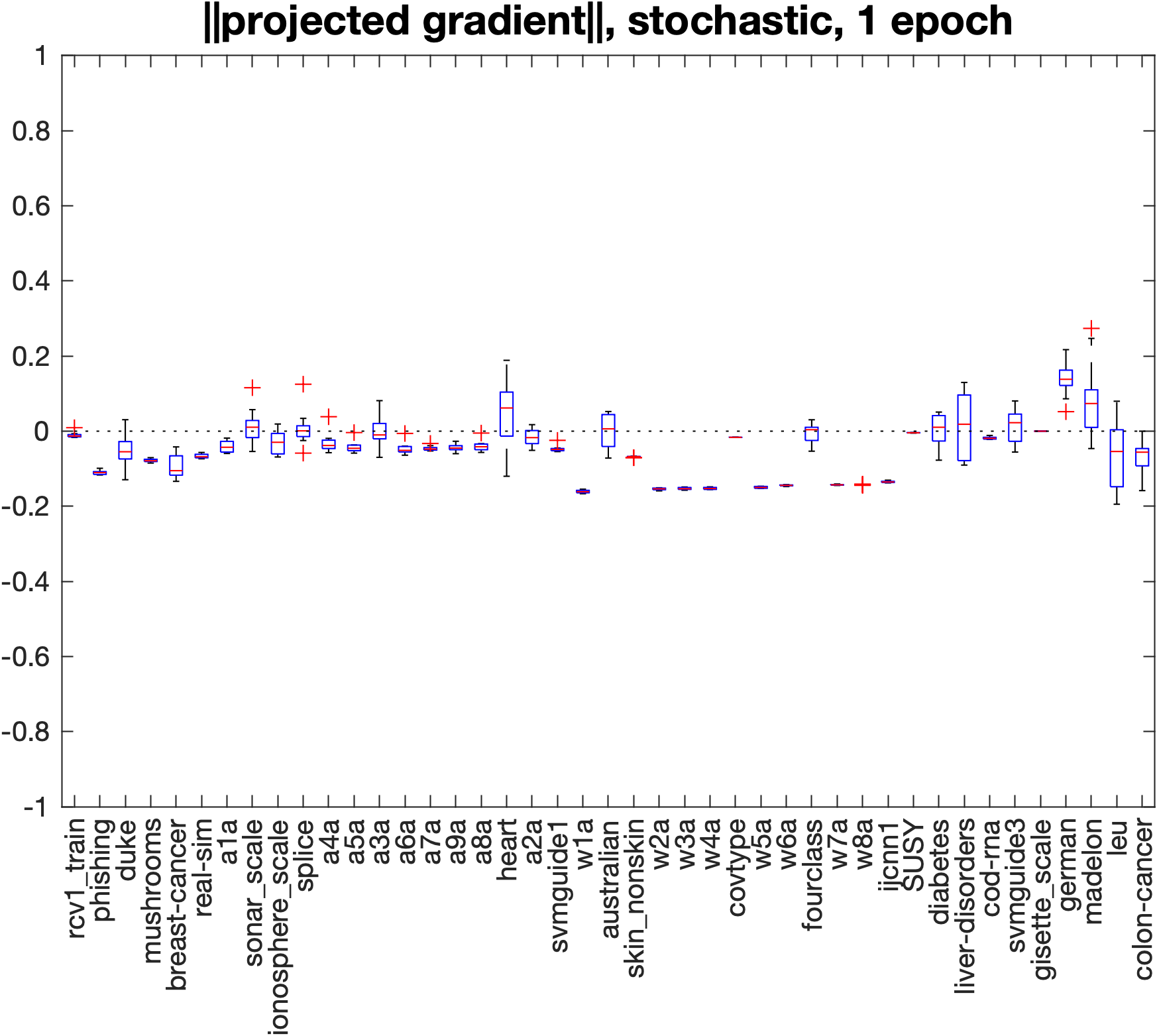}
     \includegraphics[width=3.7cm]{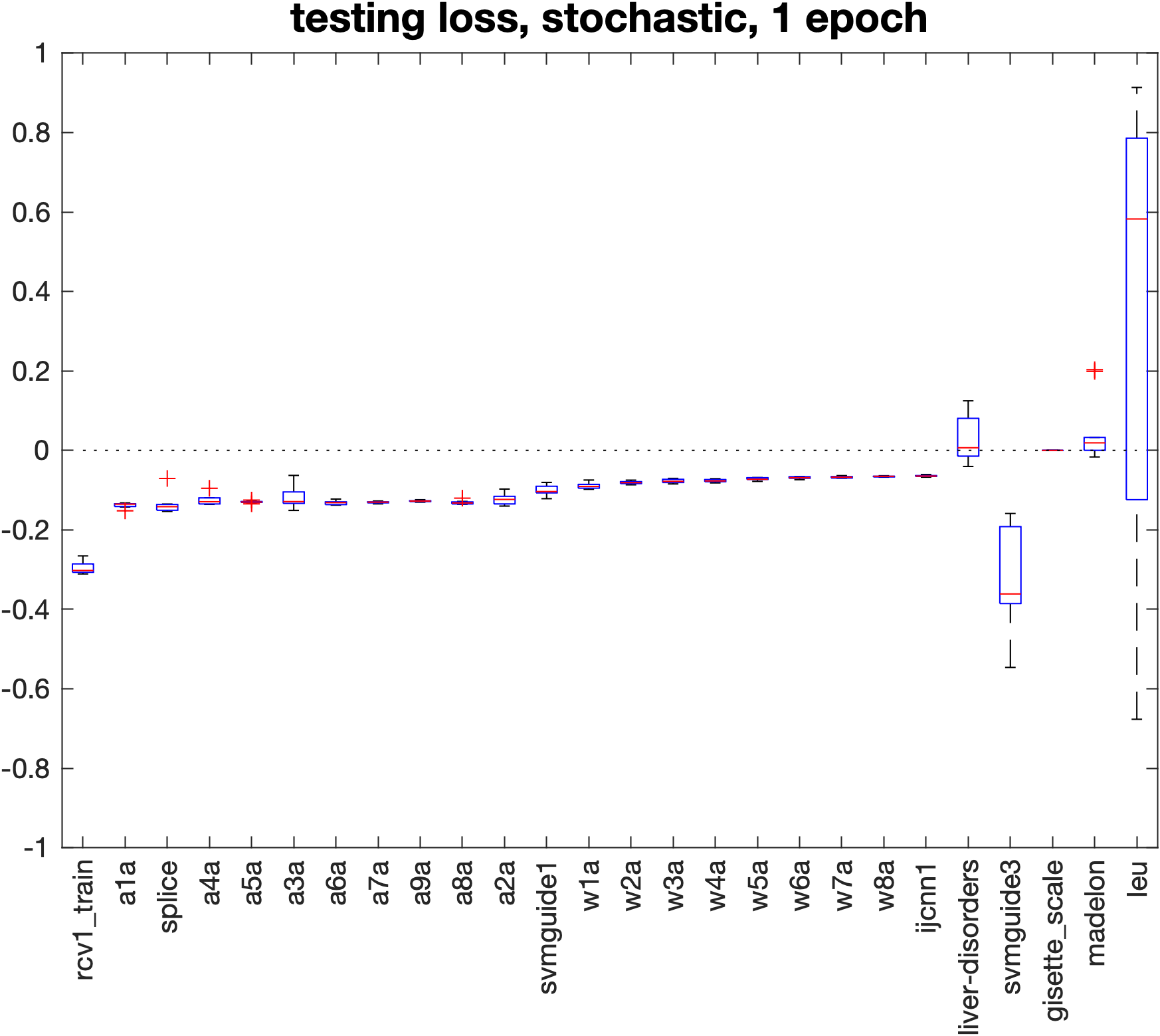}\\   
     \includegraphics[width=3.7cm]{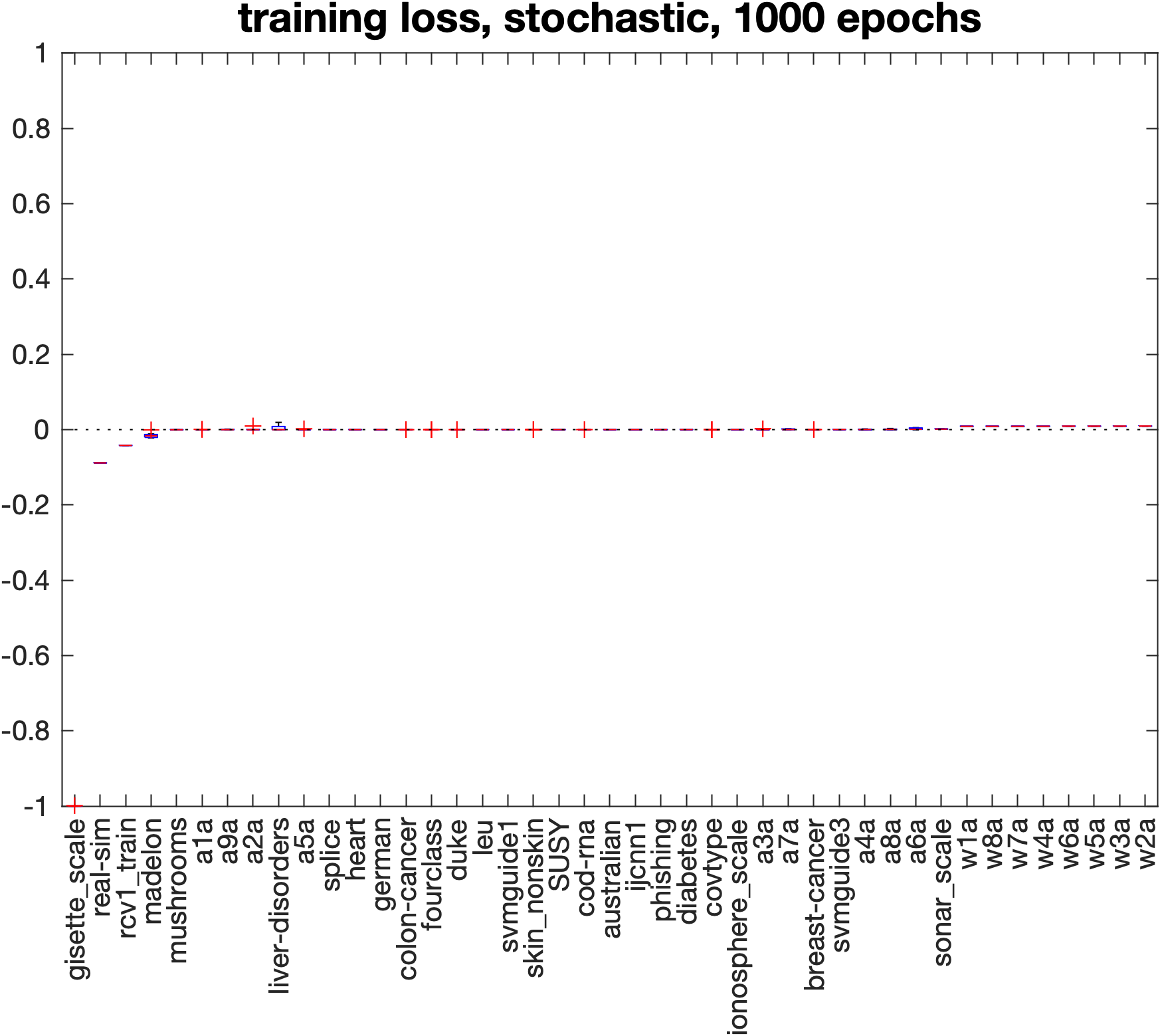}
     \includegraphics[width=3.7cm]{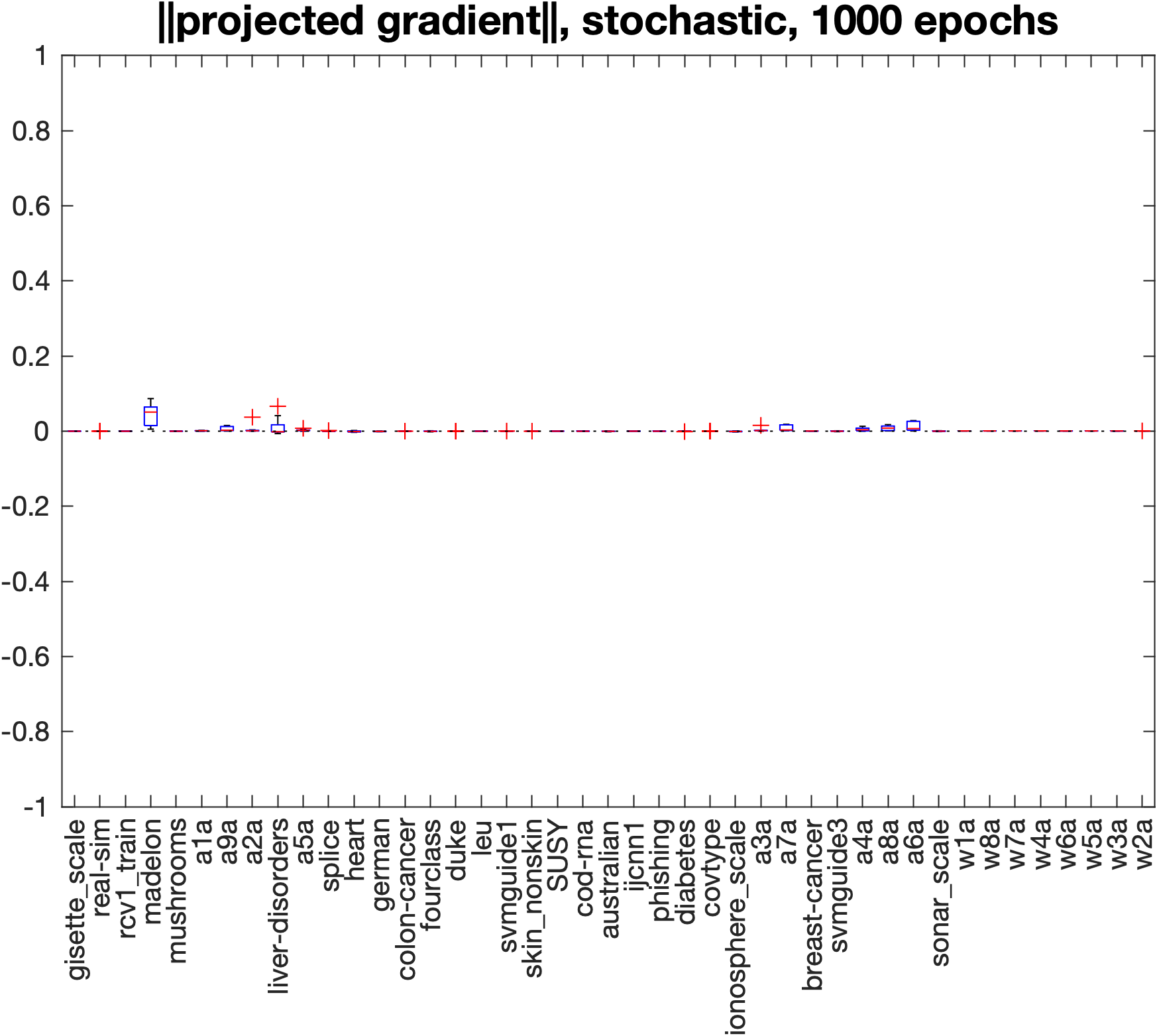}
     \includegraphics[width=3.7cm]{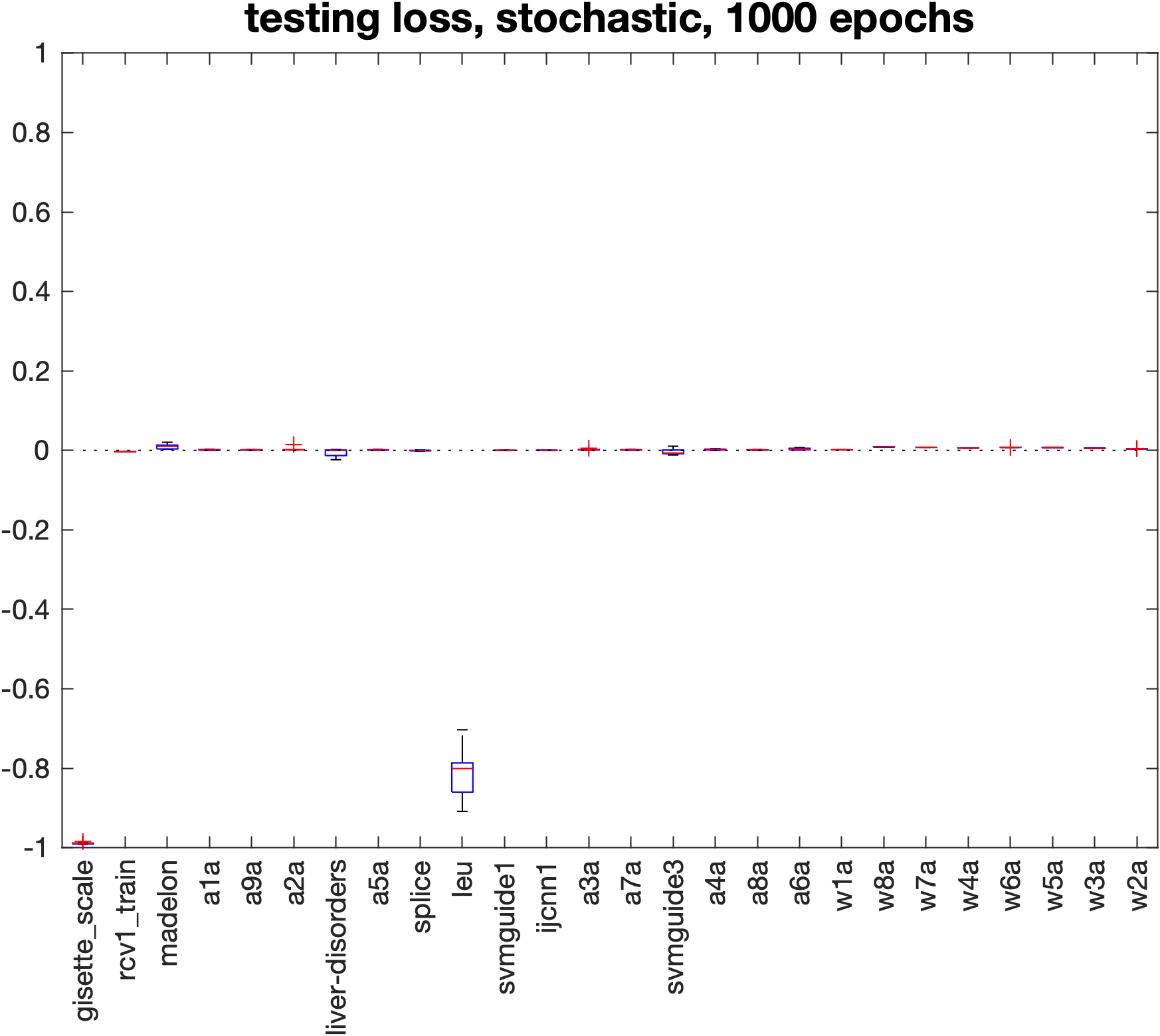}     
     \caption{
     Relative performance of SIPM and PSGM in the stochastic setting $($over 10 runs for each problem$)$ when solving logistic regression problems.  Among the 43 datasets considered for our test problems, there are 26 with corresponding testing datasets $($see last column$)$}
     \label{fig:stochastic relative performance logistic}
 \end{figure}

\begin{figure}
     \centering
     \includegraphics[width=3.7cm]{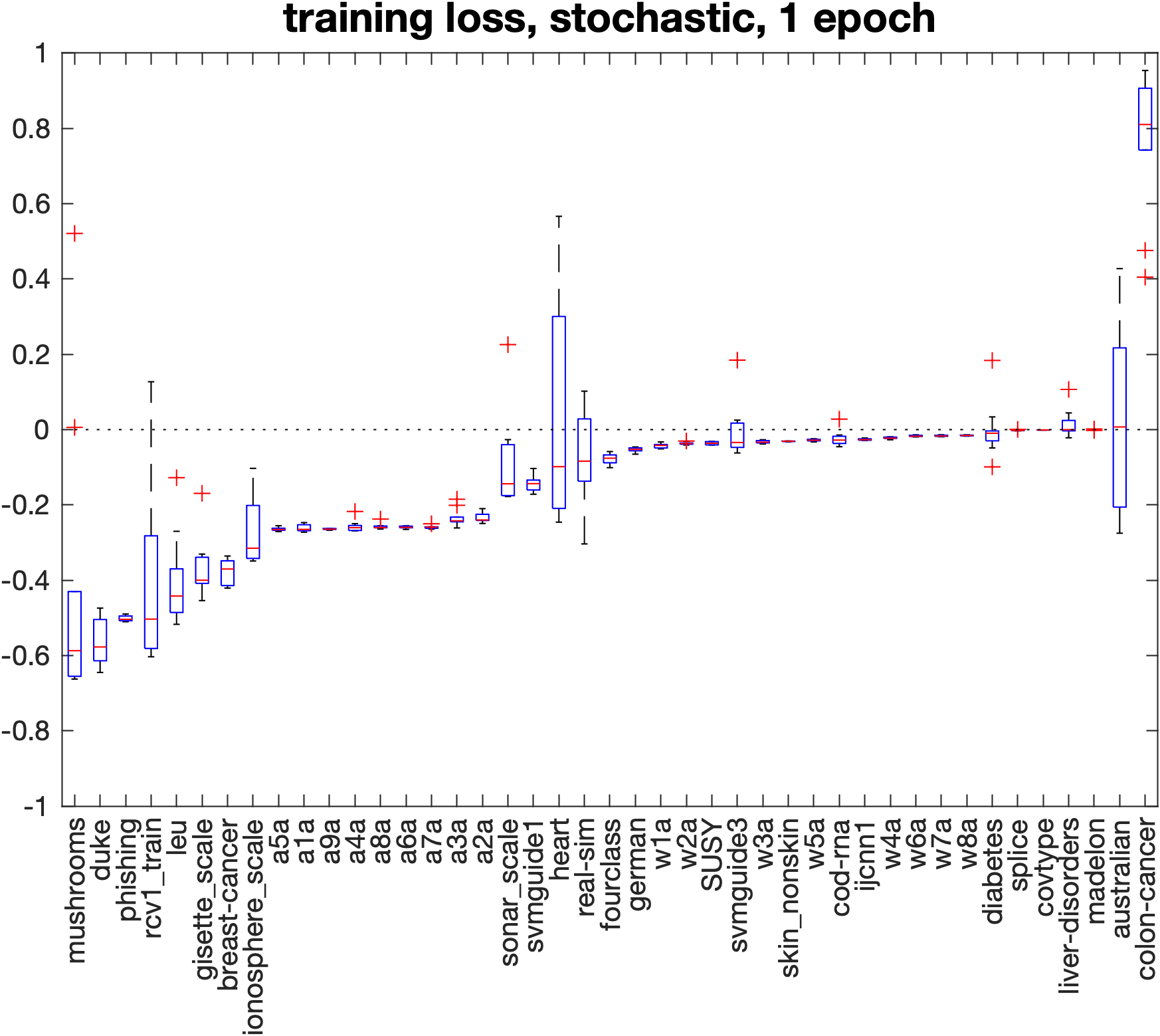}
     \includegraphics[width=3.7cm]{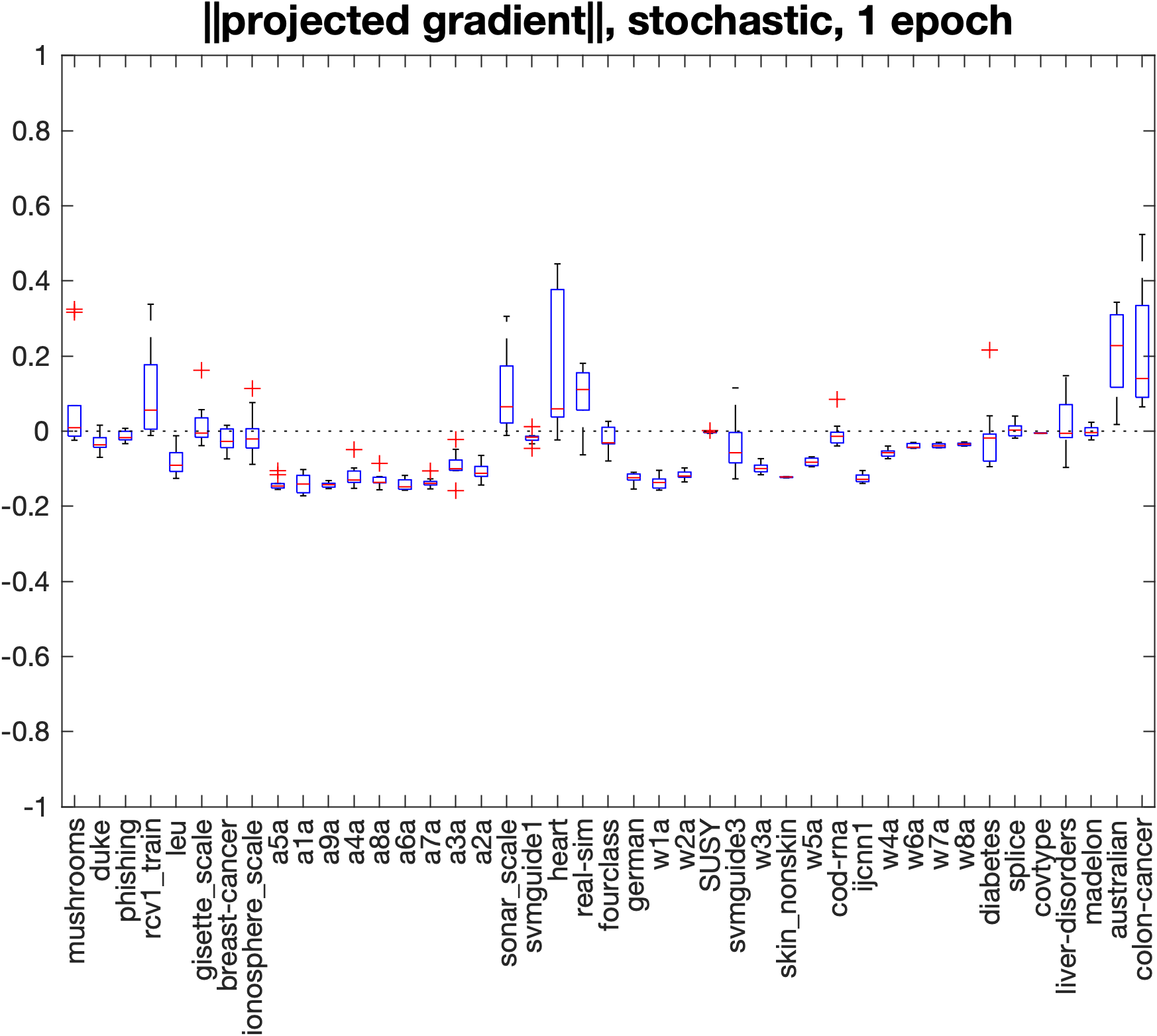}
     \includegraphics[width=3.7cm]{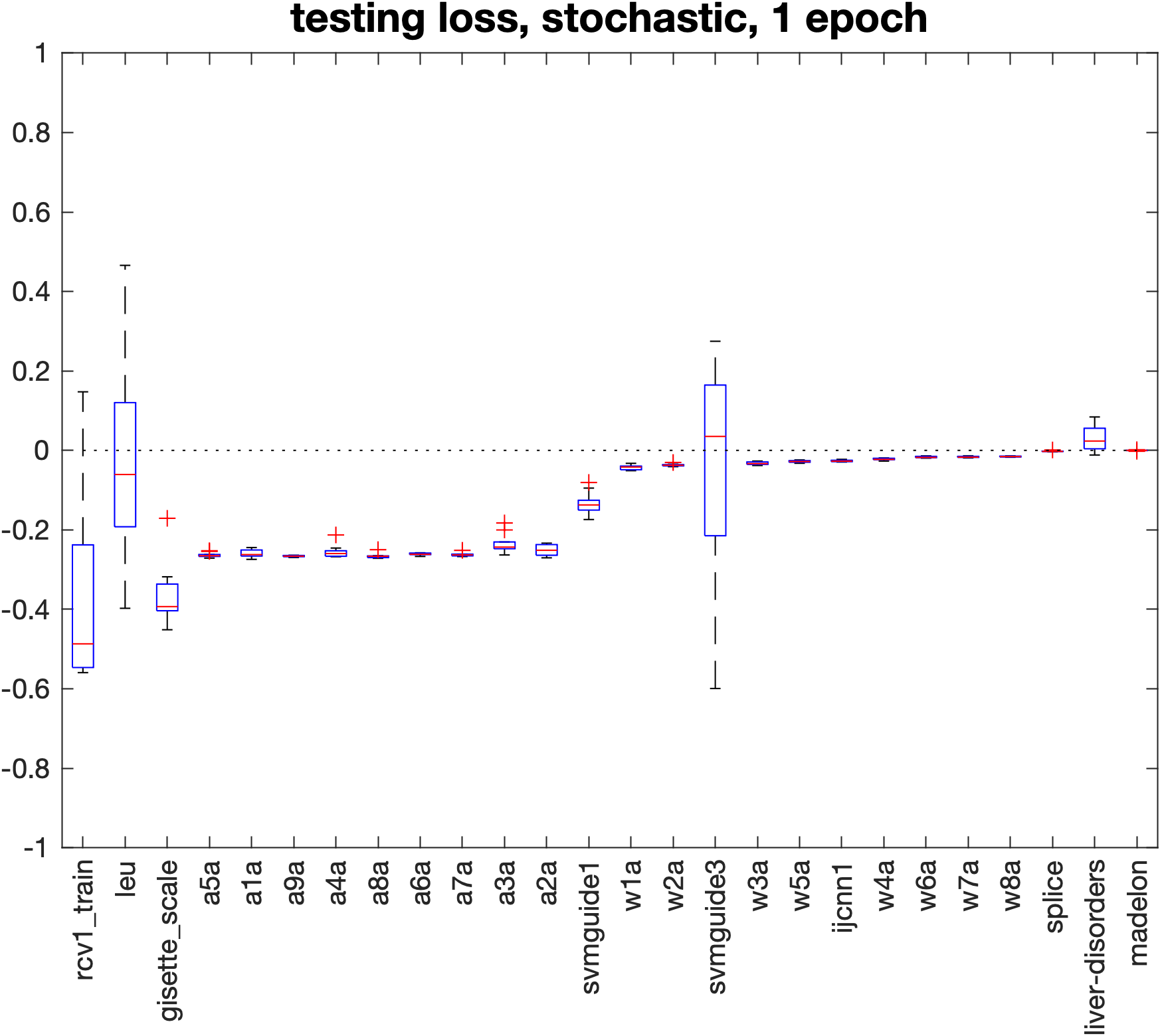} \\ 
     \includegraphics[width=3.7cm]{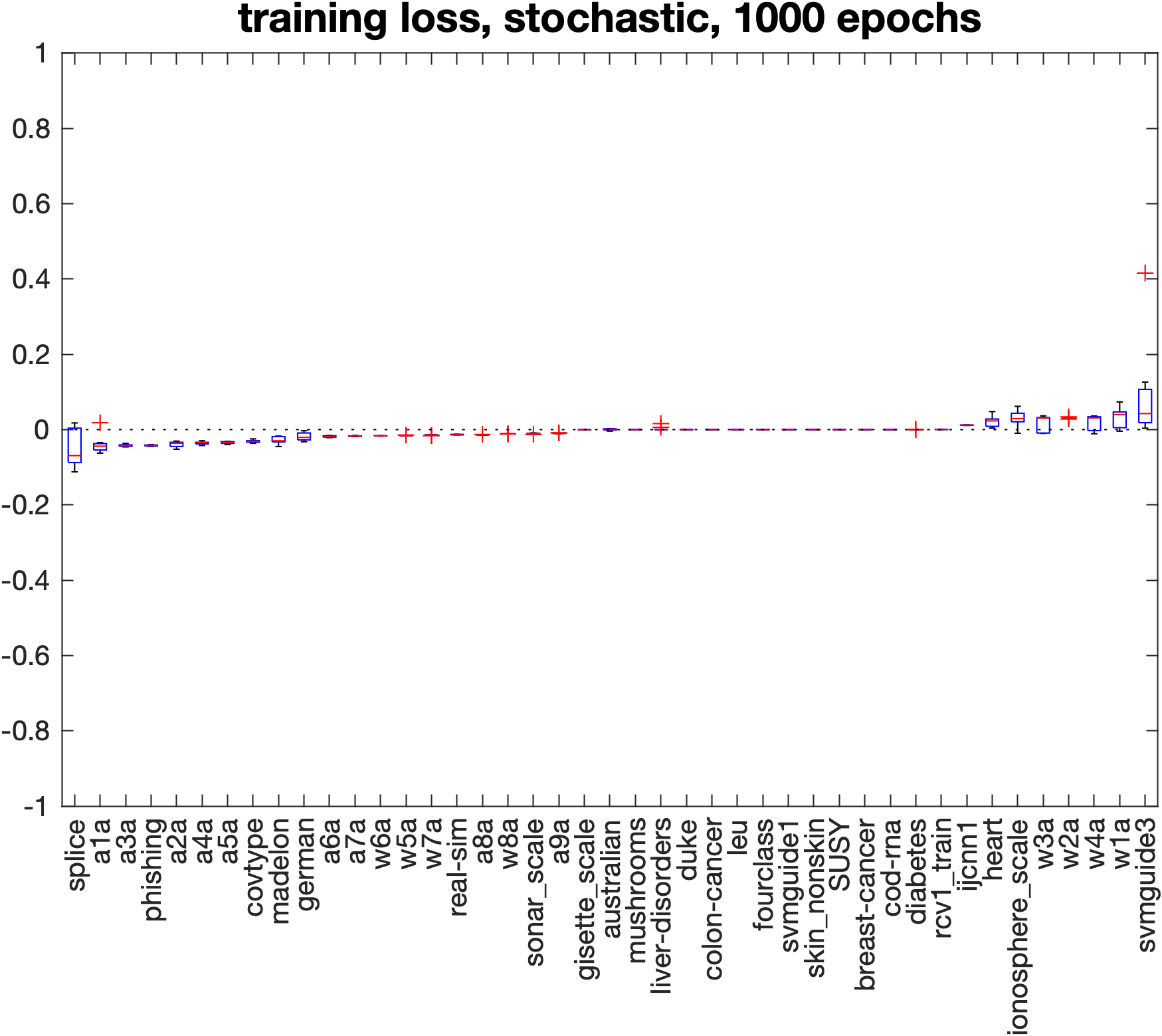}
     \includegraphics[width=3.7cm]{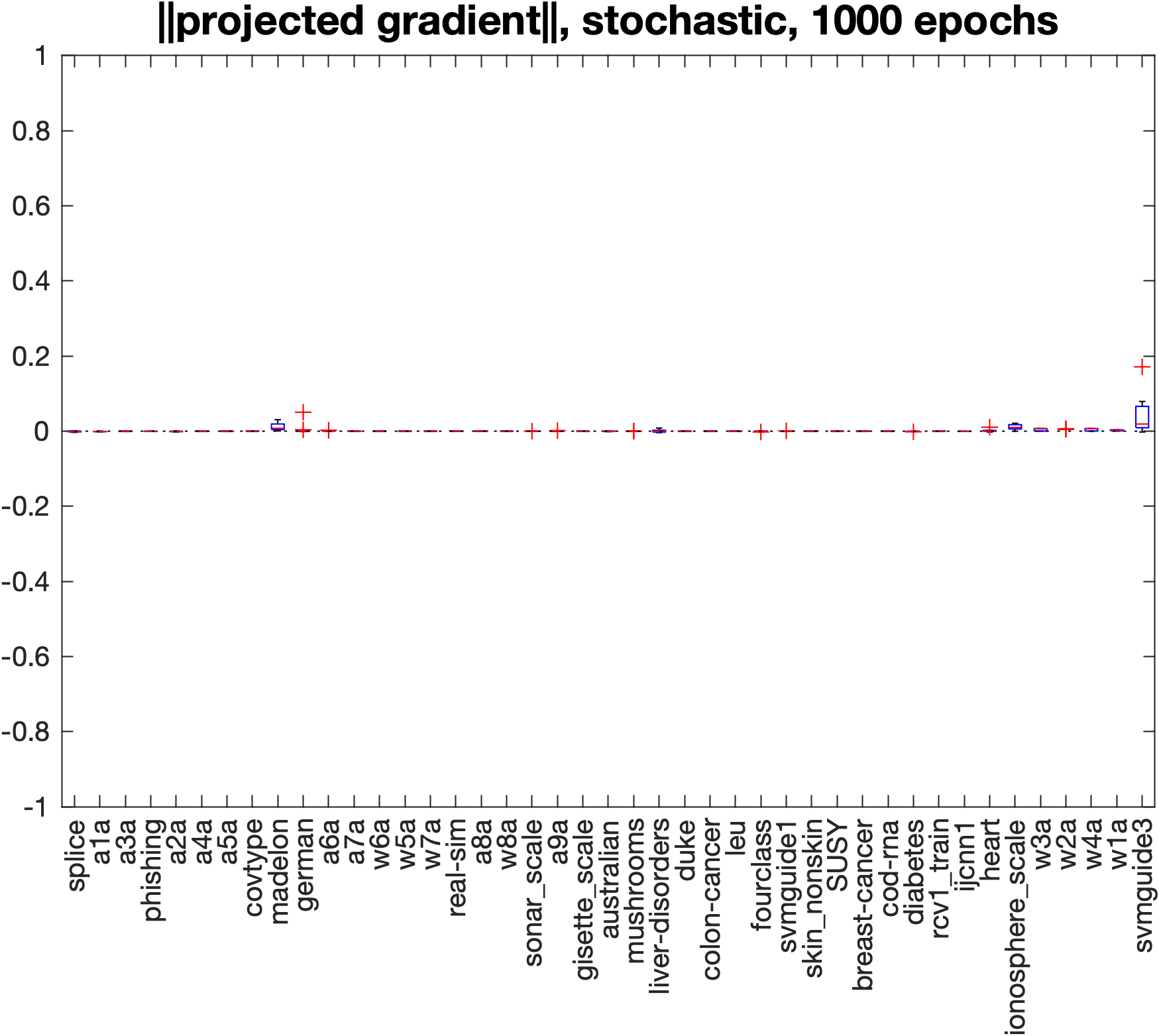}
     \includegraphics[width=3.7cm]{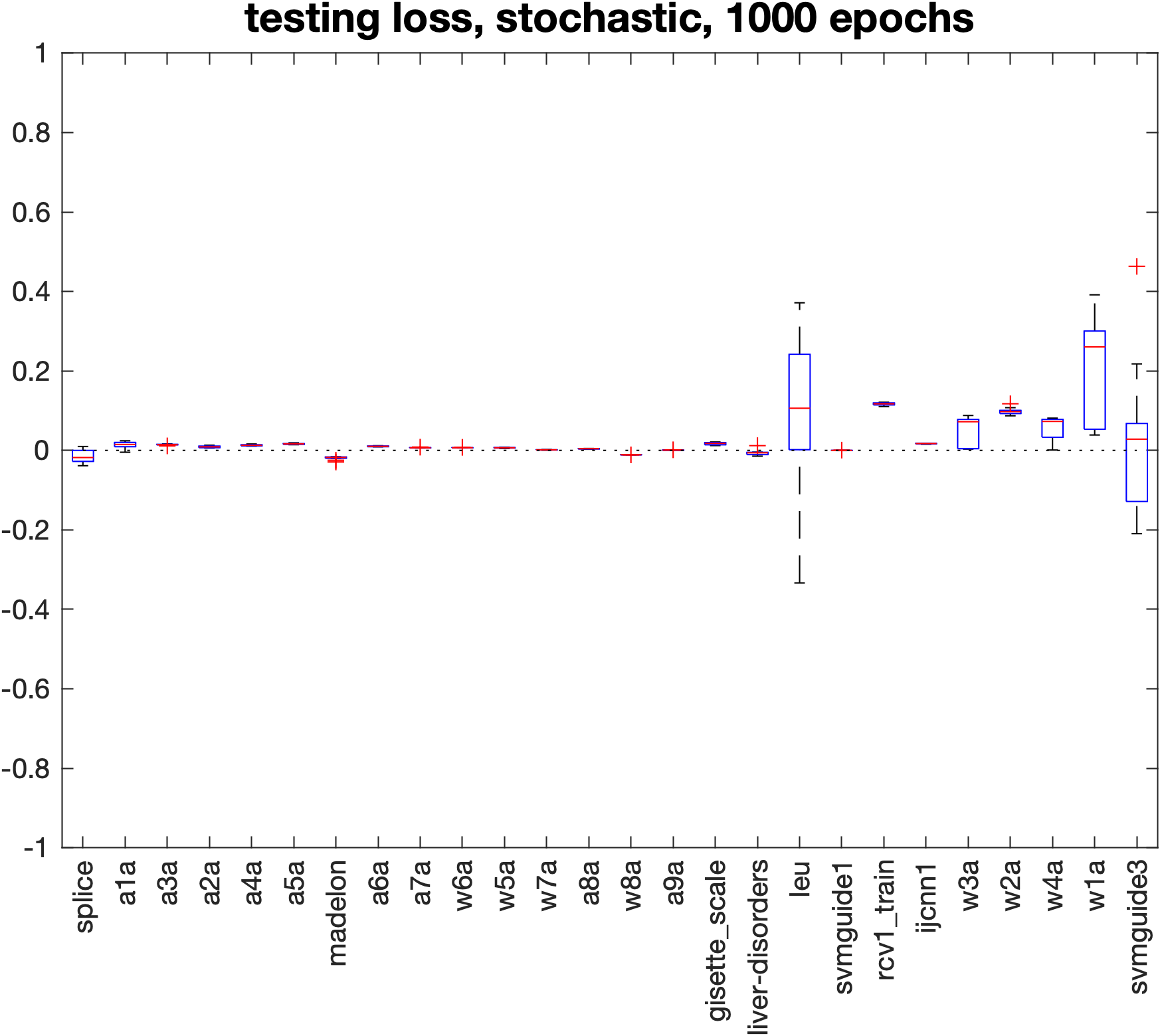}     
     \caption{Relative performance of SIPM and PSGM in the stochastic setting $($over 10 runs for each problem$)$ when training neural network models $($with one hidden layer$)$ with cross-entropy loss.}
     \label{fig:stochastic relative performance twolayerdnn}
 \end{figure}

Figures~\ref{fig:stochastic relative performance logistic} and \ref{fig:stochastic relative performance twolayerdnn} provide results for the stochastic setting in the form of box plots when each algorithm is employed to solve each problem 10 times.  The first rows in each figure consider runs over 1 epoch while the second rows consider runs over 1000 epochs.  The first columns are for training loss, the middle columns are for projected-gradient norms over the training data, and the third columns are for testing loss.  Corresponding to the goals of our experiments, these results show that SIPM performs well compared to PSGM when the budget is relatively small, and is as reliable as PSGM when the budget is large.

\section{Conclusion}\label{sec.conclusion}

We have proposed, analyzed, and provided the results of numerical experiments with a stochastic interior-point method for solving continuous bound-constrained optimization problems.  The algorithm is unique in various aspects (see Section~\ref{sec.contributions}).  In future work, it will be interesting to pair the algorithmic strategies proposed in this paper with stochastic approximation strategies for solving equality-constrained problems toward the complete design of stochastic interior-point methods for solving generally constrained optimization problems.

\section*{Acknowledgments}

This work is supported by the U.S.~National Science Foundation (NSF) award CCF-2139735, U.S.~Office of Naval Research award N00014-21-1-2532, and OP VVV project CZ.02.1.01/0.0/0.0/16\_019/0000765.  Portions were conducted on Lehigh University's Research Computing infrastructure partially supported by the NSF award 2019035.

\bibliographystyle{plain}
\bibliography{references}

\end{document}